\documentclass[12pt]{amsart}
\usepackage{latexsym,amssymb,amsmath} 

\theoremstyle{plain}
 \newtheorem{theorem}{Theorem}[section]
 \newtheorem{proposition}[theorem]{Proposition}
  \newtheorem{lemma}[theorem]{Lemma}
  \newtheorem{corollary}[theorem]{Corollary}
  \newtheorem{conjecture}[theorem]{Conjecture}
 \newtheorem{definition}[theorem]{Definition}
  \newtheorem{example}[theorem]{Example}
  
  \newtheorem{question}[theorem]{Question}
  \newtheorem{remark}[theorem]{Remark}
  \newtheorem{notation}[theorem]{Notation}
  \newtheorem{condition}[theorem]{Condition}

\numberwithin{equation}{section}

\def\HH{{\mathbb H}}

\def\RR{{\mathbb R}}

\begin{document}

\title
{Regular cell complexes in total positivity}
\author{Patricia Hersh}
\email{plhersh@ncsu.edu}


\begin{abstract}
Fomin and Shapiro conjectured that the link of  the 
identity in the Bruhat stratification of the totally nonnegative real part of the 
unipotent radical of a Borel subgroup in a semisimple, simply connected algebraic group
defined and split over $\RR $  is a regular CW complex homeomorphic to  a ball.
The main result of this paper is a proof of this conjecture.  This completes the 
solution of the question of Bernstein of identifying regular CW complexes arising 
naturally 
from representation theory having the  (lower) 
intervals of Bruhat order as their 
closure posets.  
A key ingredient  is a new criterion for determining
whether a finite CW complex is regular with respect to a choice of
characteristic maps;
it most naturally applies   
to images of maps from regular CW complexes  
and is based on an interplay of combinatorics of the closure poset 
with codimension one 
topology.  
\end{abstract}

\maketitle


\section{Introduction}
\label{definition-section}

In this paper, the
following conjecture of Sergey Fomin and  Michael Shapiro from  \cite{FS} is proven.

\begin{conjecture}\label{main-conjecture}
Let $Y$ be the link of the
identity in the totally nonnegative real part of the unipotent radical of a Borel subgroup $B$ in a 
semisimple, simply connected algebraic group defined and split over $\RR $. 
Let $B_u = B^- u B^-$ for $u$ in the Weyl group $W$.  
Then the stratification of $Y$ into Bruhat cells 
$Y\cap B_u$ 
is a regular CW decomposition. 
Moreover,  for each $w\in W$, 
$Y_w = \cup_{u\le  w} Y\cap B_u$  is a regular CW 
complex homeomorphic to a ball, as is the link of each of its cells.
\end{conjecture}

This is done in  
Theorem ~\ref{fs-theorem}.    The result
includes, for instance, the 
special case 
that the link of the identity in the Bruhat stratification of the space of upper
triangular matrices with 1's on the diagonal  whose minors are all nonnegative
is a regular CW complex homeomorphic to a ball;  more specifically, this is the 
collection of upper triangular, totally nonnegative matrices 
with 1's on the diagonal and entries immediately above the diagonal summing to a positive
constant, stratified according to which minors are strictly positive and which are 0.  
The poset (partially ordered set) of closure relations is Bruhat order.

This positively 
answers the question below regarding synthetic Schubert varieties
which appeared in a paper of Bj\"orner (see \cite{Bj}), but was actually posed by
Joseph Bernstein (personal communication, Anders Bj\"orner):

\begin{question}\label{bernstein-question}
It would be of considerable interest to know which (CW) posets can be reasonably
interpreted as face posets of cellular decompositions of complex algebraic varieties, 
and whether there is a synthetic construction for doing so.  In particular, can 
`synthetic Schubert varieties'  be naturally associated with the (lower) Bruhat
intervals of any Coxeter group.
\end{question}

Results of Bj\"orner \cite{Bj} combine with results of Bj\"orner and Wachs
\cite{BW} to imply that each interval of Bruhat order is the closure poset
of a regular CW complex.  This is what led to 
Question ~\ref{bernstein-question}. Fomin and Shapiro proved 
in \cite{FS} that the closure relations for
$Y_w = \cup_{u\le  w} Y\cap B_u$ are exactly those of the lower interval 
$(1,w)$ in Bruhat order, they obtained substantial homological
results regarding this space 
(especially in type A), and they formulated Conjecture ~\ref{main-conjecture}.  Lusztig 
interpreted $Y_w$ in \cite{Lu}
as the image of a continuous map $f_{(i_1,\dots ,i_d)}$ given by a reduced
word $(i_1,\dots ,i_d)$ for $w$, which is the viewpoint we will take as well.  The 
tools we develop in order to prove Conjecture ~\ref{main-conjecture}  
give a new approach to  the general question of how to prove that the image of a map 
from a polytope (or slightly more general regular CW complex)
which restricts to a  homeomorphism on the interior but  not 
necessarily  on the boundary is  a regular CW complex homeomorphic to a ball.

Motivation for studying these totally nonnegative parts  $Y_w$ of  varieties  
comes from a relationship 
observed by Lusztig to his theory of canonical bases.  The change of coordinates 
map resulting from applying a braid move to a reduced word specifying a
canonical basis is a tropicalized version of the corresponding change of coordinates
for a totally
nonnegative variety (cf. \cite{Lu}, \cite{Lu2}, \cite{Lu3}).  Trying to 
understand such changes of coordinates
was also an inspiration  for the theory of
cluster algebras (see \cite{BZ}, \cite{FZ}).  Our topological
approach seems to give a somewhat new perspective on these changes of 
coordinates.  The  explicit collapsing maps
we develop and use later will give quite explicit information about fibers of a map 
whose inverse has  been the subject of considerable study (see e.g. \cite{FZ}
and the preprint: A. Postnikov, Total positivity, Grassmannians, and networks, 
arXiv:math.CO/0609184).  For instance, we deduce connectedness of fibers
indirectly, using our Theorem ~\ref{sufficient-theorem} below.

Our starting point 
was the following new criterion for 
determining whether a  finite CW complex 
is  regular 
with respect to a choice of characteristic maps.  See Section ~\ref{background-section}
for a review of the requisite definitions.
This result below
gave us a route through which to approach Conjecture ~\ref{main-conjecture}.  

Conditions 1 and 2 below imply that the closure poset is graded by cell
dimension, ensuring that the subsequent conditions make sense.   Condition 3 is a 
combinatorial condition which (together with Condition 2) 
enables injectivity of attaching maps to be proven by an
induction on difference in dimensions.  Condition 4 gives the base case for
this induction.  
Notably absent 
is a more general  requirement of injectivity  for the attaching maps.

\begin{theorem}\label{sufficient-theorem}
Let $K$ be 
a finite CW complex  with characteristic maps $f_{\alpha }: B^{ \dim e_{\alpha }} \rightarrow
\overline{e_{\alpha } }$.
Then $K$ is regular with  respect to these characteristic maps  $\{ f_{\alpha } \} $ 
if and only if the following conditions hold:
\begin{enumerate}
\item
For each $\alpha $, $f_{\alpha }(B^{\dim 
e_{\alpha }})$ is a union of open cells 
\item 
For each $f_{\alpha }$, the preimages of the
open cells of dimension $\dim e_{\alpha } - 1$ form a  dense subset of the boundary
of $B^{\dim e_{\alpha }}$. 
\item
The closure poset  of $K$ is thin, i.e., each closed interval $[u,v]$ with 
$rk(v)-rk(u)=2$  has exactly four
elements. 
Additionally, each open interval $(u,v)$ with $rk(v) - rk(u) > 2$
is connected.
\item
For each $\alpha $, the
restriction of  
$f_{\alpha }$  to the preimages of the 
open cells of dimension exactly one less than $e_{\alpha }$ is an injection.
\item
 For each $e_{\sigma }\subseteq \overline{e_{\alpha }}$, 
 $f_{\sigma }$ factors as an
embedding
 $\iota : B^{\dim e_{\sigma } } \rightarrow B^{\dim e_{\alpha }}$ followed by $f_{\alpha }$.
\end{enumerate}
\end{theorem}

Theorem ~\ref{sufficient-theorem} is proven in Section ~\ref{tool-section}.
Examples are also given in Section ~\ref{tool-section}
demonstrating that each of Conditions 2, 3, 4, and 5 is not redundant.
Condition 5 makes Theorem ~\ref{sufficient-theorem}  seem likely to be
applicable  primarily to  images of maps from regular CW complexes, 
which is  indeed how we will  use Theorem ~\ref{sufficient-theorem}.  

Our proof of Conjecture ~\ref{main-conjecture} relies critically upon the fact that
condition 4, codimension one injectivity,   follows easily in our setting from the 
exchange axiom for Coxeter groups (which is 
reviewed in Section ~\ref{background-section}).  The analogous 
Coxeter-theoretic statement is not 
true in higher codimensions,  seemingly demonstrating the 
efficacy of Theorem  ~\ref{sufficient-theorem}.
One reason for interest in proving stratified spaces to be 
regular CW complexes 
is the appealing feature of regular CW complexes that 
their topological structure (homeomorphism type) is determined 
by the combinatorics  of their posets of closure relations.

The proof of Conjecture ~\ref{main-conjecture}
also involves the development of a  
combinatorial topological toolkit for performing a series of collapses on a convex polytope (which
in our case is a simplex)
in a manner that preserves regularity and homeomorphism type at each step.  Each collapse
reduces the number of
cells by eliminating some cells and identifying other
cells with each other.   These collapses that we introduce in Theorem ~\ref{topol-collapse2} (and
its extension in Corollary ~\ref{more-general-topol-collapse2})
are much in the spirit of elementary collapses, but with a 
tighter control on the maps which enables us to preserve not just homotopy type but homeomorphism type as well.  They resemble the process of Bing shrinking (cf. \cite{Bi}) 
in that we extend collapsing maps across 
collars by giving paths of homeomorphisms
deforming each collapsing map  to the identity map.
We specifically needed to develop a class of collapses that would only 
identify  points having the 
same image under $f_{(i_1,\dots ,i_d)}$, while restricting ourselves to operations where we 
could control homeomorphism type and regularity.  To this end, we collapse
cells across families of curves which seem typical  
enough of fibers of maps of interest arising e.g. in  combinatorial representation theory 
to be likely to  be useful for other examples of interest as well.

While these  collapses  
are topological in nature, 
we  have gone to considerable effort to  make the 
criteria one must check in order to use them
as combinatorial as possible.  This is done 
not only to  help us with the proof  of Conjecture ~\ref{main-conjecture},
but also to facilitate 
possible future applications to other stratified spaces of interest  in combinatorics and  
representation theory 
such as the double Bruhat decomposition for the totally nonnegative part of the 
Grassmannian or the totally nonnegative part of the flag variety, as discussed 
briefly in Section ~\ref{flag-section}.  
 Taken together, Theorem ~\ref{sufficient-theorem} 
 and Theorem ~\ref{topol-collapse2} provide a fairly combinatorial 
 general approach  to proving that images of sufficiently nice  
 maps from polytopes are regular CW complexes homeomorphic to balls.  
 
Another crucial ingredient in the proof of Conjecture ~\ref{main-conjecture}  
is the 0-Hecke algebra associated to a  Coxeter group $W$.
The  relations of the 0-Hecke algebra capture in a completely natural way 
which faces of a simplex (indexed by subwords of a reduced word)
map to the same cell
under $f_{(i_1,\dots ,i_d)}$, and in fact provide a dictionary from the topology of 
point identifications
to the combinatorics of cell identifications through suppression of a parameter. 
Checking the  requirements for our collapses
thereby translates to an analysis of properties of 
reduced and nonreduced words in this 0-Hecke algebra.

 While 
$f_{(i_1,\dots ,i_d)}$  itself  
is not a homeomorphism, we use  the aforementioned
collapses to construct a quotient space upon which the induced map 
$\overline{f_{(i_1,\dots ,i_d)}}$ will act homeomorphically, allowing us to understand its
image based upon our understanding of this quotient space.
These collapses eliminate exactly the faces of a simplex 
indexed by the subwords of a reduced word $(i_1,\dots ,i_d)$ that are themselves
non-reduced.   Theorem 
~\ref{sufficient-theorem}  gives a way then  to prove that the induced map
$\overline{f_{(i_1,\dots ,i_d)}}$ on the resulting
quotient space is a homeomorphism.  This is what guarantees 
that no further identification is 
necessary once we have performed all of the identifications which the 
non-reduced subwords  necessitate.

The remainder of the introduction
gives a more thorough overview of the main ideas going into the proof of the Fomin-Shapiro
Conjecture, including the new tools leading up to it, then briefly discusses 
other possible future applications of our approach.   Section ~\ref{background-section}
provides background and terminology in topology, topological combinatorics, Coxeter groups
and their 0-Hecke algebras, and in total positivity theory, respectively; 
readers might find it useful to read one or more of these background sections
even prior to reading the  proof overview.  Sections ~\ref{tool-section} and 
~\ref{topol-collapse-section} establish the key topological 
tools, namely Theorem ~\ref{sufficient-theorem} (our regularity criterion for CW 
complexes) and
Theorem ~\ref{topol-collapse2} (our  method for collapsing cells), respectively.
Section ~\ref{0-hecke-section} develops combinatorial properties of the 0-Hecke 
algebra. 

 Then Section ~\ref{application-section} pulls this
all together in the (mainly combinatorial) 
proof of  Conjecture ~\ref{main-conjecture},
with the most difficult 
 combinatorics  appearing in Lemma ~\ref{stay-regular}. 
Theorem  ~\ref{homeom2-lemma} assembles the various lemmas which together imply that
the complex resulting from our series of collapses is indeed regular and homeomorphic
to a ball.
Finally,  the Fomin-Shapiro Conjecture  is  proven in Theorem ~\ref{main-theorem}.  
Throughout the paper, we deliberately  include a high 
level of detail,  so as  to help readers bridge between the  combinatorics, 
topology, and 
representation theory.

 \subsection{Proof overview}\label{overview-section}
 
Following  Lusztig \cite{Lu},
we  realize the stratified spaces $Y_w$ from Conjecture ~\ref{main-conjecture}
as images of 
maps $f_{(i_1,\dots ,i_d)}$  from polytopes (which in our case 
are simplices) to spaces of matrices.

 Let $(i_1,\dots ,i_d)$ be a 
reduced word for $w\in W$.  Consider the surjective map 
$f_{(i_1,\dots ,i_d)} : \RR_{\ge 0}^d \rightarrow Y_w$ sending $(t_1,\dots ,t_d)$ to the 
product of matrices $x_{i_1}(t_1)\cdots x_{i_d}(t_d)$ where $x_i(t) = I_n + t E_{i,i+1}$ in 
type A, and more generally  $x_i(t) = exp(te_i)$ for $e_i$ a Chevalley generator.  Lusztig
proved that $f_{(i_1,\dots ,i_d)}$ applied to $\RR^d_{>0}$ is a homeomorphism.  
On the other hand, $f_{(i_1,\dots ,i_d)}$ is far from injective on $\RR_{\ge 0}^d$, due 
to the  relations (a) $x_i(u)x_i(v) = x_i(u+v)$ and (b) the type A braid relations
\begin{enumerate}
\item
$x_i(u)x_j(v) = x_j(v)x_i(u)$ for $|j-i|>1$ 
\item
$x_i(a)x_{i+1}(b)x_i(c) = x_{i+1}(\frac{bc}{a+c})x_i(a+c)x_{i+1}(\frac{ab}{a+c})$ for 
$a,b,c>0$ 
\end{enumerate}
and similar relations  $x_i(t_1)x_j(t_2)x_i(t_3)\cdots = x_j(t_1')x_i(t_2')x_j(t_3')\cdots $
of degree $m(i,j)$ in other types, where $m(i,j)$ is the order of $s_is_j$ and 
$(t_1',\dots ,t_{m(i,j)}')$ is obtained from $(t_1,\dots 
,t_{m(i,j)})$  by a change of coordinates 
map as in \cite{Lu}.
We study the image of  $f_{(i_1,\dots ,i_d)}$ restricted to the intersection of 
$\RR_{\ge 0}^d$  with 
the hyperplane $\sum t_i = 1$, denoting this domain by 
$\RR_{\ge 0}^d \cap S_1^{d-1}$.  This has the benefit of being compact while
already reflecting the  full structure of the image of $f_{(i_1,\dots ,i_d)}$  on domain 
$\RR_{\ge 0}^d$.  
This domain is a simplex, with faces  specified 
by which parameters $t_i$ are positive and which are 0.  Since $x_i(0)$ is the identity, 
it is  natural to index 
the faces of the simplex  by the subwords of $(i_1,\dots ,i_d)$.
Lusztig's result  for $\RR_{>0}^d$ together with the above relations 
implies that  $f_{(i_1,\dots ,i_d)}$
restricted to the interior of a face is injective if and only if  the
subword of $(i_1,\dots ,i_d)$ indexing that face is a reduced word.

The above relations will enable us 
to construct
for any face of the simplex indexed by a nonreduced word a family of curves covering the face such that each curve lives in a single fiber of 
$f_{(i_1,\dots ,i_d)}$.  
These curves result from the relations $x_i(u)x_i(v) = x_i(u+v)$ either 
directly or after a suitable series of (long and short) braid moves.
We will collapse each such 
non-reduced face across a family of such curves.  While every 
nonreduced expression will admit a series of braid moves leading to such a  ``stutter'' 
$x_i(u)x_i(v)$, 
a serious challenge to be overcome is that the requisite 
long braid moves give 
change of coordinate maps which a priori 
are not even well-defined on the closures of the cells to be collapsed, let alone 
homeomorphisms on them.  

We get around this by doing certain  other collapses earlier than
a collapse requiring long braid moves.  Specifically, we 
choose the collapsing order so that earlier 
identifications in the boundary of a cell requiring long braid moves will ensure that the 
change of coordinates map will in fact be a well-defined homeomorphism on each 
closed cell to be  collapsed just prior to its collapse.  
The key conceptual lemma behind these change of coordinates maps being
homeomorphisms  is Lemma 
~\ref{single-braid},
while the technical details 
are handled in Lemma ~\ref{collapse-order}. 
 
To see which faces of the simplex should be identified with each other 
in this manner, 
we suppress parameters, replacing $x_i(t)$ by  $x_i$ for each $t$ that is positive, omitting
the letters where $t$ is 0.
We thereby associate a 
so-called $x$-expression to each face.  
An examination of which $x$-expressions correspond to faces having the same image
under $f_{(i_1,\dots ,i_d)}$  
yields the  relations (1) $x_i^2 = x_i$  together  with the  braid
relations (2) $x_ix_j = x_jx_i$ for $|j-i|>1$ and
(3) $x_ix_{i+1}x_i = x_{i+1}x_ix_{i+1}$ in type A.  Going beyond type A, we replace (2) and 
(3) by  analogous long braid 
relations $x_ix_jx_i\cdots = x_jx_ix_j\cdots $ of degree $m(i,j)$ 
for each pair of Coxeter group generators
$\{ s_i, s_j \} $.  In this manner, the (unsigned)  0-Hecke algebra of the Coxeter group $W$ emerges.  Two faces of the simplex will have the same image
under $f_{(i_1,\dots ,i_d)}$ exactly when their  
$x$-expressions represent the same element of the 0-Hecke algebra, or equivalently,
in the language of \cite{KM} if they have the same Demazure product.

In Section ~\ref{topol-collapse-section}, we  introduce a general class of collapsing 
maps  which  may be performed sequentially on a polytope,  
preserving homeomorphism type and regularity on the resulting quotient cell complexes 
at each step.  
 Each such map is defined by first covering a polytope face with a family of parallel line 
 segments across which the face is collapsed, or more generally in subsequent steps 
 covering a cell to be collapsed with a family of curves which we call parallel-like 
 (see Definition ~\ref{parallel-like-def}),  due to  their
 being the image of a family of parallel line segments 
 under a map $g$ with certain convenient properties 
 (such as being  a 
 homeomorphism on the interior of the cell to be collapsed). 
 The following simple example already captures much of the idea of our collapses.

\begin{example}\label{explicit-homeom}
Let
$\Delta_2$ be the convex hull of $(0,0), (1,0), (0,1/2)$ in $\RR^2$, and let
$\Delta_1$ be the convex hull of $(0,0)$ and $(1,0)$ in $\RR^2$.  
We will construct a surjective, continuous function $h: \RR^2 \rightarrow \RR^2$ that
acts homeomorphically on 
$\RR^2 - \Delta_2 $ sending it to $\RR^2 - \Delta_1$.  The idea is to map parallel, vertical line 
segments covering $\Delta_2$ onto their endpoints in $\Delta_1$, then take a 
neighborhood $N$ of $\Delta_2$, specifically a collar for 
$\overline{\RR^2 - \Delta_2}$, and define $h$ in such a way that it
stretches $N$ to cover $\Delta_2$ by mapping extensions of the parallel line segments 
surjectively onto the extended segments.
For $0\le x \le 1$ and $0\le y \le -x/2 + 1/2$, 
let $h(x,y) = (x,0)$. 
For $0\le x \le 1$ and $ - x/2 + 1/2\le y \le 1$, 
let $h(x,y) = (x, \frac{y-1/2+x/2}{1/2 + x/2} )$.
For $-1\le x \le 0$ and $0\le y \le -x/2 + 1/2$, let $h(x,y) = (x,y\frac{-x}{-x/2 + 1/2})$.
For $-1\le x \le 0$ and $-x/2 + 1/2 \le y \le 1$, let $h(x,y) =  (x, -1+2y)$.
Let  $h$
act as the identity outside  $R = \{ (x,y) : -1\le x\le 1, 0\le y\le 1\} $.
\end{example}

\begin{remark}\label{spanier-remark}
See p. 42-43 in Spanier \cite{Sp} for a closely related, though fundamentally different,  homeomorphism also 
given by explicit maps. 
\end{remark}

 Each curve in a family of parallel-like
 curves will have one endpoint in a closed cell $G_1$ in the boundary of the cell $F$ to be
 collapsed,
 and the other endpoint in another closed boundary cell 
 $G_2$.  In Example ~\ref{explicit-homeom}, this is the segment
 from $(0,1/2) $ to $(1,0)$, and  
  the segment from $(0,0)$ to $(1,0)$, respectively. 
 The  collapse given by this family of curves will  
 map each curve to its endpoint in $G_2$, stretching a collar for the closed 
 complement of $F$ within the boundary of a cell of dimension one higher than $F$
 so as to homeomorphically cover  
 $F$ by the part of this closed collar given by $G_1\times [0,1]$.  Convexity
 of each face of the polytope whose image  we are studying 
 will enable this stretching to be accomplished by a continuous map.  Corollary 
 ~\ref{more-general-topol-collapse2} generalizes our collapsing maps somewhat beyond
 polytopes to help accommodate requisite changes of coordinates, using that these 
 collapsing maps may be transferred from one regular CW ball to a homeomorphic one,
 provided that both have the same cell structure on the closed cell to be collapsed, with the 
 homeomorphism of regular CW complex restricting to a cell structure preserving one 
 on this closed cell.

To relate this to our main application, notice e.g. that 
$\{ (t_1,t_2,t_3) \in \RR_{\ge 0}^3 | t_1+t_2 = k_{1,2}\hspace{.05in}{\rm and}\hspace{.05in} t_3=k_3 \} $ 
for the various choices of constants 
$k_{1,2} , k_3$ adding to 1
give parallel line segments covering the simplex  $\{ (t_1,t_2,t_3)\in \RR_{\ge 0}^3| 
\sum t_i=1 \} $ and comprising exactly 
the fibers of the map
$(t_1,t_2,t_3)\mapsto x_1(t_1)x_1(t_2)x_2(t_3)$.
 All of our 
  families of parallel-like curves  will result 
 either directly from stuttering pairs $x_i (u) x_i (v)$ of 
 consecutive letters in a non-reduced expression yielding such parallel line segments, 
 or  as families 
 of curves obtained from such parallel line segments 
 by the composition of a series of 
 earlier collapsing maps  
(that restrict to homeomorphisms on the interior of the 
cell now under consideration)  with the change of 
coordinate homeomorphisms given by the  long braid moves used to create a 
stuttering pair. 

We require  the following properties for parallel-like curves:
 \begin{enumerate}
 \item
 Distinct initial points condition (DIP): the endpoints of the parallel-like curves in $G_1$ are distinct.  
  \item
 Distinct endpoints condition (DE): for each nontrivial
 curve in the collection, 
 its endpoints in $G_1$ and $G_2$ are distinct, by virtue of  open cells in
 $G_1$ and $G_2$, respectively, not having already been identified with each other. 
 \end{enumerate}
 We also  require a combinatorial property of the collapses themselves:
 \begin{itemize}
 \item
 Least upper bound condition (LUB): whenever two cells are identified by the collapse of a cell 
 that is a least   upper bound for  the pair of  them just prior to 
 the collapse, then all cells that are least upper bounds for them just prior to the collapse also 
 get collapsed in that same  step.  
 \end{itemize}

The condition (DIP) is needed 
 for the collapsing map across  parallel-like 
 curves to be well-defined, since the endpoint of each 
 curve in
 $G_1$ (along with the rest of the curve) is mapped onto the other endpoint of the curve,
 which is in $G_2$.
 Condition (DE) will allow us to extend the parallel-like curves across a collar for the 
 closed complement of the cell being collapsed within the boundary of a cell of dimension
 one higher, once we prove such a collar exists.
 The proof that regularity is preserved under our collapses relies heavily on (LUB).

Collapses meeting the  more precisely formulated versions of these
conditions given later will automatically  meet a further  
condition we call the inductive manifold condition, namely  that the closure of the complement of
an $i$-cell within the boundary of an $(i+1)$-cell is a compact topological
manifold with boundary, 
hence has a collar.  This will allow us to extend the collapsing map for 
a low dimensional cell $F$ from a low-dimensional subcomplex where it is  most
naturally defined to our entire complex; this extension process in based upon 
the existence of requisite collars
together with 
the fact that  our particular  collapsing maps  
admit approximations by homeomorphisms;  in fact, we use that 
each comes with a path of homeomorphisms to the identity map, 
enabling each collapsing map to be 
deformed to the identity map across the layers of a collar.

Checking the above conditions for a family of curves  covering a cell to be collapsed
in our main application will rely on a combinatorial analysis of which cells have been 
identified with each other  at the time of each collapse.  
In preparation for these combinatorial arguments,
we develop in Section ~\ref{0-hecke-section}  properties of the 0-Hecke 
algebra, based mainly on  the following convenient
notion: we say that a pair of letters $\{ i_r,i_s \} $  for $r<s$ in a non-reduced 
word $(i_1,\dots ,i_d)$ is a {\it deletion pair} if $(i_r,\dots ,i_{s-1})$ and 
$(i_{r+1},\dots ,i_s)$ are both reduced while $(i_r,\dots ,i_s)$ is non-reduced. 
The 0-Hecke algebra lacks a cancellation law, adding to the challenge of working with
it.  However,  focusing on deletion pairs will enable 
some critical properties of Coxeter groups to be transferred to the 0-Hecke algebra, 
using the fact  that  reduced expressions
in the 0-Hecke algebra are exactly the reduced expressions in the associated 
Coxeter group.

 Section ~\ref{application-section} gives  a particular series of  cell collapses,
 performed sequentially on a simplex $\RR^d_{\ge 0}\cap S_1^{d-1}$,  to 
 produce a regular CW complex  $(\RR_{\ge 0}^d \cap S_1^{d-1})/\sim $
 homeomorphic
 to a ball.   The guiding principle behind our  choice of ordering  
 is that we want a collapsing order amenable to  proof by induction on $d$.  This will
 allow us to 
 assume all our results for smaller $d$; this guarantees when we perform a long
 braid move on the expression associated to a face to be collapsed
 that all possible identifications based on subexpressions of the expression to 
be braided will have already been done.  To this end, we collapse 
faces in an order consistent with linear order on the position of the 
right endpoint of the  leftmost deletion pair in the (highest priority) 
$x$-expression representing that face.  
This implies at the time of the collapse of 
a face whose leftmost deletion 
pair is $\{ i_r,i_s\} $ that all possible point identifications based on letters strictly 
to the left of  position $s$ will have already been done, so that we may apply braid moves 
to the segment from positions $r$ through $s-1$ so as to create a stutter with the letter at 
position $s$. 
 
 In this manner,  we collapse away all faces of the simplex indexed by 
non-reduced words based on point and cell identifications  which are clearly 
necessary.  However, this still 
leaves the challenge of proving that these identifications are enough, that there are 
no remaining instances of two points mapping to the same place under the induced
quotient space map $\overline{f}_{(i_1,\dots ,i_d)}$.  
This is where we turn to Theorem ~\ref{sufficient-theorem},  to prove that  
$\overline{f}_{(i_1,\dots ,i_d)}$ is indeed a homeomorphism
 from $(\RR_{\ge 0}^d \cap S_1^{d-1})/\sim $ to  $Y_w$.

Corollary ~\ref{regular-build} 
sets up the framework in which  Theorem ~\ref{sufficient-theorem}
will be used,  both in this paper and most likely in other applications as well.  
It focuses on images of maps from regular CW complexes having a maximal
cell, with the further requirement that
the map be a homeomorphism on this open big cell.  It singles out 
 conditions 3 and 4 of Theorem ~\ref{sufficient-theorem} to be checked for the particular
 application,  with conditions 1, 2 and 5 then following automatically from the general  set-up.

In the setting of the Fomin-Shapiro Conjecture, condition 3 is immediate from the 
result of Bj\"orner and Wachs that Bruhat order is thin and shellable (\cite{BW}).
The  idea we will use in Lemma ~\ref{injective-condition}
to verify condition 4 is as follows.
Given a reduced expression $s_{i_1}\cdots s_{i_d}$,
the Coxeter group element obtained by deleting a single letter $s_{i_u}$ cannot be the 
same as the Coxeter group element obtained by deleting a letter $s_{i_v}$ for $u\ne v$.  
This combines with  Lusztig's result that $f_{(i_1,\dots ,i_d)}$ restricts to a homeomorphism
on the interior of a cell indexed by reduced word $(i_1,\dots ,i_d)$  to allow 
us to verify codimension one injectivity.   The point is that approaching
the boundary of a cell by letting a single parameter $t_i$ go to 0 as opposed  to 
approaching the boundary by letting a different, individual parameter $t_j$ go to 
0 must give points in distinct cells, hence distinct points.

We will use the map 
$\overline{f}_{(i_1,\dots ,i_d)}: (\RR_{\ge 0}^d \cap S_1^{d-1})/\!\!\sim \hspace{.05in} \rightarrow Y_w$ 
induced from $f_{(i_1,\dots ,i_d)}$ along with its 
restriction to various closed cells 
as the characteristic maps with respect to which we will prove that $Y_w$  is 
a regular CW complex.  
 Our results will imply that $\sim $ identifies exactly those points 
 having the same  image under Lusztig's map 
 $f_{(i_1,\dots ,i_d)}:(t_1,\dots ,t_d)\rightarrow x_{i_1}(t_1)\cdots x_{i_d}(t_d)$ given by any
 reduced word $(i_1,\dots ,i_d)$ for  $w\in W$.    
 We will prove  in Theorem ~\ref{main-theorem}
 that $\overline{f}_{(i_1,\dots ,i_d)}$ is a  homeomorphism from 
 $(\RR_{\ge 0}^d \cap S_1^{d-1})/\!\!\sim $ to $Y_{s_{i_1}\cdots s_{i_d}}$ with 
 $\overline{f}_{(i_1,\dots ,i_d)}$ 
 sending the open cells of $(\RR_{\ge 0}^d \cap S_1^{d-1})/\!\!\sim $ to the cells $Y_u^o$ with 
 $u\le w = s_{i_1}\cdots s_{i_d}$, completing the proof of Conjecture ~\ref{main-conjecture}.

\subsection{Potential further applications} 
\label{flag-section}

Lusztig and Rietsch have studied a 
combinatorial decomposition for the 
totally nonnegative part of a flag variety, namely the decomposition 
into double Bruhat cells (cf. \cite{Lu} and \cite{Ri}).   Lusztig proved 
contractibility of the entire space  in \cite{Lu} while Rietsch and Williams
proved contractibility of each cell closure in \cite{RW}.
Williams conjectured in \cite{Wi} that this is  a regular CW 
complex homeomorphic to a ball.  It seems quite plausible that 
Theorem ~\ref{sufficient-theorem} together with 
tools from Section ~\ref{topol-collapse-section}
could also be used to prove that conjecture, though we believe that
significant further new ideas would also be needed.

Rietsch determined the closure poset of this decomposition 
 in \cite{Ri}.
Williams proved in \cite{Wi} that this poset is shellable and thin,
implying it meets Condition 3 of Theorem ~\ref{sufficient-theorem}.  
Postnikov, Speyer and Williams
proved in \cite{PSW} in the case of the Grassmannian that its 
double Bruhat decomposition
is a CW decomposition; Rietsch and Williams subsequently
generalized this to all flag varieties in
\cite{RW}.  In each case, it remains open whether these CW complexes are regular
and whether the spaces themselves
are homeomorphic to balls. 

\begin{remark}
Williams' conjecture  is related to Conjecture ~\ref{main-conjecture}
in that the stratified spaces we prove to be regular CW complexes
arise as links of cells in the double Bruhat 
stratification of the flag variety.  However, Williams' conjecture does not 
imply Conjecture ~\ref{main-conjecture} since links of cells in regular CW
complexes are not always themselves  regular.  
Consider e.g. the 
double suspension of a Poincare homology 3-sphere with a big cell glued in (personal
communication, Anders Bj\"orner).
\end{remark}

In the case of an arbitrary flag variety, the preimage polytope that has been constructed in 
\cite{RW} is quite abstract, guaranteed to exist by properties of canonical bases.  
In the special case of the Grassmannian, much more  explicit
combinatorics is known about Postnikov's polytope of plabic graphs as well as the 
map from this polytope to the totally nonnegative part of the Grassmannian (see: 
A. Postnikov, Total positivity, Grassmannians, and networks, 
arXiv:math.CO/0609184).


\section{Background and Terminology}\label{background-section}
Now we collect together  basic terminology and
facts  
from topology, combinatorics, the theory of Coxeter groups, and total positivity theory 
that will be essential to this paper.
See e.g. \cite{Bj}, \cite{BB}, \cite{Br}, \cite{Hu}, \cite{Hu2}, 
\cite{Lu}, \cite{Mu}, \cite{RS}, \cite{Sp}, or \cite{St} for further details.

\subsection{Background in topology}\label{top-bg} 
\begin{definition}\label{cw-def}
{\rm 
A {\it CW complex}  is a space $X$ and a collection of disjoint open cells $e_{\alpha }$
whose union is $X$ such that:
\begin{enumerate}
\item
$X$ is Hausdorff.
\item
For each open $m$-cell $e_{\alpha }$ of the collection, there exists a continuous map
$f_{\alpha } :B^m \rightarrow X$ that maps the interior of $B^m$ homeomorphically onto
$e_{\alpha }$ and carries the boundary of $B^m$ into a finite union of 
open cells, each of dimension less than $m$.  
\item
A set $A$ is closed in $X$ if $A \cap \overline{e}_{\alpha }$ is closed in $\overline{e}_{\alpha }$
for each $\alpha $.
\end{enumerate} }
\end{definition}

An {\it open $m$-cell} is any topological space which is homeomorphic to the interior of 
an  $m$-ball $B^m$, with an open $0$-cell being a point.  
The restriction of a characteristic map $f_{\alpha }$ 
to the boundary of $B^m$ is an {\it attaching map}.   Denote the closure of a cell $\alpha $
by $\overline{\alpha }$.
A {\it finite CW complex} is a CW complex with finitely many open cells.

\begin{definition}{\rm 
A  CW complex is {\it regular} with respect to $\{ f_{\alpha } \} $
if additionally each  $f_{\alpha }$ restricts to
a homeomorphism
from the boundary of $B^m$ onto a finite union of lower dimensional open cells. }
\end{definition}

The following (which appears as Theorem 38.2 in \cite{Mu}) 
will enable us to build CW complexes by  induction on dimension.  

\begin{theorem}\label{cw-build}
Let $Y$ be a CW complex of dimension at most $p-1$, let $\sum B_{\alpha }$ be a topological
sum of closed $p$-balls, and let $g: \sum Bd (B_{\alpha }) \rightarrow Y$ be a 
continuous map.  Then the adjunction space $X$ formed from $Y$  and $\sum B_{\alpha }$ by
means of $g$ is a CW complex, and $Y$ is its $(p-1)$-skeleton.
\end{theorem}

\begin{definition} \label{ident-def}
{\rm 
Let $g: X\rightarrow Y$ be a continuous, surjective function.  Then 
the {\it quotient topology} on $Y$ is the topology whose open sets
are the sets whose inverse images are open in $X$.    Say that $g$ is an {\it identification map}
if the topology on $Y$ is the quotient topology  induced by $g$.
}
\end{definition}

The requirement  needed for a continuous, surjective $g:X\rightarrow Y$ to be
an identification map is automatic if $X, Y$ are compact and Hausdorff, which will always
hold for our upcoming collapsing maps. 

\begin{remark}\label{ident-remark}
Given an identification map $g: X\rightarrow Y$ and a function 
$f: X \rightarrow Z$ such that 
$g(x)=g(y)$ implies $f(x)=f(y)$, then Proposition 13.5 of [Br] implies that 
$f$ is continuous iff the induced
function $\overline{f}:Y \rightarrow Z$ satisfying
$f = \overline{f}\circ g$ is continuous.   
\end{remark}

 Remark ~\ref{ident-remark} will allow
us to use continuity of $f_{(i_1,\dots ,i_d)}$ to deduce
continuity of the induced map 
$\overline{f_{(i_1,\dots ,i_d)}}$ on the quotient space after a series of collapses, 
each of which is given by an identification map.

\begin{definition}\label{manifold-def}
{\rm 
A {\it topological $n$-manifold} is a Hausdorff space $M$ having a countable basis of
open sets, with the property that every point of $M$ has a neighborhood homeomorphic
to an open subset of $\HH^n$, where $\HH^n$ is the half-space of points $(x_1,\dots ,x_n)$
in $\RR^n$ with $x_n\ge 0$.  The {\it boundary} of $M$, denoted $\partial M$, is the set
of points $x\in M$ for which there exists a homeomorphism of some neighborhood of $x$
to an open set in $\HH^n$ taking $x$ into $\{ (x_1,\dots ,x_n)| x_n = 0\}  = \partial \HH^n$.
}
\end{definition}

Next we review the notion of a collaring, since this will be critical to our general 
construction of cell collapses 
in Theorem ~\ref{topol-collapse2}.

\begin{definition}\label{collar-def}
{\rm
Given a topological manifold $M$ with boundary,  a {\it collar} or {\it collaring}
for $M$ is a closed neighborhood
$N$ of $\partial M$ contained in $M$ that is homeomorphic to 
$\partial M \times [0,1]$ with $\partial M$ mapping to $\partial M \times \{ 0 \} $.
}
\end{definition}

A proof of the following 
appears  in \cite{Co} and in  Appendix II of  \cite{Vi}.

\begin{theorem}\label{collar-theorem}
If $M$ is a compact,
topological manifold with boundary $\partial M$, then $M$ has a collar.
\end{theorem}

We will  extend collapsing maps from low-dimensional cell boundaries in which they are most naturally defined to higher dimensional cells by proving existence of requisite 
collars and showing that our collapsing maps can be approximated by 
homeomorphisms, in fact constructing a continuous path of homeomorphisms from our 
collapsing map (which  itself 
is not a homeomorphism) to the identity map.

Finally, we briefly recall from \cite{GM} the notion of link for Whitney stratified spaces, and we
refer the reader to  \cite{GM} for further details.  Let
$Z$ be a Whitney stratified subset  of a smooth manifold $M$, let $N'$ be a smooth submanifold
which contains a given point $p$ of $Z$ and which is transverse to each stratum  of $Z$
containing  $p$.  Let  $B_{\delta }(p)$ be the closed ball of radius $\delta $ centered at $p$.
Then  the {\it link}, denoted $L(p)$,  of a stratum $S$ at the point $p$
is the set $L(p) = N' \cap Z \cap \partial D_{\delta }(p)$.

\subsection{Background in topological combinatorics}

\begin{definition}{\rm 
The {\it closure poset} of a finite CW complex is the partially ordered set (or poset) 
of open cells with ]
$\sigma \le  \tau $ iff $\sigma \subseteq \overline{\tau }$.  By convention, we adjoin a 
unique minimal element $\hat{0}$ which is covered by all the 0-cells, which may be 
regarded as representing the empty set. }
\end{definition}

 Let $\partial \tau $ denote $\overline{\tau } \setminus \tau $, i.e. the boundary of 
 $\overline{\tau }$. 

\begin{definition}{\rm 
The {\it order complex} of a finite partially set is the simplicial complex whose $i$-dimensional
faces are the chains $u_0 < \cdots < u_i$ of $i+1$ comparable poset elements.  
}
\end{definition}

A poset is {\it graded} if for each $u\le v$, all saturated chains $u = u_0 \prec u_1 \prec 
\cdots \prec u_k = v$ involve the 
same number $k$ of covering 
relations $u_i \prec u_{i+1}$ (i.e. $u_i < u_{i+1}$ such
that $u_i \le v \le u_{i+1}$ implies $v=u_i$ or $v=u_{i+1}$).  In this case, we say that
the poset interval $[u,v]$ has {\it rank} $k$.  
Recall that a finite, graded poset with unique minimal and maximal elements 
is {\it Eulerian} if  each interval $[u,v]$ 
has equal numbers of elements at even and odd ranks.  This
is equivalent to its
M\"obius function satisfying $\mu (u,v) = (-1)^{rk(v)-rk(u)}$ for each pair
$u<v$, or in other words the order complex of each open interval $(u,v)$ having the same 
Euler characteristic as that of a sphere  $S^{rk(v)-rk(u)-2}$.  A finite, graded poset is
{\it thin } if each rank two  closed interval $[u,v]$ has exactly four elements, in other words if
each such interval is Eulerian.  

\begin{remark}
The order complex of the closure poset of a  finite regular CW complex $K$ 
(with $\hat{0} $ removed) 
is the first 
barycentric subdivision of $K$, hence is homeomorphic to $K$.  In particular, this 
implies that the order complex for any open interval $(\hat{0},v)$ in the closure poset of $K$
will be homeomorphic to a sphere $S^{rk(v)  -2}$.
\end{remark}

In \cite{Bj}, Bj\"orner characterized
which finite, graded posets are closure posets of regular CW complexes, calling such
posets {\it CW posets}:

\begin{theorem}[Bj\"orner]
A finite, graded poset with unique minimal element $\hat{0}$ is the closure poset of a 
regular CW complex if and only if (1) it has at least one additional interval, and (2)
each open 
interval $(\hat{0}, u)$ has order complex homeomorphic to a sphere $S^{rk(u) -2}$.
\end{theorem}

Results of Danaraj and Klee in \cite{DK}  
give a convenient way to verify (2) for a finite, graded poset $P$, 
namely  by proving $P$
is thin and shellable.  

\begin{remark}
Two finite CW complexes may
have the same closure poset in spite of having very different topological structure,  so 
proving that the closure poset of a stratified space is a CW poset gives evidences that the 
stratified space is a regular CW complex, but is not enough to 
determine topological structure of the stratified space itself.
\end{remark}

\begin{definition}
{\rm 
A {\it convex polytope} is the convex hull of a finite collection of points in $\RR^n$, or
equivalently it is an intersection of closed half spaces that is bounded.
}
\end{definition}

For simplicial complexes and polytopes, the closure poset is often called the face poset.
Let $[\sigma ,\tau ]$ denote the subposet consisting of elements $z$ such that 
$\sigma \le z \le \tau $, called the {\it closed interval} from $\sigma $ to $\tau $.  Likewise,
the {\it open interval} from $\sigma $ to $\tau $, denoted $(\sigma ,\tau )$,  is the subposet of
elements $z$ with $\sigma < z < \tau $.  A cell $\sigma $ {\it covers} a cell $\rho $, denoted
$\rho \prec \sigma $, if $\rho < \sigma $ and each $z$ with $\rho \le z \le \sigma $ must
satisfy $z= \rho $ or $z=\sigma $.

For a regular cell
complex in which the link of any cell is also regular,  
$\Delta(u,v)$ is homeomorphic to the link of $u$ within the 
boundary of $v$, hence 
is homeomorphic  to $S^{\dim (v)-\dim (u)-2}$.

\begin{remark}
If each closed interval $[u,v]$ of a finite poset $P$ is Eulerian and shellable, then 
each open
interval has order complex homeomorphic
to a sphere $S^{rk(v)-rk(u)-2}$,  implying 
condition 3 of Theorem ~\ref{sufficient-theorem}.   
\end{remark}

The stratified spaces we consider in our main application have closure posets that are
the intervals of Bruhat order, which were proven to be thin and 
shellable  by Bj\"orner and Wachs in \cite{BW}.

\subsection{Background on Coxeter groups and their associated 0-Hecke algebras}
\label{0-hecke-bg-section}

Let $s_i$ denote the adjacent transposition $(i,i+1)$ swapping the 
letters $i$ and $i+1$ in type A, and more generally denote a  member of a minimal 
set of generators called
the {\it simple reflections} of a Coxeter $W$ 
group by $\{ s_i | i\in I \} $.  Its relations are all of the form $(s_is_j)^{m(i,j)} = e$ with 
$m(i,i)=2$ for all $i$ and $m(i,j)\ge 2$ otherwise.  Finite Weyl groups are all Coxeter groups.

An {\it expression} for a 
Coxeter group element $w$ is a way of writing it as a product of simple reflections $s_{i_1}\cdots 
s_{i_r}$.  An expression is {\it reduced} when it minimizes $r$ among all expressions for
$w$, in which case $r$ is called the {\it length} of $w$.  
An expression $s_{i_1}\cdots s_{i_d}$ may be represented more compactly
by its {\it word}, namely by $(i_1,\dots ,i_d)$.
Breaking now from standard
terminology, we also speak of the {\it wordlength} of a (not necessarily reduced)
expression $s_{i_1}\cdots s_{i_r}$, by which we again mean $r$.  Given simple reflections
$s_i,s_j$, define $m(i,j)$ to be the least positive integer such that  $(s_is_j)^{m(i,j)} = 1$.

The following basic lemma will be key to our proof that the complexes $Y_w$
satisfy Condition 4 in our CW complex regularity criterion:

\begin{lemma}\label{red-lemma}
Given a reduced word $s_{i_1}s_{i_2}\cdots s_{i_r}$ for a Coxeter group element $w$, 
any two distinct subwords of
length $r-1$ which are both themselves reduced must give rise to distinct Coxeter group
elements.
\end{lemma}

We include a short proof of this vital fact for completeness sake.

\begin{proof}
  Suppose deleting $s_{i_j}$ yields the same 
Coxeter group element which we get by deleting $s_{i_k}$ for some pair $1\le j<k\le r$.
This implies $s_{i_j}s_{i_{j+1}}\cdots s_{i_{k-1}} = s_{i_{j+1}}\cdots s_{i_{k-1}} s_{i_k}$.
Multiplying on the right by $s_{i_k}$ yields
$$s_{i_j}s_{i_{j+1}} \cdots s_{i_{k-1}} s_{i_k}  
= s_{i_{j+1}}\cdots s_{i_{k-1}} (s_{i_k})^2 = s_{i_{j+1}}\cdots s_{i_{k-1}},$$
contradicting the fact that the original expression was reduced.
\end{proof}

The expression $s_1s_2s_1$ in the symmetric group demonstrates
that the statement of the above lemma no longer holds if we replace $r-1$ by $r-i$ for
$i>1$.  Thus, it really seems to be quite essential to our proof of the Fomin-Shapiro 
Conjecture that  Theorem ~\ref{sufficient-theorem}
enables us to focus  on codimension one cell incidences.

\begin{lemma}[Exchange Condition, \cite{Hu}]\label{exchange-lemma}
Let $w=s_1\cdots s_r$ (not necessarily reduced) where each $s_i$ is a simple reflection.
If $l(ws) < l(w)$ for some simple reflection $s = s_{\alpha }$, then there exists index $i$
for which $ws = s_1\cdots \hat{s_i}\cdots s_r$.  In particular, $w$ has a reduced 
expression ending in $s$ if and only if $l(ws) < l(w)$.
\end{lemma}

Given a (not necessarily reduced) expression 
$s_{i_1}\cdots s_{i_d}$ for a Coxeter group element $w$, define
a {\it braid-move}  to be the replacement of $s_is_js_i\cdots  $
by $s_js_is_j \cdots $ yielding a new expression for $w$ by virtue of a 
braid relation $(s_is_j)^{m(i,j)} = 1$ with $i\ne j$.
Define a {\it nil-move} to be the replacement of  a substring $s_is_i$ appearing
in consecutive positions by $1$.  We call braid moves with $m(i,j)=2$ {\it commutation
moves} and those with $m(i,j)>2$ {\it long braid moves}. 

\begin{theorem}[\cite{BB}, Theorem 3.3.1]\label{word-property}
Let $(W,S)$ be a Coxeter system consisting of Coxeter group $W$ and  minimal
generating set of simple reflections $S$.  Consider $w\in W$.  
\begin{enumerate}
\item 
Any expression $s_{i_1}s_{i_2}\cdots 
s_{i_d}$ for $w$ can be transformed into a reduced expression for $w$ by a sequence of
nil-moves and braid-moves.
\item
Every two reduced expressions for $w$ can be connected via a sequence of braid-moves.
\end{enumerate}
\end{theorem}

 The {\it Bruhat order} is the partial order on the elements of a Coxeter group $W$ having $u\le v$
 iff there are reduced expressions $r(u),r(v)$ for $u,v$ with $r(u)$ a subexpression
 of $r(v)$.   Bruhat order is also the closure order on the cells $B_w = B^- w B^-$ of the Bruhat 
 stratification of the reductive algebraic group having $W$ as its Weyl group.

Associated to any  Coxeter system $(W,S)$ is a 0-Hecke algebra, with generators 
$\{ x_i | i\in S\}$ and the following relations: for each braid relation 
$s_is_j\cdots = s_js_i\cdots $ in $W$, there is an analogous 
relation $x_ix_j\cdots = x_jx_i\cdots$, again of degree $m(i,j)$;
there are also  relations $x_i^2 = -x_i$ for each $i\in S$.  In our set-up, we will
need relations $x_i^2=x_i$, but this sign change is inconsequential in our setting, so
refer to the algebra with relations $x_i^2 = x_i$ as the (unsigned)
0-Hecke algebra of $W$.  This variation on the usual 0-Hecke algebra has previously 
arisen in work on Schubert polynomials (see e.g. \cite{FSt} or \cite{Ma}).  We refer to 
$x_i^2 \rightarrow x_i$ as a {\it modified nil-move}.  
It still makes sense to speak of reduced and non-reduced expressions, and
many properties (including Lemma ~\ref{red-lemma} and Theorem ~\ref{word-property})
carry over  to the 0-Hecke algebra by virtue of having the same braid moves; 
there are important differences though too, largely resulting from the lack of inverses and a  
 cancellation law.  
 
\subsection{Background in
total positivity theory}\label{tot-pos}

Recall that a real matrix is totally nonnegative (resp. totally positive) if each minor is
nonnegative (resp. positive).
The totally nonnegative part of $SL_n(\RR )$ consists of the matrices in
$SL_n(\RR )$ whose minors are all nonnegative.
Motivated by connections to canonical bases, 
Lusztig generalized this dramatically in \cite{Lu} as follows.
The totally nonnegative  part of a reductive algebraic group
$G$  defined and 
split over $\RR $ is the semigroup generated by the sets
$\{ x_i(t) | t\in \RR_{>0} , i\in I \} , \{ y_i(t) | t\in \RR_{>0} , i \in I \} , $ and 
$\{ t \in T | \chi (t) > 0 \text{ for all } \chi \in X^* (T) \} $, for $I$ indexing the simple roots.  
In type A, we have $x_i(t) = I_n + tE_{i,i+1}$, namely the identity
matrix modified to have the value $t$ in position $(i,i+1)$, and likewise, $y_i(t) = 
I_n + tE_{i+1,i}$.
More generally, $x_i(t) = exp(te_i)$ and $y_i(t) = exp(tf_i)$ for $\{ e_i, f_i | i\in I \} $  the
Chevallay generators.  In other words, if we let $\phi_i$ be the homomorphism of 
$SL_2 $ into $G$ associated to the $i$-th simple root, then
$$
x_i(t) = \phi_i 
\left(
\begin{array}{cc}
1 & t \\ 0 & 1 
\end{array}
\right)  \text{ and }
y_i(t) = \phi_i 
\left(
\begin{array}{cc}
1 & 0 \\ t & 1
\end{array}
\right) .$$

  Let $B^+, B^-$ be opposite Borel subgroups with $N^+ $ (or simply $N$) and
  $N^-$  denoting their unipotent radicals.  In type
  A, we may choose $B^+, B^-$ to consist of the upper and 
  lower triangular matrices in $GL(n)$, respectively.  In this case, $N^+,N^-$ are the
  matrices in $B^+, B^-$ with diagonal entries all equalling one.  The totally nonnegative part of 
  $N^+$, denoted $Y$,  is the submonoid generated
  by $\{ x_i(t_i) | i\in I, t_i \in \RR_{>0} \} $.  
  Let $W$ be the Weyl group of $G$.  One obtains a 
  combinatorial decomposition of $Y$  by taking
  the usual Bruhat decomposition of $G$
  and intersecting each open Bruhat cell $B_w = B^- w B^-$ for $w\in W$ with $Y$ to obtain an 
  open cell $Y_w^o := Y \cap B_w$ in $Y$.  
  We  follow \cite{Lu} in using the standard
 topology on $\RR $ throughout this paper. 
  
\begin{theorem}[Lusztig] For $(i_1,\dots ,i_d)$ any reduced word for $w$  that  
  the map $f_{(i_1,\dots ,i_d)} $ sending $
 (t_1,\dots ,t_d) $ to $ x_{i_1}(t_1)\cdots x_{i_d}(t_d)$ is a homeomorphism 
 from $\RR_{> 0}^d $ to $Y_w^o$ (see Proposition 2.7 in \cite{Lu}). 
 \end{theorem}
  
  The closure of $Y_w^o$, denoted 
 $Y_w$, is  the image of this same map applied to $\RR_{\ge 0}^d$.  Since $x_i(0)$ is the
 identity matrix, the cells in the closure of $Y_w^o$ are obtained by choosing subwords of 
 $(i_1,\dots ,i_d)$, 
 so $Y_w = \cup_{u\le w} Y_u^o$ for $u\le w$ in Bruhat order on $W$. 
 Fomin and Shapiro 
 suggested for each $u < w$ in Bruhat order 
 that the link of the open cell $Y_u^o $ within $Y_w$ 
 should serve as a good geometric model for the Bruhat interval
 $(u,w]$, namely as a naturally arising regular CW complex with $(u,w]$ as its closure 
 poset.  They define $lk(u,w)$ as we describe next.

 Fomin and Shapiro introduced the following projection map
 $\pi_u : Y_{\ge u} \rightarrow Y_u^o$.
    Letting $N(u) = u^{-1}Bu \cap N $ and $N^u = B^-uB^- \cap N$, 
   Fomin and Shapiro proved  that each $x\in Y_{\ge u}$ 
   has a unique expression as $x = x_ux^u$ with 
   $x_u\in N^u$ and $ x^u \in N(u)$.  In light of results in [FS],  $\pi_u(x)$ may be
   defined as equalling    $x_u\in N^u$.
 They defined 
 $lk(u,w)$ as $ (\pi_u^{-1}(x_u))\cap Y_{[u,w]} \cap S_{\epsilon }(x_u)$ for $x_u$ an 
 arbitrary point in $Y_u^o$ and $S_{\epsilon }(x_u)$ a
 sufficiently small sphere about $x_u$ (cf. p. 11 
 in \cite{FS}).    Thus, points of $lk(u,w)$ 
 belong to cells $Y_{u'}$ for $u < u' \le w$, and closure relations are inherited from $Y_w$.   
They proved that each of the proposed
open cells in $lk(u,w)$ is indeed homeomorphic to 
$\RR^n $ for some $n$.

 Recall from  \cite{FZ},  \cite{Lu},   
the relations (1) $x_i(t_1) x_j(t_2) = x_j(t_2)x_i(t_1)$ for any  $s_i,s_j$ which commute,
and  (2) $x_i(t_1)x_j(t_2)x_i(t_3) = x_j(\frac{t_2t_3}{t_1+t_3} )x_i(t_1+t_3)x_j(\frac{t_1t_2}{t_1+t_3}) $
for any $s_i,s_j$ with $(s_is_j)^3 = e$ and any $t_1+t_3 \ne 0$.  These are not 
difficult to verify directly.   In \cite{Lu}, it is proven that there are 
more general relations of a similar nature 
for each braid 
relation  $(s_is_j)^{m(i,j)} = e $ in $W$, i.e., relations 
$x_i(t_1)x_j(t_2)\dots = x_j(t_1')x_i(t_2')\dots $
of degree $m(i,j)$
for $t_1',\dots t_{m(i,j)}'$ rational functions of $t_1,\dots t_s$ each mapping $\RR_{>0}^d$ to $\RR_{>0}$.

\begin{lemma}\label{same-sum}
The new parameters after applying a braid relation will 
have the same sum as the old ones; moreover, this preservation of sum
refines to the subset of parameters for any fixed $x_i$. 
\end{lemma}

\begin{proof}
This follows from the description
of $x_i(t)$ as $exp(te_i)$, simply by comparing the linear terms in the expressions
$x_i(t_1)x_j(t_2)\dots = x_j(t_1')x_i(t_2')\dots $ appearing in a braid relation.
\end{proof}

Thus, our  description of $\RR_{\ge 0}^d \cap S_1^{d-1}$ 
may be used  even after a change of coordinates (as in Lemma ~\ref{collapse-order})
resulting from  a braid relation.

\section{A new regularity criterion for CW complexes}
\label{tool-section}

Before proving Theorem
~\ref{sufficient-theorem}, the new regularity criterion,
we first give a few examples demonstrating the need for  its various hypotheses.  It 
seems likely that this new criterion will mainly apply to images of regular CW complexes.

\begin{example}
The CW complex consisting of
a 2-cell with its entire boundary attached to a 0-cell violates 
condition 2 of Theorem ~\ref{sufficient-theorem}.   Condition 2 is 
designed also to 
preclude examples such as  a CW complex
whose 1-skeleton is the simplicial complex comprised of the
faces $\{ v_1,v_2,v_3,e_{1,2},e_{1,3},e_{2,3} \} $, also having
a  two cell with a nontrivial closed interval of its boundary 
mapped to $v_2$ and the 
remainder of its boundary mapped homeomorphically to the rest of the 1-skeleton.  
\end{example}

\begin{remark}
In the latter example above, one may choose a different characteristic map which is a 
homeomorphism.  Whether  this can always be done for finite
CW complexes graded by cell dimension and satisfying conditions 1, 3, 4, and 5
seems subtle at best, in light of examples such as the Alexander horned sphere: a ball
which 
cannot be contracted to a point without changing the topology of its complement, since
that is not simply connected. 
\end{remark}

The next two examples give non-regular CW complex
satisfying conditions 1, 2,  4, and 5  of Theorem
~\ref{sufficient-theorem}, but violating condition 3.  The first example violates thinness,
while the second one violates the requirement that open intervals of rank at least 3 be
connected.

\begin{example}\label{non-example}
Let $K$ be a 2-dimensional CW complex whose 1-skeleton is the simplicial 
complex with maximal faces $\{ e_{1,2}, e_{2,3}, e_{1,3}, e_{3,4}, e_{4,5}, e_{3,5} \} $ and
which has a unique $2$-cell $\sigma $.  The boundary of $\sigma $
is mapped by $f_{\sigma }$
to the 1-cycle $ (e_{3,1}, e_{1,2}, e_{2,3}, e_{3,4}, e_{4,5}, e_{5,3})$.
The attaching map $f_{\sigma }$ sends to different points of the boundary of
$\sigma $ to $v_3$.
\end{example}

\begin{example}
For a non-regular CW complex  
satisfying conditions 1, 2, 4, and 5 of Theorem ~\ref{sufficient-theorem}
as well as thinness of the closure poset, but violating the connectedness requirement
for open intervals of rank at least 3, take a 3-dimensional cube and glue together a 
pair of antipodal vertices.
\end{example}

One might ask if condition 3  
could be replaced by  the requirement that the closure poset be Eulerian, i.e. whether
this could replace the connectedness part of the requirement.  Closure posets do have 
the feature that open intervals $(\hat{0},u)$ with $rk(u)>2$ are connected, by virtue of
the fact that the image of a continuous map from a sphere $S^d$ with $d>0$ is 
connected.    However, there are Eulerian closure posets of  
CW complexes having 
disconnected intervals $(u,v)$ with $rk(v) - rk(u) >2$ (personal communication, Hugh
Thomas).   Still, it seems plausible  that
condition 3 in Theorem ~\ref{sufficient-theorem} might be replaceable by  the
Eulerian property for the closure poset, i.e. that this together with conditions 1, 2, 4, and 5 
could imply regularity.

Next is a non-regular CW decomposition of $\RR P_2$ 
satisfying conditions 1, 2, 3, and 5 of Theorem
~\ref{sufficient-theorem}, but failing condition 4.

\begin{example}\label{rp2}
Let $K$ be the CW complex having as its  1-skeleton 
the simplicial complex with maximal faces $e_{1,2}, e_{2,3},e_{1,3}$.  Additionally,
$K$ has a 
single 2-cell whose boundary is mapped to the 1-cycle which goes twice around the 
1-cycle $(v_1,v_2,v_3)$.  Notice that this CW decomposition of $\RR P_2$ has the same
closure poset as a $2$-simplex, but the attaching map for the $2$-cell is a 2 to 1 map onto
the lower dimensional cells. 
 \end{example}

Finally, we give an example (due to David Speyer) of a CW complex with characteristic
maps meeting conditions 1, 2, 3 and 4,  but failing 
condition 5, though this CW complex is regular with 
respect to a different choice of characteristic maps.
David Speyer also helped with the formulation of condition 5.

\begin{example}
Let the 2-skeleton be the boundary of a pyramid.  Now attach a 
3-cell which is a triangular prism by sending an entire edge of one of the rectangular 
faces to the unique vertex of 
degree 4 in the pyramid, otherwise mapping the boundary of the prism 
homeomorphically to the boundary of the pyramid.
\end{example}  

\begin{proposition}
Conditions 1 and 2 of Theorem ~\ref{sufficient-theorem}
imply that the closure poset is graded by cell dimension.
\end{proposition}

\begin{proof}
Consider any $e_{\rho } \subseteq \overline{e_{\sigma }}$ with $\dim (e_{\sigma }) - \dim (e_{\rho }) > 1$.
Choose a point $p$ in $e_{\rho }$  expressible as $f_{\sigma}(x) $ for some $x\in S^{\dim
 (e_{\sigma } )-1}$.
If we take an infinite series of smaller and smaller open sets about $x$, by Condition 2 each 
must include a point sent by $f_{\sigma }$ to an open cell of higher dimension than $e_{\rho }$; 
finiteness of the CW complex then implies some such open cell $e_{\tau }$ is mapped into 
infinitely often, implying $p \in \overline{e_{\tau }}$.  Thus,
$e_{\rho }< e_{\sigma }$ for $\dim (e_{\sigma }) - \dim (e_{\rho } ) > 1$ implies there exists 
$e_{\tau }$ with $e_{\rho } < e_{\tau } <  e_{\sigma }$.
\end{proof}

This motivates us to say that 
a finite CW complex is {\it dimension-graded} whenever it meets 
conditions 1 and 2 of Theorem ~\ref{sufficient-theorem}.
Now to the proof of Theorem ~\ref{sufficient-theorem}.

\begin{proof}
Conditions 1, 2, and 4 are each necessary tautologically.  The necessity
of 3 follows easily from the fact that a regular CW complex is homeomorphic to
the order complex of its closure poset.   To see that 5 is also necessary, 
note that if $K$ is regular with respect to the characteristic maps $\{ f_{\alpha } \} $, then
$e_{\sigma } \subseteq \overline {e_{\tau }}$ implies that $f_{\sigma }$ factors as 
$f_{\tau } \circ f_{\tau }^{-1}|_{\sigma } \circ f_{\sigma }$ where 
$ f_{\tau }^{-1}|_{\sigma } \circ f_{\sigma }$ is the desired embedding.

Now  to the
sufficiency of these five conditions.  
We must prove that each attaching map $f_{\sigma }$ 
is a homeomorphism from $\partial (B^{\dim \sigma })$
to the set of open 
cells comprising  $\overline{e_{\sigma }} \setminus e_{\sigma }$.  Since $K$ is a 
CW complex in which the closure of each cell is a union of cells,  $f_{\sigma }$ must be 
continuous and surjective onto a union of lower dimensional
cells, leaving us to  
prove injectivity of $f_{\sigma }$ and continuity of $f_{\sigma }^{-1}$. 
However, once we  prove injectivity, we may use the 
fact that any bijective, continuous map from a compact set to a Hausdorff space is a 
homeomorphism to conclude continuity of the inverse, so it suffices to prove
injectivity.

If the attaching maps for $K$ were not all injective, then we could 
choose open cells $e_{\rho },e_{\sigma }$ with $ \dim (e_{\sigma }) - \dim (e_{\rho } )$
as small as possible such that $e_{\rho } \in \overline{e_{\sigma }}$ and 
 $f_{\sigma }$ restricted to the  preimage of $e_{\rho }$ is
not 1-1.  Then we could choose  a point $z \in e_{\rho }$ with
$| f_{\sigma }^{-1}(z)| = k$ for some $k>1$.  
By condition 4, $\dim (e_{\sigma }) - \dim (e_{\rho } )$ must be at least 2.
We will now show 
that the open interval $(e_{\rho },e_{\sigma })$ in the closure poset has 
at least $k$ connected components, which by condition 3
forces $[e_{\rho },e_{\sigma }]$ to have rank exactly two.   The point is to show for
each point  $p_i \in f^{-1}_{\sigma }(z)$
that there is an open cell $e_{\tau_i } \subseteq \overline{e_{\sigma }}$ such that
$p_i \in \overline{\iota ( B^{\dim e_{\tau_i }})} $, 
 and then to show for distinct points $p_i, p_j \in f^{-1}_{\sigma }(z)$
 that the open cells $e_{\tau_i },e_{\tau_j }$ are incomparable in the closure poset.
 To prove the first part, 
take an infinite sequence of smaller and smaller balls about $p_i$, which by
condition 2 must each intersect  $f^{-1}_{\sigma }(e_{\tau })$ for some 
$e_{\tau }  < e_{\sigma }$ with $\dim e_{\sigma }- \dim e_{\tau }= 1$; 
by finiteness of $K$, the preimage of some such  $e_{\tau_i } $ is hit infinitely
often, implying $p_i \in \overline{f^{-1}_{\sigma }(e_{\tau_i })}$, hence $e_{\rho } 
\subseteq \overline{e_{\tau_i }}$.  
We prove next 
that the collections of cells whose closures
contain the various points in $f^{-1}_{\sigma }(z) $ must belong to distinct 
components of $(e_{\rho },e_{\sigma })$, yielding the desired $k$ components in the open
poset interval. 
Consider $p_1\ne p_2$ with $p_i \in \overline{f^{-1}_{\sigma }(e_{\tau_i })}$ for $i=1,2$.
If  $e_{\tau_i } < e_{ \tau_j}$ in the closure poset for $\{ i, j \} =  \{ 1,2 \} $, then
condition 5 would 
imply $ \overline{f^{-1}_{\sigma }(e_{\tau_i})} \subseteq \overline{f^{-1}_{\sigma }(e_{\tau_j}) }$,
and hence $p_1,p_2 \in  \overline{f^{-1}_{\sigma }(e_{\tau_j} )}$, contradicting the fact that
$f_{\tau_j }$ restricted to the preimage of $\rho $ is a homeomorphism.  Thus, $(e_{\rho } ,
e_{\sigma })$ has no comparabilities between cells whose preimages under $f_{\sigma }$ 
have closures containing distinct points of $f^{-1}(z)$; in particular, $(e_{\rho } ,e_{\sigma })$ has at least $k$ connected components, hence must be  rank two.

Finally, we show 
that $(e_{\rho },e_{\tau })$ has at least $2k $ elements,  forcing $k$ to be $1$, by the thinness
requirement in condition 3.  This will contradict our assumption that $k$ was strictly larger
than 1.  Lemma ~\ref{extra-lemma} provides the desired $2k$ elements by  showing that
 for each of the $k$ preimages of $z$, there are  at least two open cells $e_{\tau } $
in $(e_{\rho } ,e_{\sigma })$ with $\overline{f^{-1}_{\sigma }(e_{\tau })}$ 
containing that particular preimage of $z$.
\end{proof}

\begin{lemma}\label{extra-lemma}
If a CW complex $K$ meets the  conditions of Theorem ~\ref{sufficient-theorem},
then it also satisfies the  following condition:
for each open cell 
$e_{\tau }$ and each $x \in \overline{e_{\tau }} \setminus e_{\tau }$
with $f_{\tau }(x)$  in an open cell $e_{\rho }\subseteq \overline{e_{\tau }}$ with
$\dim e_{\tau }- \dim e_{\rho } = 2$, there exist distinct open 
cells $e_{\sigma_1 },e_{\sigma_2 }$ with $\dim e_{\sigma_i }= 1 + \dim e_{\rho }$ 
and $x \in  \overline{f^{-1}_{\tau }(e_{\sigma_i})}$ for $i=1,2$. 
\end{lemma}

\begin{proof}
 Condition 2 ensures that the boundary of $B^{\dim e_{\tau }}$ does not include
any open $(\dim e_{\tau } - 1)$-ball, all of 
whose points map are mapped by 
$f_{\tau }$ into
$e_{\rho }$.  In particular, each such ball containing $x$ includes points
not sent by $f_{\tau }$ to  $e_{\rho }$.  Since $K$ is finite, there must be some particular cell
$e_{\sigma_1 }$ such that points arbitarily close to $x$ within the boundary of $B^{\dim \tau }$ 
map into  $e_{\sigma_1}$, implying $x \in \overline{e_{\sigma_1}}$, with $\dim e_{\rho } < \dim 
e_{\sigma_1} <
\dim e_{\tau }$.   Thus,  $e_{\rho }\subseteq \overline{e_{\sigma_1 }}$ and
$\dim e_{\sigma_1 } = \dim e_{\rho } + 1$, just as needed.

Now let us find a suitable $e_{\sigma_2 }$.   
Here we
use the fact that removing the boundary of $e_{\sigma_1}$ from a sufficiently small ball 
$B^{\dim e_{\tau } - 1}$ about
$x$ yields a disconnected region, only one of whose components may include points from 
$e_{\sigma_1 }$.  This forces the existence of the requisite 
open cell $e_{\sigma_2}$ which includes points of the
other component and has $x$ in its closure.  
\end{proof}

We will actually use Theorem ~\ref{sufficient-theorem} within the following framework:

\begin{corollary}\label{regular-build}
Let $K$ be a finite, regular CW complex of dimension $p$ 
and let $f$ be a continuous function from $K$ to a Hausdorff space $L$.  Suppose
that $f$ is a 
homeomorphism on the interior of each open cell and on the closure of each cell of the 
$(p-1)$-skeleton of $K$.
 Then $f(K)$ is a finite CW complex satisfying conditions 1, 2, and 5 of Theorem
 ~\ref{sufficient-theorem}, with the restrictions of $f$ to various closed cells in
 $K$ serving as the characteristic maps.
 \end{corollary}

\begin{proof} 
The  restrictions of $f$ to a collection of 
closures of
cells of the $(p-1)$-skeleton give the characteristic maps needed to prove that
the $(p-1)$-skeleton of $f(K)$ is a finite CW complex.
Now we use Theorem ~\ref{cw-build} to attach the $p$-cells and 
deduce that $f(K)$ is a finite CW complex with characteristic
maps given by the various restrictions of $f$.  

Conditions 1 and 2 are immediate  from our assumptions on $f$.
If there are two 
open cells $\sigma_1,\sigma_2$  in $K$ (of dimension at most $p-1$) with identical 
image under $f$, then the fact that $\overline{\sigma }_1 $ and $\overline{\sigma }_2 $ are
both regular with isomorphic closure posets gives  
a homeomorphism from $\sigma_1$ to $\sigma_2$ preserving cell structure, namely
the map sending each $x$ to the unique $y$ with $f(y)=f(x)$.
This allows us to use the embedding of
either $\overline{\sigma }_1$ or $\overline{\sigma }_2$ 
in the closure of any higher cell of $K$ to deduce Condition 5.  We use 
that $L$ is Hausdorff and that we have finitely many cells 
to deduce requirements (1) and (3) of  CW complexes.  
\end{proof}

Although the requirements of Corollary ~\ref{regular-build} may seem quite demanding,
Corollary ~\ref{regular-build}
is well-suited to proving for a family of regular CW complexes  that their
images under $f$ are also  regular CW complexes 
by  an induction on dimension.  
We will use Corollary ~\ref{regular-build} in our main application  in the 
proof of Theorem ~\ref{main-theorem},  a key inductive result.

\section{Topological collapsing lemmas} 
\label{topol-collapse-section}

In this section, we introduce certain types of collapses that may be
performed sequentially on a 
convex polytope, and we prove that these
preserve homeomorphism type as well as the property of being a 
finite regular CW complex, though they do not preserve polytopality.  
In Theorem ~\ref{topol-collapse2}, we  explicitly define the maps accomplishing
these  collapses.  Then we give a relaxation of the requirements of Theorem 
~\ref{topol-collapse2} so as to enable the transfer of a parametrization function on curves
in  a closed
cell $\sigma $ 
meeting the requirements of Theorem ~\ref{topol-collapse2} to a closed cell with the same
cell structure in a different, homeomorphic regular CW complex which does not necessarily
have the same cell structure outside of $\sigma $.    

In preparation for these results, we first 
introduce some helpful properties a topological 
space or a map may have, with the commonality 
that in practice 
these properties may be verified using primarily combinatorics.
 In what follows, we will typically have a topological space $X$ endowed with the
 structure of a finite, regular CW complex $K$.  We denote this by $X_K$
 and  call $X_K$ a regular CW space.
  Denote by 
$K/(\ker g)$ the quotient space obtained from an identification map $g$ on a 
topological space $K$ by 
setting $x\sim y$ iff $g(x)=g(y)$.

We begin  by establishing 
the following convenient general notion of collapse which is much in the spirit of the 
concept of elementary collapse, but in fact will preserve homeomorphism type rather
than just homotopy type.   This  section will then be devoted to 
developing a specific set of checkable conditions which will yield such a collapse.

\begin{definition}\label{collapse-map-defn}
{\rm 
Given a finite regular CW complex $K$ on a set $X$ 
and an open  cell $L$ in $K$, 
define a {\it face collapse} or {\it cell collapse} of
$\overline{L}$ onto $\overline{\tau }$ for $\tau $ an open
cell contained in $\partial {L} $ to be an identification map $g: X\rightarrow X$ 
such that:
\begin{enumerate}
\item
Each open cell of $\overline{L}$ is mapped
surjectively onto an open cell of $\overline{\tau }$ with $L$ mapped onto $\tau $
\item
$g$ restricts to a homeomorphism from $K \setminus \overline{L}$ 
to $K  \setminus \overline{\tau }$ and acts homeomorphically on $\overline{\tau }$.
\item
The images under $g$ of the cells of $K$ form a regular CW complex
with new characteristic
maps obtained by composing the original characteristic maps of $K$ with 
$g^{-1}:X_K \rightarrow X_K$ for
those  cells of $K$ contained either in $\overline{\tau  }$ or in $ (K \setminus \overline{L} )$.
\end{enumerate}
We call such a  map  $g$ a {\it collapsing map}.
Remark ~\ref{collapse-preserve-homeom} will
show that these collapses preserve homeomorphism
type.
}
\end{definition}

The collapse of a cell $L$ often will also collapse other cells in its closure in the process.

\begin{remark}\label{collapse-preserve-homeom}
The induced map $\overline{g} : X_K/(\ker g) \rightarrow X$ 
is continuous by Remark ~\ref{ident-remark},  and it is  
bijective by how it is defined.  
Since $K$ is Hausdorff and $K/(\ker g)$ is compact, we may 
conclude that $\overline{g}$ is a homeomorphism.  Thus, face collapses 
preserve homeomorphism type.
\end{remark}

\begin{remark}\label{stay-collapse}
If $g$ is a collapsing map on $X_K$, and $h$ is a homeomorphism from the
underlying space $X$ to itself, 
then $g\circ h$, composing functions right to left, is a collapsing map as well.
\end{remark}

For example, Remark ~\ref{stay-collapse} may enable a change of coordinates for 
convenience prior to a collapse.
The following 
may also be helpful for controlling how such
a homeomorphism such as a change of coordinates map may 
interact with subsequent cell  collapses. 

 \begin{lemma}\label{quotient-lemma}
Suppose $K_1$ and $K_2$ are topological spaces, 
$f$ is a homeomorphism from $K_1$ to $K_2$, 
$\pi_1:K_1\rightarrow K_1$ and $\pi_2: K_2\rightarrow K_2$ are identification
maps giving rise to
quotient spaces $K_1/\sim $ and $K_2/\sim '$.  If we also have
$x\sim y$ in $K_1$ iff $f(x)\sim ' f(y)$ in $K_2$, then $K_1/\sim $ is homeomorphic to
$K_2/\sim '$ under  the induced map $\overline{f}$. 
\end{lemma}

\begin{proof}
This follows easily  from Proposition 13.5 in Ch. 1 of \cite{Br} by constructing a suitable 
commutative diagram.   
\end{proof}

\begin{definition}\label{length-fn}{\rm 
Given a collection of parallel line segments $\mathcal{C}$
covering a face $F$ of a polytope $P$, define a length function
$len: F\rightarrow \RR $ by letting $len(x)$ be the length of the element of 
$\mathcal{C}$ containing $x$.  Now define a {\it parametrization} 
$p: F \rightarrow [0,1]$ by letting the 
restriction of $p$ to any $c\in \mathcal{C}$ be the linear function from 
$c$ to $[0,1]$. }
\end{definition}

\begin{remark}
Convexity of $P$ implies $len $ is continuous, which in turn 
implies that $p$ is also continuous on $int(F)$ and everywhere on $F$ except at points
comprising parallel line segments consisting of just single points, though here nearby 
parallel line segments also approach length 0.
Moreover, if $h$ is any homeomorphism from $int(F)$ to another topological
space, then $p\circ h^{-1}$ is also continuous.  
\end{remark}

Continuity of $\it len $ will allow us to perform 
collapses across families of parallel-like curves,  as defined next, 
according to a continuous
function that stretches a collar to cover a face being collapsed.  This
stretching function
is provided in Theorem ~\ref{topol-collapse2}.

\begin{definition}\label{parallel-like-def}{\rm 
Let $K_0$ be  a  convex polytope, and let $\mathcal{C}_i^0$ be a 
family of parallel line
segments covering a closed  face $L_i^0$ in $\partial K_0$ with the elements
of $\mathcal{C}_i^0$ given by linear functions $c: [0,1]\rightarrow L_i^0$.
Suppose that there
is a pair of closed faces $G_1,G_2$ in $\partial L_i^0$ with 
$c(0)\in G_1$ 
and $c(1)\in G_2$ for each $c\in \mathcal{C}_i^0$ and
there is a  composition $g_i\circ \cdots \circ g_1$ of face collapses  (cf. Definition
~\ref{collapse-map-defn}) on
$K_0$ such that:
\begin{enumerate}
\item
$g_i\circ \cdots \circ g_1$ acts homeomorphically on $int(L_i^0)$.
\item
For each  $c\in \mathcal{C}_i^0$, 
$g_i \circ \cdots \circ g_1$ either sends $c$ to a  single 
point or acts homeomorphically on $c$.
\item
Suppose $g_i\circ \cdots \circ g_1(c(t)) = g_i \circ \cdots \circ g_1(c'(t'))$ 
for $c\ne c'\in \mathcal{C}_i^0$ and some  $(t,t')\ne (1,1)$.   Then
$t=t'$, and for each $t\in [0,1]$ we have 
 $g_i\circ\cdots \circ g_1(c(t)) = g_i \circ\cdots\circ g_1(c'(t))$.
\end{enumerate}
Then call  $\mathcal{C}_i = \{ g_i\circ\cdots\circ g_1(c) | c\in \mathcal{C}_i^0 \} $ 
a family of {\it parallel-like} curves on the closed cell 
$L_i = g_i\circ\cdots\circ g_1(L_i^0)$ of the finite regular CW complex 
$K_i = g_i\circ \cdots\circ g_1(K_0)$. }
\end{definition}

\begin{remark}
These conditions are designed so that they  only need to be checked  just prior to the 
$k$-th collapsing step
for the curves covering a cell to be  collapsed in the $k$-th collapse.
\end{remark}

Notice that  Definition ~\ref{parallel-like-def}, part (3), 
 implies the curves are nonoverlapping except perhaps at their endpoints 
 in $g_i\circ\cdots\circ g_1(G_2)$.  Verifying (3) mainly requires showing 
 curves have distinct endpoints in $g_i\circ\cdots \circ g_1 (G_1)$, leading to the following:
 
 \begin{condition}\label{dip-def}
 Let us  call  Definition ~\ref{parallel-like-def}, part (3), the 
 {\it distinct initial points condition} (DIP). 
 \end{condition}

\begin{remark}\label{why-four}
In practice, we will verify Definition ~\ref{parallel-like-def}, part (2), by verifying 
Condition ~\ref{distinct-endpoint-condition} below; this suffices  because 
$g_i\circ \cdots\circ g_1$ acts homeomorphically on each open cell not
collapsed by any $g_j$ for $j\le i$.
\end{remark}

We say that a curve is {\it nontrivial} if it includes more than one point.

\begin{condition}
\label{distinct-endpoint-condition}{\rm 
A collection $\mathcal{C}$
of curves covering the closure of a cell $F$ satisfies the {\it distinct endpoints
condition} (DE) if for any nontrivial
curve $c\in \mathcal{C}$, the two endpoints of $c$ live
in distinct cells in $\overline{F}$. }
\end{condition}

The fact that each collapse restricts to a homeomorphism on the interior of each cell
it does not collapse allows us to reduce what actually must be checked to show that
a family of curves is parallel-like to the following:

\begin{remark}
To  verify the requirements of Definition ~\ref{parallel-like-def} for the family of curves 
$\mathcal{C}_i$ covering a cell $L_i$ to be collapsed in the $i$-th collapsing step, it
suffices to show:
\begin{enumerate}
\item
Each earlier collapse restricts to a  homeomorphism on each open cell which is 
not collapsed prior to the $i$-th collapsing step
\item
(DIP) for $\mathcal{C}_i$ holds just prior to the $i$-th collapse
\item
(DE) for $\mathcal{C}_i$ holds just prior to the $i$-th collapse
\item
Consistency of parametrizations when entire curves in $\mathcal{C}_i$ that are 
in $\partial L_i$  are identified with each other in earlier collapses.
\end{enumerate}
\end{remark}

In our main application, and quite possibly in other future applications as well, the last 
condition will follow immediately 
from our set-up; more specifically, it will follow from how each collapse replace a 
pair of parameters $t_{r_k-1}'$ and $t_{r_k}'$ by a single parameter equalling their
sum.   Next we give the ingredients that will allow us to 
extend our collapsing maps through collars, culminating in 
Condition ~\ref{manifold-property} and Lemma ~\ref{bigger-collapse}.

 \begin{definition}\label{interpolating-def}{\rm
Given a continuous, surjective function $g_{i+1}: X\rightarrow X$, 
define an {\it interpolating family} of maps 
$\{ g_{i+1,t} | t\in [0,1] \} $ from  $X$ to $X$ as a collection of maps with  
$g_{i+1,0} = g_{i+1}$ and  $g_{i+1,1} = id $, 
requiring (1) for each $t\in (0,1]$
that $g_{i+1,t}$ be a homeomorphism from $X$ to $X$ 
and (2) continuity of   the map $h_{i+1}: X \times [0,1] \rightarrow X$ defined by 
$h_{i+1}(x,t) = g_{i+1,t}(x)$. }
\end{definition}
 
 \begin{remark}
 This notion of interpolating family of maps is
 almost exactly the topological notion of isotopy, except that our 
 initial map $g_{i+1,0}$ is not a homeomorphism.  In our setting, 
 $g_{i+1,0}$ is a 
 collapsing map which will admit approximations   by homeomorphisms.
 \end{remark}
 
\begin{lemma}\label{interp-from-curves}
Suppose a collapsing map  $g$ collapses 
across a family of parallel-like curves $\mathcal{C}$  where 
each $c\in \mathcal{C}$ is sent to itself  by a  monotonically 
increasing, piecewise linear function
$g: [0,1]\rightarrow [0,1]$. 
Then $g$ gives rise to an interpolating family of maps.
\end{lemma}

\begin{proof}
By definition, we must have $0 =  a_1 < \cdots < a_k =1$ and 
$0 = b_1 \le \cdots \le b_k = 1$ such that $g$ maps $[a_i,a_{i+1}]$ to 
$[b_i,b_{i+1}]$ by a linear map for each $1\le i \le k-1$.  Then define 
$g_t$ instead to map $[a_i,a_{i+1}]$ linearly to 
$[ta_i + (1-t)b_i, ta_{i+1} + (1-t)b_{i+1}]$.
\end{proof}

\begin{lemma}\label{extend-interior}
If a collapsing map $g_{i+1}$ on a regular CW space $X_K$ 
homeomorphic to a sphere  gives 
rise to an interpolating family of maps $\{ g_{i+1,t} | t\in [0,1]\} $ on 
$X_K$,  then $g_{i+1}$ extends to a collapsing
map on any ball $B$ having $X_K$ as its boundary, implying 
$g_{i+1}(B)$ is  a regular CW complex
homeomorphic to a ball.
\end{lemma}

\begin{proof}
Choose a collar $X\times [0,1]$ for $B$, which exists
by Theorem ~\ref{collar-theorem}.
For each $(x,t)\in X\times [0,1]$, let $g_{i+1}(x,t) = g_{i+1,t}(x)$.
Let $g_{i+1}$ act as the identity map on $B\setminus (X\times [0,1])$.
\end{proof}

\begin{condition}\label{manifold-property}{\rm 
A  finite
regular CW complex $K$ has the {\it inductive manifold condition} (IM) if for each
open $d$-cell $\tau $ and each $(d+1)$-cell $\sigma $ such that 
$\tau \subseteq \partial \sigma $, 
$\overline{\partial \sigma \setminus \tau } $
is a compact manifold with boundary.}
\end{condition}

Next we show how to use collars 
to extend a collapsing map from the boundary of 
a low-dimensional closed cell  to an entire cell complex.

 \begin{lemma}\label{bigger-collapse}
Let $K$ be a regular CW complex having a unique maximal cell
and satisfying Condition ~\ref{manifold-property} (IM).
Let $\tau $ be an $i$-cell  in the boundary of an $(i+1)$-cell $\sigma $.  Let
$g_{i+1}$  be a collapsing map on $\partial\sigma $ that collapses the cell $\tau $.  
If $g_{i+1}$ gives rise to an interpolating family of
maps, then $g_{i+1}$ may be extended to a collapsing map on $K$.
 \end{lemma}
 
 \begin{proof}
First choose a series of cells $\tau = \sigma_1 , \sigma = \sigma_2 , 
\sigma_3 , \dots ,
\sigma_k $ such that $\sigma_j \subseteq \overline{\sigma_{j+1}} $ with
$\dim \sigma_{j+1} = \dim \sigma_j + 1 $ for each $j$, letting $\sigma_k$ 
be the unique 
maximal cell of $K$.  We are given $g_{i+1}$ defined on $\partial \sigma $
and will now 
describe for each $2\le j\le k-1$ how to extend $g_{i+1}$ from $\partial \sigma_j$ to 
$\overline{\sigma_{j}}$ and then to $\partial \sigma_{j+1}$.

Lemma ~\ref{extend-interior} enables us  to extend $g_{i+1}$ from 
$\partial \sigma_j $ to  $\overline {\sigma_j } $. 
 If $\sigma_j $ is the big cell, we are done.  Otherwise,
 choose $\sigma_{j+1}$
 with  $\overline{\sigma_j }\subseteq \partial \sigma_{j+1}$ and 
 take a collar for $\overline{\partial \sigma_{j+1} \setminus \sigma_j }$ within
 $\partial \sigma_{j+1} $, which exists by Condition ~\ref{manifold-property} and
 Theorem ~\ref{collar-theorem}.  Use an
 interpolating family to extend $g_{i+1}$ from $\overline{\sigma_j }$ to 
  $\partial \sigma_{j+1} $, defining the interpolating family as follows.  
 
Assuming $g_{i+1}$ has been defined on the first $r-1$ collars  in 
the above alternation,  let $g_{i+1}(x,t_1,t_2,\dots ,t_r) 
= (g_{i+1,1-(1-t_1)\cdots (1- t_r)}(x), t_2,\dots ,t_r)$ for each 
 $(t_1,\dots ,t_r) \in [0,1]^r$,  and  let $g_{i+1}$ act  as the identity on all points of
 $\overline{\sigma_r} $ (resp.  $\partial \sigma_{r+1} \setminus \overline{\sigma_r} $) not
 in our collar for $\overline{\sigma_r }$ 
 (resp. $\partial \sigma_{r+1}\setminus \overline{\sigma_r}$) as well as for points that are
 in our collar but sit over points that were not in the collar at an earlier stage.
 \end{proof}

Now to a  condition that will be used to prove that our upcoming collapses preserve regularity: 

\begin{condition}
\label{lub-condition}
{\rm
Let $g$ be an identification map on a regular CW
space $X_K$ 
such that $g$ maps cells onto cells,  maps an open cell $F$ 
onto one of its boundary  cells, and acts homeomorphically on 
$X_K \setminus \overline{F} $.  Then 
$g$ satisfies the {\it least upper bound condition}  (LUB)
if for any pair of open cells $A,B\subseteq \overline{F}$ such that $g(A)=g(B)$
and any face $F'$ that is a least upper bound for $A$ and $B$
in the closure poset just prior to the application of $g$,  $F'$ is
also mapped onto one of its boundary cells by $g$. }
\end{condition}

\begin{remark}
At each collapsing step, 
there will be  
one cell among those getting collapsed that has all other cells 
getting collapsed in the same step in its closure.  (LUB)  then  implies for any pair of cells
$F, F'$ which are both least upper bounds for  cells $\sigma $ and $ \sigma '$ 
just prior to the 
collapse of $F$ such that this collapse identifies  $\sigma $ with $\sigma '$, 
then this step also collapses a larger cell having both $F$ and $F'$ in its closure,
doing so in a way that induces the collapses of  $F$ and $F'$. 
\end{remark}

Next is the main result of this section, a topological construction 
showing how under
certain (mainly combinatorial) conditions a regular CW ball admits a cell collapse 
(in the sense of 
Definition ~\ref{collapse-map-defn}), hence admits an identification 
map preserving homeomorphism type and regularity.  The result is phrased as an
inductive statement so as 
to allow the performance of a series of such collapses  by showing that after
each collapse the conditions are preserved that are needed to apply the theorem
again.  The proof is largely devoted to 
defining explicitly a suitable continuous, surjective map based on a collection of 
parallel-like curves covering the cell to be collapsed.  Figures 
~\ref{collapse-figure} and ~\ref{collar-collapse-figure} provide pictures that may be 
helpful to seeing what these maps are doing.

\begin{theorem}\label{topol-collapse2}
Let $K_0 $ be  a convex polytope. 
Let  $g_1,\dots ,g_i$ be collapsing maps 
with  $g_j: X_{K_{j-1}}\rightarrow X_{K_j} $ for regular CW complexes 
$K_0,\dots ,K_i$ all having underlying space $X$.
Suppose  $K_i$ satisfies Condition ~\ref{manifold-property} (IM)
and that there is an open cell $L_i^0 $ in $\partial K_0$ upon which $g_i\circ\cdots\circ g_1$
acts homeomorphically and a collection $\mathcal{C} = \{ g_i\circ
\cdots \circ  g_1(c)| c\in 
\mathcal{C}_i^0\} $ of parallel-like curves covering $\overline{L_i}$ for 
$L_i = g_i\circ \cdots \circ g_1(L_i^0) \in K_i$. 
Then there is an identification map $g_{i+1}:X_{K_i} \rightarrow X_{K_{i+1}}$ 
specified by $\mathcal{C}$.  If
$g_{i+1}$ also satisfies Condition ~\ref{lub-condition} (LUB), 
then $g_{i+1}$ is a collapsing map and $K_{i+1}$ is a regular CW complex also
satisfying Condition ~\ref{manifold-property} (IM). 
\end{theorem}

\begin{proof}
We may assume $L_i$ is top dimensional in $\partial K_i$, 
because otherwise  we may 
choose a cell 
$L'$ in $K_i$ with $\dim L' =\dim L_i + 1$ and $L_i\subseteq \partial L'$, 
define the collapsing map 
on $\partial L'$ as described below, then use Lemma 
~\ref{bigger-collapse} to extend $g_{i+1}$ to the entire complex.  We will define $g_{i+1}$
on $\partial L'$ in such a way that Lemma ~\ref{extend-interior} will enable its 
extension to $\overline{L'}$.

We will construct a continuous, surjective 
function $g_{i+1}$ that maps entire curves in
$\mathcal{C}$ to points in $G_2$ (cf. Definition ~\ref{parallel-like-def}), 
thereby collapsing $\overline{L_i}$ onto $G_2$, and
in the process identifying each point of $G_1$ with a point of $G_2$; $g_{i+1}$ restricted to
$K_i\setminus \overline{L_i}$ will be a homeomorphism. 
First we define an auxiliary family $\mathcal{C'}$ 
of curves that covers not only $\overline{L_i}$ but also a collar just 
outside its boundary.  We will stretch these curve extensions from the collar to
cover $\overline{L_i}$; the introduction of additional curves  within the collar will enable 
interpolation from the action of $g_{i+1}$ on 
$\overline{ L_i}$ to the identity map  outside this collar.  Now to the details.

 First consider any 
$c\in \mathcal{C}$ with $c\cap \partial (L_i) = \{ x,y\} $ for  points 
$x\in \overline{G_1}$ and $y\in \overline{G_2}$.
 Extend $c$ to include all points $ (y,t)$ and $(x,t)$ for 
 $t\in [0,1]$ to obtain a lengthened curve $c'$.   
 Since $c_1\ne c_2$ for  $c_1,c_2\in\mathcal{C}$ 
 implies that $c_1$ and $c_2$ have distinct 
 endpoints in $G_1$ by Definition ~\ref{parallel-like-def}, part 3, 
 the curve extensions $\{ (x,t)| t\in [0,1]\} $ given by the various
 points $x\in G_1$ are nonoverlapping.  
 It will not matter if distinct $c_1,c_2\in \mathcal{C}$ have the same
 endpoint $y \in G_2$, 
  because $g_{i+1}|_{G_2 \times [0,1]} = ID$. 
   In this situation, let $y\times [0,1]$ be part of both $c_1'$ and
 $c_2'$.   Definition ~\ref{parallel-like-def}, part 2,  guarantees $x\ne y$ for
 each nontrivial $c\in \mathcal{C}$.    For each $c\in \mathcal{C}$  where
 $c$ is a single point in $\partial (L_i)$, extend to
 $c' = \{ (c,t)| t\in [0,1] \} $, and let $g_{i+1}|_{c'} = ID$. 
  For each nontrivial curve
 $c$ with $c \subseteq \partial (L_i)$, we  create a family 
 $F_c$ of curves in $N_i $ (see Figure ~\ref{collar-collapse-figure})
 so that $  F_c$ covers exactly $ \{ c\} \times [0,1]  = 
 \{ (x,t) | x\in c; t\in [0,1] \} $.
  We make  one such curve $c_t\in F_c$ for each $t\in [0,1 ]$, doing this in such a way
  that we have  $c\subseteq c_0$.  
 Letting $c = \{ c(t) | t\in [0,1] \} $ with $c(0) \in G_1$ and $c(1) \in G_2$, then 
  for each $t\in [0,1]$, we define 
  $c_t$ as 
$  \{ (c(t/2),t')| t'\ge t \} \cup  \{ (c(t''),t) | t''\ge t/2 \} \subseteq \{ c \} \times [0,1]$.  Now 
$\mathcal{C'}$ is comprised 
 of the union of these families of curves $F_c$ 
 for each nontrivial curve $c\subseteq \partial L_i$, along 
 with an extended curve $c'$ resulting from each $c\in \mathcal{C}$  which is
 trivial or only intersects
 $\partial (L_i)$ in its endpoints.

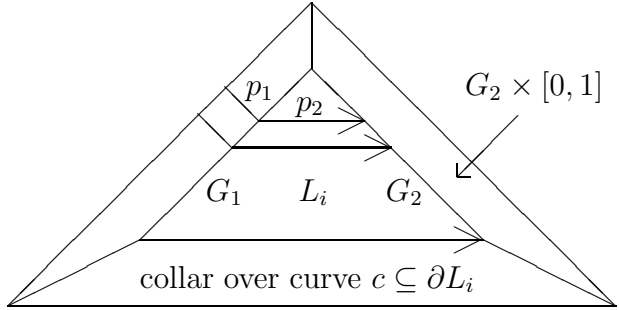
\begin{figure}[h]
\begin{picture}(100,115)(-340,-70)
\put(-360,-45){\line(1,0){130}}
\put(-315,0){\line(1,0){40}}
\put(-325,-10){\line(1,0){60}}
\put(-325,-10){\line(-1,1){13}}
\put(-328,13){\line(1,-1){13}}
\put(-360,-45){\line(1,1){65}}
\put(-230,-45){\line(-1,1){65}}
\put(-265,-10){\line(-2,1){10}}
\put(-265,-10){\line(-2,-1){10}}
\put(-275,0){\line(-2,1){10}}
\put(-275,0){\line(-2,-1){10}}
\put(-231,-45){\line(-2,1){10}}
\put(-231,-45){\line(-2,-1){10}}
\put(-295,20){\line(0,1){25}}
\put(-295,45){\line(1,-1){115}}
\put(-295,45){\line(-1,-1){115}}
\put(-410,-70){\line(1,0){230}}
\put(-301,4){$p_2$} 
\put(-320,10){$p_1$} 
\put(-300,-30){$L_i$}
\put(-335,-30){$G_1$}
\put(-267,-30){$G_2$}
\put(-237,10){$G_2 \times [0,1] $}
\put(-217,2){\line(-1,-1){23}}
\put(-240,-21){\line(0,1){5}}
\put(-240,-21){\line(1,0){5}}
\put(-410,-70){\line(2,1){50}}
\put(-180,-70){\line(-2,1){50}}
\put(-360,-63){${\rm collar}\hspace{.05in}  
{\rm over} 
\hspace{.05in}
{\rm curve } \hspace{.05in} c\subseteq \partial L_i$} 
\end{picture}
\caption{Schematic for  collapsing map} 
\label{collapse-figure}
\end{figure}

Now we define
$g_{i+1}:K_i \rightarrow K_i$  by specifying how it maps
each $c'\in \mathcal{C'}$ surjectively onto itself.  
First consider any 
$c'\in \mathcal{C'}$ obtained by extending some nontrivial 
$c\in \mathcal{C}$.  Represent the points of $c'$ as 
$\{ c(t) | t\in [-1,2 ]  \} $, where $ [-1,0]$ gives $p_1 = c' \cap  (G_1 \times [0,1] )$, i.e. the part
of the collar sitting over the endpoint of $c$ in $G_1$,
whereas $[0,1]$ specifies the points in $p_2 = c' \cap L_i$, and
$[1,2]$ gives the points in $p_3 = c' \cap  (G_2 \times [0,1] )$.
Let $g_{i+1}(c(t)) = c(1)\in \overline{G_2}$ 
for $t\in [0,1]$, i.e. let $g_{i+1}$ map the entire segment $p_2$ in Figure ~\ref{collapse-figure}
to its endpoint in $\overline{G_2}$;  let 
$g_{i+1}(c(t)) = c(2t+1)$ for $t\in [-1,0]$, i.e. stretch the segment $p_1$ in 
Figure ~\ref{collapse-figure} to cover $p_1\cup p_2$; and let
$g_{i+1}(c(t)) = c(t)$ for $t\in [1,2]$.

Next consider any  family $F_c$ of elements of $\mathcal{C'}$ covering $c\times [0,1]$ for
some $c\in \mathcal{C}$ with $c\subseteq \partial(L_i)$.  The map here is designed so as 
to  interpolate 
from the collapsing map needed in $\partial L_i$ to the identity map outside
the collar. 
Points are represented as 
ordered pairs $(c(t_1),t_2)$ for $t_1,t_2\in [0,1]$.  
For each $t\in [0,1]$,
the map $g_{i+1}$ sends $s_2 \cup s_3 = \{ (c(t'),t) | t/2 \le t' \le 1 \} $ to 
$ s_3 = \{ (c(t'),t)| 1-t/2 \le t' \le 1\} $ by appropriate scaling of the parameter, and $g_{i+1}$ 
stretches $s_1 = \{ (c(t/2),t')| t'\ge t/2 \}   $ to cover $s_1\cup s_2$   for 
$s_2 = \{ (c(t'),t) | t/2 \le t' \le 1-t/2 \} $, again  
by reparametrization by a suitable  linear
scaling factor.  See Figure
~\ref{collar-collapse-figure}.  In other words, 
$g_{i+1}$ sends $\{ (c(t/2),t')| \frac{1+t}{2}\le t'\le 1 \} $ to $\{ (c(t/2),t') | t\le t'\le 1\} $
and sends $\{ (c(t/2),t')| t\le t'\le \frac{1+t}{2} \} $ to $\{ (c(t'),t)| t/2\le t' \le 1-t/2 \} $. 

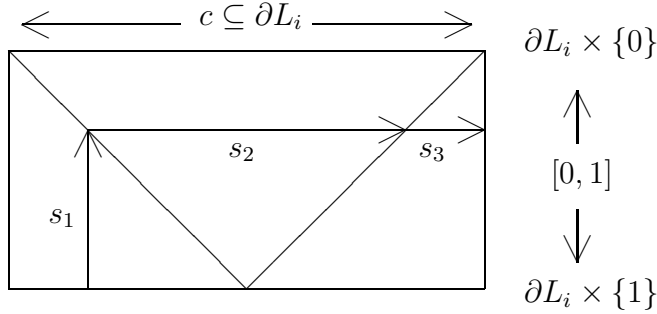
\begin{figure}[h]
\begin{picture}(100,110)(-340,-40)

\put(-385,-40){\line(1,0){180}}

\put(-385,50){\line(1,0){180}}

\put(-385,-40){\line(0,1){90}}
\put(-205,-40){\line(0,1){90}}
\put(-385,50){\line(1,-1){90}}
\put(-205,50){\line(-1,-1){90}}

\put(-355,-40){\line(0,1){60}}
\put(-355,20){\line(1,0){150}}

\put(-355,20){\line(1,-2){5}}
\put(-355,20){\line(-1,-2){5}}

\put(-235,20){\line(-2,1){10}}
\put(-235,20){\line(-2,-1){10}}

\put(-205,20){\line(-2,1){10}}
\put(-205,20){\line(-2,-1){10}}

\put(-302,10){$s_2$}
\put(-370,-15){$s_1$}
\put(-230,10){$s_3$}

\put(-313,60){$c\subseteq \partial L_i$}

\put(-180,0){$[0,1]$}
\put(-170,15){\line(0,1){20}}
\put(-170,-10){\line(0,-1){20}}

\put(-170,35){\line(1,-2){5}}
\put(-170,35){\line(-1,-2){5}}
\put(-170,-30){\line(1,2){5}}
\put(-170,-30){\line(-1,2){5}}

\put(-270,60){\line(1,0){60}}
\put(-320,60){\line(-1,0){60}}

\put(-210,60){\line(-2,1){10}}
\put(-210,60){\line(-2,-1){10}}
\put(-380,60){\line(2,1){10}}
\put(-380,60){\line(2,-1){10}}

\put(-190,50){$\partial L_i \times \{ 0\}  $}
\put(-190,-45){$ \partial L_i \times \{ 1\} $}
\end{picture}
\caption{Defining $g_{i+1}$ on 
collar portion over curve $c\subseteq \partial L_i $} 
\label{collar-collapse-figure}
\end{figure}

Note that $g_{i+1}$ acts as the identity on $\partial(N_i)$
and  acts injectively on  $N_i\setminus \overline{L_i}$. 
We describe next how to 
choose the parametrizations of various curves in $\mathcal{C'}$ in a way that makes 
$g_{i+1}$ continuous.  We use that the relative interiors of
the curves in 
$\mathcal{C}$ have preimages that are  a family of parallel line segments covering a 
convex region (a face of a polytope), enabling us to
choose parametrizations within this simplex 
which induce suitable ones for $\mathcal{C}$; the point is to use Remark ~\ref{length-fn}
to give a continuous function from $L_0^i$ to $[0,1]$ recording for each point how it is 
parametrized within a curve in $\mathcal{C}$.
From this, 
we obtain suitable parametrizations for $\mathcal{C'}$ by virtue of the collar that is also 
covered by $\mathcal{C'}$ 
being   homeomorphic to $\partial (L_i)\times [0,1]$.

Surjectivity and continuity of $g_{i+1} $ imply it induces a  continuous,
bijective function $\overline{g_{i+1}}$ 
from  $K_i/(\ker g_{i+1} )$  to $K_{i+1}$.  Continuity of $(\overline{g_{i+1}})^{-1}$ is then 
immediate, because any bijective, continuous function from a compact set to a Hausdorff
space has continuous inverse.  Thus, 
$g_{i+1}(X_{K_i}) = X_{K_{i+1}}$ is homeomorphic
to $X_{K_i}/(\ker g_{i+1})$ via $\overline{g_{i+1}}$. 
 It is straightforward to see 
that the equivalence relation under $\ker g_{i+1} $ gives a closed subset of
$X_k\times X_k$, implying $X_{K_i}/(\ker g_{i+1})$ is Hausdorff.  

Now let us check that regularity is also preserved under $g_{i+1}$. 
Condition  ~\ref{lub-condition} (LUB) implies  for any cell
$G$ not collapsed by $g_{i+1}\circ \cdots \circ g_1$ that  any 
cells $\sigma_1,\sigma_2 \subseteq \partial (G)$ 
identified by $g_{i+1}$ must have some least upper bound $A \subseteq \partial (G) $ 
 which  is also  collapsed by $g_{i+1}$.
 Thus, our homeomorphism
 $\overline{g_{i+1}}$ will restrict to  $\overline{G}/(\ker g_{i+1})$, 
 enabling us to define the
 attaching map for $G $ as a composition of three maps, 
 first applying  the attaching map for $G$ within $K_i $,  then composing this with $g_{i+1}^{-1}|_{\overline{G}_{final}} $, 
 regarded as a map from $X_{K_i}$ to 
 $X_{K_i}$, then composing with $g_{i+1}: X_{K_i}\rightarrow X_{K_{i+1}}$, the 
 point of the second map being to send points the cells of $K_i$ to the cells of 
 $K_{i+1}$; 
  here we let  $\overline{G}_{final} $ denote the set of cells
 mapped homeomorphically to themselves by $g_{i+1}$ or in other words the cells that 
the various fibers of $g_{i+1}$ are mapped onto.
 Lemmas ~\ref{interp-from-curves} and   ~\ref{extend-interior} 
 show that this may be extended  to yield a characteristic
map for all of 
$\overline{G}$, making it a regular CW complex homeomorphic to a ball.

Lemma ~\ref{preserve-manifold} will 
verify Condition ~\ref{manifold-property} (IM) 
for $K_{i+1}$. 
\end{proof}

\begin{lemma}\label{preserve-manifold}
Collapses 
as in Theorem ~\ref{topol-collapse2}
preserve Condition ~\ref{manifold-property}, i.e. the inductive manifold condition.
\end{lemma}

\begin{proof}

Consider any pair of cells $\tau \subseteq \overline{\sigma }$ in $K_{i+1}$
with $\dim \sigma = \dim \tau + 1$.  We must show that $\overline{\partial \sigma \setminus
\tau }$ is a compact manifold with boundary.   By our definition of collapsing map,
 there must be
cells $\tau_i,\sigma_i$ in $K_i$ with $g_{i+1}$ mapping $\tau_i$ homeomorphically to
$\tau $ and $\sigma_i$ homeomorphically to $\sigma $ with 
$\tau_i \subseteq \overline{\sigma_i}$. The proof of Theorem ~\ref{topol-collapse2} shows
that 
$\overline{g_{i+1}}$ is a homeomorphism from $\overline{\sigma_i }/(\ker g_{i+1})$ to 
$\overline{\sigma }$.  In particular, this implies that
$\overline{g_{i+1}}$ gives a bijection from 
$(\overline{\partial \sigma_i \setminus \tau_i })/(\ker g_{i+1})$
to $\overline{\partial \sigma \setminus \tau }$.  By definition, $\overline{g_{i+1}}$ 
is continuous and 
$\ker ( \overline{g_{i+1}}) = \{ (x,x') | \overline{g_{i+1}}(x)=\overline{g_{i+1}}(x') \} \subseteq K_i
\times K_i$ is closed, 
implying that $\overline{g_{i+1}}((\overline{\partial \sigma_i\setminus \tau_i })/(\ker g_{i+1}))$ 
is  compact and Hausdorff by Proposition 13.8 in  chapter 1 of \cite{Br}. 
Since $\overline{g_{i+1}}$ is a continuous, bijective map from a compact set to a 
Hausdorff space, $\overline{g_{i+1}}$ also has continuous inverse.  
Hence, the property of being a compact manifold with boundary transfers as desired.
\end{proof} 

Our collapsing map as defined in the proof of Theorem ~\ref{topol-collapse2} 
is defined in terms of  curve parametrizations and their 
extensions across collars.  Polytopality of $K_0$ is only used 
to supply such a parametrization function that is continuous.  Thus, 
Theorem ~\ref{topol-collapse2},   Lemma ~\ref{preserve-manifold}, and 
Proposition ~\ref{regular-links} all hold 
in more generality than how they are stated, without requiring any modifications to their 
proofs, yielding: 

\begin{corollary}\label{very-general-topol-collapse2}
In Theorem ~\ref{topol-collapse2}, we 
may replace the  polytope $K_0$ by any regular CW complex satisfying the inductive
manifold condition and replace the parallel line segments covering a face $L_i^0$ 
by any family $\mathcal{C}_0$ of curves covering a closed cell
$L_i^0$ such that (1) all of the curves have one endpoint living in a closed cell 
$\overline{G}_1 \subseteq \partial (L_i^0)$ and the other endpoint in a closed cell 
$\overline{G}_2\subseteq \partial (L_i^0)$, 
(2) these curves in $\mathcal{C}_0$ are 
nonoverlapping except possibly at their 
endpoints  in $\overline{G}_1$,  and (3) there is a continuous function
$p$ from $L_i^0 \setminus \{  x\in L_i^0| \hspace{.03in} $x$ \hspace{.06in} {\rm comprises} \hspace{.06in} {\rm a}
\hspace{.06in} {\rm trivial } \hspace{.06in} {\rm curve } \hspace{.06in} {\rm in} \hspace{.06in} \mathcal{C}_0 \} $  to $[0,1]$ that restricts to a homeomorphism from each nontrivial  curve 
$c\in \mathcal{C}_0 $ to $[0,1]$. 
\end{corollary}

An especially useful special case is the following, which in our main application
will allow us to incorporate a change of  coordinates  homeomorphism $ch$ 
which changes the reduced word with respect to which we work
and hence may change the cell structure outside of the closed cell being collapsed:

\begin{corollary}\label{more-general-topol-collapse2}
Let $K$ be a regular CW complex  with closed cell $L$ covered by a family of curves 
$\mathcal{C}$ and  let $\phi $ be a homeomorphism  
from $K$ to a regular CW complex $K'$ 
with closed cell $L'$ covered by a family $\mathcal{C'}$ of parallel-like curves such that
$\phi $  restricted
to $L$ is a cell-structure preserving homeomorphism to $L'$ mapping each 
curve in $\mathcal{C}$ to a curve in $\mathcal{C'}$.  Then we may transfer the parametrization
function for  $L'$ to one for $L$, 
enabling the collapse of $L$ across the curves in $\mathcal{C}$ via
exactly the collapsing map given in the proof of Theorem ~\ref{topol-collapse2}.
Specifically, it does not matter if $K$ has an entirely different cell structure from $K'$ 
outside of the closed cell $L$. 
\end{corollary}

Next we turn to the  links of the cells.   
Following \cite{FS}, we use (essentially) the notion of link 
in the sense of Whitney stratified spaces
(as defined in \cite{GM})  for 
a series of quotient cell complexes obtained by repeated application of Theorem 
~\ref{topol-collapse2} and its extension in Corollary ~\ref{more-general-topol-collapse2}.   In the following proposition, we use a polyhedral cone $\mathcal{C}$ whose cross-sectional slices are copies of our polytope, 
as well as using the quotient spaces of $\mathcal{C}$ 
under our collapses, denoted $\mathcal{C}/\!\!\sim_k$.  
We assume that for each $\tau \subseteq \overline{\sigma }$ 
we have  a projection map $\pi_{\tau }^{\sigma }$ onto the region $R_{\tau }$ 
of $\mathcal{C}/\!\!\sim_k$ indexed by the 
cell $\tau $.  We require that 
the inverse image under $\pi_{\tau }^{\sigma }$  of each point $p\in R_{\tau }$  
lives in a transversal to the open cell $\tau $ where $\tau $ is  obtained by 
restricting  $R_{\tau }$ to the cross-sectional slice of $\mathcal{C}/\!\!\sim $
containing $p$; moreover, we require that this transversal also gives transversals through
sufficiently nearby cross-sectional slices farther from the origin as well.   
The map $\pi_{\tau }^{\sigma }$ 
is defined to have as its domain the intersection of the following  two sets $S_1$ and 
$S_2$:   (1) $S_1$ is the part of $\mathcal{C}/\!\!\sim_k$ consisting of 
the cross-sectional slice containing $p$ as well as all slices farther from the origin,
while  (2) $S_2$
is  the set of regions of $\mathcal{C}/\!\!\sim_k$ indexed by the
cells in  $\overline{\sigma }$ having $\tau $ in their closure.  

\begin{definition}\label{strat-link-def}
{\rm 
Take a point $p$ in the interior of $\tau $ (in any chosen slice of $\mathcal{C}/\!\!\sim_k$
besides the origin). 
 Define the {\it link} of 
$\tau $ in $\overline{\sigma }$, denoted $lk(\tau ,\sigma )$, 
as the intersection of the set $(\pi_{\tau }^{\sigma })^{-1}(p)$ 
with a sufficiently nearby cross-sectional slice of 
$\mathcal{C}/\!\!\sim_k$  contained in $S_1$ and not containing $p$. 
}
\end{definition}

\begin{proposition}\label{regular-links}
Suppose a cell collapse of $L_i$  which meets the conditions of 
Theorem ~\ref{topol-collapse2} or  Corollary ~\ref{more-general-topol-collapse2} 
also satisfies  
$\dim (G_1) = \dim (G_2) = 
\dim (L_i) - 1$.   Suppose additionally that there are projection maps $\pi_{\tau }^{\sigma }$ 
giving
rise to links as defined above both after the current collapsing step and at all earlier steps.
Also suppose that the preimage under the current collapse of 
each transversal  
is the closure of a disjoint union of such transversals given by the correponding projection maps
on the preimage restricted to $int(L_i)$ and that this closure 
is contained in the disjoint union of the  transversals for the 
projection maps on $\overline{L_i}$.
Then the collapse will  preserve the property that
the link of each cell 
is a closed ball with induced cell decomposition a regular CW decomposition.
\end{proposition}

\begin{proof}
Let $F_{i+1}$ be the maximal cell collapsed at this step, and let 
$g_{i+1}$ be the collapsing map.  The result is 
immediate for cells whose links (prior to the collapse) 
do not intersect $\overline{F}_{i+1}$, since 
$g_{i+1}$ acts homeomorphically  everywhere 
except on $\overline{F}_{i+1}$, leaving such links unchanged.

Letting $\overline{G}_2$ be the closed cell onto which 
$g_{i+1}$ maps 
$\overline{F}_{i+1}$,  we now check the result for the link of the open cell $G_2$.  Let 
$\overline{G}_1$ be the closed cell containing the other endpoints of the parallel-like 
curves across which $\overline{F}_{i+1}$ is collapsed.
The point is to show that  
$g_{i+1}^{-1}(lk(p))$ for $p\in G_2$  is a ball and 
that homeomorphism type is preserved under 
application of $g_{i+1}$.
Our use of parallel-like curves  allows us to 
decompose 
$g_{i+1}^{-1}(lk(p))$  into three pieces, 
namely  its restrictions to  (a) points sitting over $\overline{G}_1$, (b) points sitting over 
$\overline{G}_2$ and (c) points  sitting over $\overline{F}_{i+1}
\setminus (\overline{G}_1 \cup \overline{G}_2)$.    The fact that earlier collapses were 
performed successfully yielding links that were regular CW  balls 
implies that (a) and (c) are each regular CW balls and that (b) is the product of a regular 
CW ball with the open interval  $(0,1)$. The parallel-like curves glue 
these together in a natural way that makes the union also  a regular CW ball, since (b) is 
$(0,1)\times lk_{\sigma /\sim_i} F_{i+1}$ where the desired link is being taken in $\sigma $ 
and $\sim_i$ is the equivalence relation resulting from the first $i$ collapses,
while (a) and (c) are each homeomorphic to $lk_{\sigma /\sim_i }F_{i+1}$.   Moreover,
applying $g_{i+1}$ glues the ball sitting over $\overline{G}_2$ to the ball sitting over
$\overline{G}_1$ by identifying the endpoints of curves sitting over  $\overline{F}_{i+1}
\setminus (\overline{G}_1\cup \overline{G}_2)$, yielding a ball.
This would be more subtle  without our  assumption that 
$\dim (G_1) = \dim (G_2)$. 

The same approach  applies to the link of any open cell contained in 
$\overline{G}_2$ which likewise is not collapsed by $g_{i+1}$ 
but also has some  cell of dimension one higher than it  collapsed onto it by $g_{i+1}$.  
Finally, consider the link of any cell $\sigma $ that is not collapsed by $g_{i+1}$, but is in the 
closure of a cell that is collapsed by $g_{i+1}$; the only remaining such case is for a face
$\sigma $ that is covered by parallel-like curves that are each just a single point.  
In this case, the result follows from the fact 
that our homeomorphism from $X/(\ker g_{i+1}) $ to $X$ naturally restricts to some
neighborhood of any point in $\sigma $ by choosing a neighborhood whose boundary
is a union of points which are fixed by $g_{i+1}$ together with entire curves 
from our family of parallel-like curves across which $F_{i+1}$ is collapsed. 
The regular cell structure is obtained by restriction of cells of $X/(\ker g_{i+1})$  to
the resulting ball.
\end{proof}

\section{Combinatorial $0$-Hecke algebra lemmas}\label{0-hecke-section}
 
 The relations $x_ix_i =  x_i$ will yield a 
0-Hecke algebra variant on the deletion exchange property for Coxeter
groups, namely our upcoming
notion of  ``deletion pair''.  In preparation, we first discuss a poset map from
a Boolean algebra to Bruhat order.

 It is natural (and will be helpful) to
associate a  Coxeter group element $w(x_{i_1}\cdots x_{i_d})$ to any 0-Hecke algebra
expression  $x_{i_1}\cdots x_{i_d}$.  This is done 
by applying braid moves and modified
nil-moves to $x_{i_1}\cdots x_{i_d}$ 
to obtain a new expression $x_{j_1}\cdots x_{j_s}$ such that $s_{j_1} \cdots s_{j_s} $ is  
reduced, then letting $w(x_{i_1}\cdots x_{i_d}) = s_{j_1}\cdots s_{j_s}$.  The fact that this
does not depend on the choice  of braid moves and modified nil-moves will follow
from the geometric description for 
$w(x_{i_1}\cdots x_{i_d})$ given next  
in Proposition ~\ref{geometric-w}.

\begin{proposition}\label{geometric-w}
Lusztig's map $f_{(i_1,\dots ,i_d)}$ sends $R_S = \{ (t_1,\dots ,t_d) \subseteq  
\RR_{\ge 0}^d | t_i > 0 \hspace{.1in} {\rm iff} \hspace{.1in} i\in S \} $
with $S = \{ j_1,\dots ,j_k\} $  to the open cell
$Y_u^o$ for $u = w(x_{i_{j_1}}\cdots x_{i_{j_k}})$.
\end{proposition}

\begin{proof}
This follows from Theorem ~\ref{word-property}, which ensures the existence of a 
series of braid moves and modified nil-moves which may be applied to $x_{i_{j_1}}
\cdots x_{i_{j_k}}$ mapping the points of $R_S$ onto the points of some cell 
$R_T$ indexed by a reduced expression, sending each $x\in R_S$ to some $y\in R_T$ 
with the property that $f_{(i_1,\dots ,i_d)}(x)=f_{(i_1,\dots ,i_d)}(y)$.
\end{proof}

\begin{corollary}\label{indep-w}
The Coxeter group element $w(x_{i_{j_1}}\cdots x_{i_{j_k}})$ does not depend on the 
series of braid moves and modified nil-moves used to convert $x_{i_{j_1}}\cdots x_{i_{j_k}}$
into a reduced expression.
\end{corollary}

\begin{corollary}\label{w-bruhat}
If $A = x_{j_1}\cdots x_{j_r}$ and $B=x_{k_1}\cdots x_{k_s}$ with 
$\{ j_1,\dots ,j_r\} \subseteq
\{ k_1,\dots ,k_s \} $, 
then $w(A) \le_{Bruhat} w(B)$.   Thus,
$w$ is a poset map from a Boolean algebra to Bruhat order.
 \end{corollary}

\begin{proof}
 $A$ is obtained from $B$ by setting some parameters to 0, hence the open cell
to which $A$ maps is in the closure of the open cell to which $B$ maps.  But Bruhat order is
the closure order on cells of $Y_w$.
\end{proof}
See \cite{AH} for additional properties of this poset map $w$ from a Boolean
algebra to Bruhat order.

The following 0-Hecke algebra notion, that of deletion pair,  will play a key role in
various lemmas in Section ~\ref{explicit-section} , i.e. in checking the combinatorial
conditions needed in our proof of the Fomin-Shapiro Conjecture in order to use the
cell collapses developed in Section ~\ref{topol-collapse-section}.

\begin{definition}\label{del-def}
Define a {\it deletion pair} in  a 0-Hecke algebra expression
$x_{i_1}\cdots x_{i_d}$ to be a pair $\{ x_{i_r}, x_{i_s} \} $ such that 
the subexpression $x_{i_r}\cdots x_{i_s}$  is not reduced  but 
$\hat{x}_{i_r}\cdots x_{i_s}$ and $x_{i_r}\cdots \hat{x}_{i_s}$ are each reduced.
\end{definition}

For example, in type A 
the expression $x_1x_2x_1x_2$ has deletion pair $\{ x_{i_1},x_{i_4}\} $.

\begin{lemma}\label{del-pair-lemma}
If $\{ x_{i_r},x_{i_s} \} $ comprise 
a deletion pair, then $w(x_{i_r}\cdots x_{i_s}) = w(\hat{x}_{i_r}
\cdots x_{i_s}) = w(x_{i_r}\cdots \hat{x}_{i_s})$.
\end{lemma}

\begin{proof}
$w(x_{i_r}\cdots \hat{x}_{i_s}) \le  w(x_{i_r}\cdots x_{i_s})$ and
$w(\hat{x}_{i_r}\cdots x_{i_s}) \le  w(x_{i_r}\cdots x_{i_s})$ in Bruhat order, 
while all three of these Coxeter
group elements have the same length, so the equalities follow. 
\end{proof}

See \cite{FG} for a faithful representation 
of the 0-Hecke algebra in which the simple reflections  which do not increase length
act by doing nothing.

\begin{remark}
{\rm 
If $\{ x_{i_u}, x_{i_v} \} $ comprise a deletion pair in a word $x_{i_1}\cdots x_{i_d}$ and we apply 
a braid relation in which $x_{i_u}$ is the farthest letter from $x_{i_v}$
in the segment being braided, then 
the resulting expression will have as a deletion pair $x_{i_v}$ together with the 
nearest letter  to it in the segment that was braided.
We regard this
as  a braided version of the same deletion pair.  
}
\end{remark}

\begin{example}
Applying a braid relation 
to $x_1x_2x_1x_2$ yields $x_2x_1x_2x_2$; we regard the third and fourth letter in the new
expression as a braided version of the  deletion pair  comprised of the first and fourth
letters in the original expression.  
\end{example}

Given  a reduced expression $x_{i_1}\cdots x_{i_d}$, associate a 
Coxeter  group element    $R(x_{i_j})$ to each $x_{i_j}$ by letting 
$R(x_{i_j}) = s_{i_1}\cdots s_{i_{j-1}}s_{i_j}s_{i_{j-1}}\cdots s_{i_1}$.  For finite
Coxeter groups, these  will be the reflections.  
In the case of a 
nonreduced expression, 
if  $w(x_{i_1}\cdots x_{i_j}) = w(x_{i_1}\cdots x_{i_{j-1}})$, then 
we find the largest $j'<j$ such that $s_{i_{j'}}\cdots s_{i_{j-1}} =
s_{i_{j'+1}}\cdots s_{i_j}$ and let $R(x_{i_j}) = R(x_{i_{j'}})$.  

Our original proof of the next lemma relied on the fact that all finite Coxeter groups are 
also reflection groups.
Sergey Fomin provided us with the proof presented below which 
avoids passing to reflection groups.  In particular, this allows  
us to avoid assuming here, and thereby
throughout the paper,  that our Coxeter groups are of finite type.

\begin{lemma}\label{first-root-lemma}
Given a reduced expression $x_{i_1}\cdots x_{i_m}$ in the 0-Hecke algebra of a 
Coxeter group in which 
$R(x_{i_m}) = s_{i_0}$, then 
$x_{i_0}x_{i_1}\cdots x_{i_m}$ has $\{ x_{i_0}, x_{i_m} \} $ as a deletion pair.
\end{lemma}

\begin{proof}
Let us show that $x_{i_0}x_{i_1}\cdots x_{i_{m-1}}$ is reduced whereas 
$x_{i_0}x_{i_1}\cdots x_{i_m}$ is not.  By exercise 8 in [BB, Chapter 1], 
it suffices to prove $R(x_{i_j}) \ne R(x_{i_k})$
for all $0\le j < k \le m-1$ along with proving $R(x_{i_0}) = R(x_{i_m})$.  But
$x_{i_1}\cdots x_{i_{m-1}}$ is reduced, which implies
$s_{i_1}\cdots s_{i_{j-1}}s_{i_j}s_{i_{j-1}}\cdots s_{i_1} \ne s_{i_1}\cdots s_{i_{k-1}}s_{i_k}
s_{i_{k-1}}\cdots s_{i_1}$ for $1\le j < k \le m-1$.  This in turn implies
$R(x_{i_j}) \ne R(x_{i_k})$ for $1\le j < k \le m-1$ with respect to the expression
$x_{i_0}x_{i_1}\cdots x_{i_m}$, since we simply conjugate the 
previous inequalities by $s_{i_0}$ to obtain the desired inequalities.  On the other hand, 
$R(x_{i_0}) = s_{i_0} = s_{i_0}^3 =  s_{i_0}(s_{i_1}\cdots s_{i_m}\cdots s_{i_1})s_{i_0}$,
completing the proof.  
\end{proof}

\section{Proof of the Fomin-Shapiro Conjecture }\label{application-section}

In this section, we apply Theorem ~\ref{sufficient-theorem} to the
stratified spaces $Y_w$
introduced by Fomin and Shapiro in \cite{FS} to prove the following.

\begin{theorem}\label{fs-theorem}
The Bruhat decomposition  $Y_w$ of the link of the identity in  the totally
nonnegative real part of the unipotent radical of a Borel subgroup in a 
semisimple, simply connected
group 
defined and split over $\RR $ is a regular CW complex homeomorphic
to a ball.  Moreover the link of each cell
is as well.
\end{theorem}

Much of the proof will consist of first constructing a regular CW complex $K$ 
that  will be a quotient space 
$(\RR_{\ge 0}^d \cap S_1^{d-1})/\sim $
of a simplex; this is obtained from the simplex 
by a series of collapses  in Section ~\ref{explicit-section}.  Then we prove that the 
quotient space map induced from $f_{(i_1,\dots ,i_d)}$ will act on $K$ in a manner 
that meets all the requirements of Corollary ~\ref{regular-build}.
As preparation, we 
first define a much simpler equivalence relation in Section
~\ref{explicit-special-section}, denoted $\sim_C$, 
doing identifications based only on stuttering pairs  which may be obtained in 
nonreduced expressions exclusively by  applying commutation moves. Then we give
the more difficult analogous argument for $\sim $.  
Finally, Section
~\ref{final-section}  applies Theorem ~\ref{sufficient-theorem} to show
  that the induced map $\overline{f_{(i_1,\dots ,i_d)}}$ from 
$K$ to $Y_w$ is a homeomorphism. 

Let us now establish some convenient notation for the proof.  Let
$\RR_{\ge 0}^d \cap S_1^{d-1}$ denote the restriction of $\RR_{\ge 0}^d$ to the
 hyperplane with coordinates summing to 1.  
We will make extensive use of the fact that
this is a convex polytope.  Define  the {\it regions} or {\it faces} 
in $\RR_{\ge 0}^d \cap S^{d-1}_1$ 
as the sets $R_S =  \{ (t_1,\dots ,t_d) \in \RR_{\ge 0}^d \cap S_1^{d-1} |
t_i > 0 \text{ iff } i\in S \} $.
We associate the  $0$-Hecke algebra
expression $x_{i_{j_1}}\cdots x_{i_{j_k}} $ to the 
region $R_{\{ j_1,\dots ,j_k \} }$,  calling this the {\it $x$-expression} of 
$R_{ \{ j_1,\dots ,j_k \} }$, denoted $x(R_{ \{ j_1,\dots ,j_k \} })$.
  For  subexpressions $x(A)$ and $x(B)$ of 
$x_{i_1}\cdots x_{i_d}$, let $x(A)\vee x(B)$ be
the expression made of the union of the indexing positions from $x(A)$ and $x(B)$, and 
let $A\vee B$ denote the cell given by that expression.
To keep track of the positions of the 
nonzero parameters, we sometimes  also include 1's  as placeholders, 
e.g. describing the region $R_{\{ 1,3\} } $ 
given by the map $f_{(1,2,1)}$ by the expression $x_1\cdot 1  \cdot x_1$.  
We say that an $x$-expression is {\it stuttering } if it directly admits a  modified nil-move $x_i x_i \rightarrow x_i$, i.e. $x_i(u)x_i(v) = x_i(u+v)$  (cf. Section ~\ref{0-hecke-bg-section}).
An expression is {\it commutation
equivalent to a stuttering expression} if  it  admits a series of commutation moves 
yielding a stuttering expression.

In Section  ~\ref{explicit-section}, we will define 
the equivalence  relation $\sim $ on $ \RR_{\ge 0}^d \cap S_1^{d-1}$ by repeated
application of the following idea:
if the $x$-expression $x_{i_{j_1}}\cdots x_{i_{j_k}}$
associated to a point $(t_1,\dots ,t_d)\in \RR_{\ge 0}^d$ is not
reduced, then we may apply commutation moves and  long
braid moves to it, causing a coordinate change to new
coordinates $(u_1,\dots ,u_d)$ in which we may apply  a substitution 
$x_i(u_{s-1})x_i(u_s) = x_i (u_{s-1} + u_s)$. 
Each region $R_{ \{ j_1,\dots ,j_k \} }$  indexed by a non-reduced word
is collapsed by such a move, at which point  
we  say  $(u_1,\dots ,u_d) \sim (u_1',\dots ,u_d')$ for those points 
$(u_1',\dots ,u_d') \in \RR_{\ge 0}^d \cap S_1^{d-1}$ 
such that $u_{s-1}' + u_s' = u_{s-1} + u_s$ and $u_i'=u_i$ for  $i\ne s-1,s$. 
It is important to our collapsing argument that for 
each non-reduced subword $(i_{j_1},\dots ,i_{j_k})$ of $(i_1,\dots ,i_d)$, 
we choose exactly one such way of identifying points of the  open
cell $R_{ \{ j_1,\dots ,j_k \} }$
with points having strictly fewer nonzero parameters, namely the identifications 
dictated by the collapse we choose to apply to  $R_{\{ j_1,\dots ,j_k \} }$.  
Additional identifications will hold by transitivity of $\sim $.

Before turning to the details, 
let us briefly enumerate the 
main upcoming definitions, lemmas and theorems and how they fit together. 
Definition ~\ref{curve-def}  introduces 
the parallel-like curves  that will be used to induce the collapses leading to  $\sim_C$. 
Then we  verify the distinct endpoints condition (DE) in Lemma
~\ref{end-lemma},  the least upper bound condition (LUB) in
Lemma  ~\ref{inject-lemma},  the distinct initial points condition (DIP)  in 
Lemma ~\ref{G_1-ends}, 
and deduce  from all this the regularity of the 
quotient cell complex $(\RR_{\ge 0}^d \cap S^{d-1}_1)/\!\!\sim_C $ in Theorem
~\ref{homeom-lemma}.  Afterwards, we characterize exactly which faces are
identified with each other
by $\sim_C$ in Lemma ~\ref{commutation-equiv}.

Now in the general case  of $\sim $, 
we prove a similar series of lemmas, after first showing that
long braid moves may be accomplished by change of coordinates maps that are 
homeomorphisms on the closed cells to be collapsed.   The requisite 
parallel-like curves are specified in  Definition  ~\ref{general-level-def}.
The result about changes of coordinates 
is obtained through 
Lemmas ~\ref{single-braid}, ~\ref{induct-on-d}, 
~\ref{ioa},  and ~\ref{collapse-order}.
Next we verify the conditions (DIP) 
in Lemma ~\ref{G_2-ends}, (DE) 
in Lemma ~\ref{distinct},  (LUB) 
 in Lemma ~\ref{stay-regular}, and  we show the requisite equidimensionality 
 to deduce regularity of  links in Lemma 
~\ref{same-dim-lemma}.  In proving these results for a particular collapse, we assume
inductively that all earlier collapses were already performed successfully, and we also 
assume inductively that all results in the paper hold for all $d' < d$ to prove the results 
for reduced word  $(i_1,\dots ,i_d)$ of length $d$.
We then combine these ingredients to deduce homeomorphism
type and regularity of  $(\RR_{\ge 0}^d \cap S^{d-1}_1)/\!\!\sim $ in 
Theorem  ~\ref{homeom2-lemma}.

Finally, we prove that the induced map $\overline{f}_{(i_1,\dots ,i_d)}$ from this
quotient space to $Y_w$  
 is a homeomorphism that preserves cell structure, 
implying that $Y_w$ is a regular CW complex homeomorphic to a
ball. 
To this end, 
Lemma ~\ref{injective-condition} uses the exchange axiom for Coxeter groups in order
to verify the attaching map injectivity requirement of Theorem ~\ref{sufficient-theorem}, 
allowing the proof of the Fomin-Shapiro Conjecture to 
 be completed in Theorem ~\ref{main-theorem}.

\subsection{Collapsing a simplex to obtain
$(\RR_{\ge 0}^d \cap S^{d-1}_1)/\!\!\sim_C $}
\label{explicit-special-section}

In this section, we collapse those faces of $\RR_{\ge 0}^d \cap S^{d-1}_1$
whose words are commutation equivalent to 
stuttering words,   denoting  the resulting 
identifications by $\sim_C$.  We prove 
that $(\RR_{\ge 0}^d \cap S_1^{d-1})/\!\!\sim_C$ is a regular CW complex homeomorphic
to a ball by proving that regularity and homeomorphism type are preserved under 
each in a series of collapses of the type introduced in Section ~\ref{topol-collapse-section}.
 A separate proof  for $\sim_C$
is given before turning to the general case 
for two reasons: (1) it illustrates the general strategy in a much simpler setting, and  
(2) this result will be used in the proofs of
 Lemmas ~\ref{single-braid} and ~\ref{collapse-order}, two 
key ingredients to the general case.

\begin{definition}
{\rm 
An {\it omittable pair} of an $x$-expression  $x(F)$ is a 
pair  $\{ x_{i_l},x_{i_r} \}  $ of letters appearing in $x(F)$ with $i_l = i_r$  
such that there exists a series of commutation moves applicable to 
$x(F)$ placing the letters into neighboring positions, thereby
enabling a modified nil-move.   
}
\end{definition}

\begin{example}
The $x$-expression
$x_1x_3x_4x_3x_1$ in type A has omittable pair $\{ x_{i_1},x_{i_5} \} $, sometimes denoted
simply by the positions, i.e. as  $\{ 1,5 \} $.
\end{example}

Order as follows  all possible 
triples $(x(F), \{ i_l, i_r \} )$, where 
$F$ is a face in $\RR_{\ge 0}^d \cap S_1^{d-1}$ 
and  $\{ x_{i_l},x_{i_r} \} $ is an omittable pair in 
$x(F)$ with $l<r$.  Use linear order on the  index $r$, then break ties with
linear order on  $r-l$, breaking further ties by reverse linear order
on  $\dim F$,  and breaking any remaining ties arbitrarily.

\begin{example}
$(x_1x_3x_4x_3x_1,\{ 1,5\} )$ precedes $(x_1x_3\cdot x_3 x_1, \{ 1,5\} )$ in this ordering, 
while  
$(x_1x_3\cdot x_3x_1,
\{ 2,4\} )$ comes earlier than both of these.
\end{example}

We obtain from this our sequence of face 
collapses by repeatedly choosing for  the next collapse 
the earliest  triple $(x(F), \{ x_{i_l},x_{i_r}\} )$ whose face $F$ has not yet been collapsed.
Denote by 
$(x(F_m), \{ x_{i_{l_m}}, x_{i_{r_m}} \} )$  the triple chosen for the $m$-th collapse, and let   
$g_m$ be the collapsing map accomplishing this, based on the level curves from 
Definition ~\ref{curve-def}.  Our main task in this section will be to prove that these are 
parallel-like curves  and that $g_m$ meets the requirements of Theorem
~\ref{topol-collapse2}. 

\begin{remark}
 It often will happen that the step collapsing a cell
$F_m$ will also collapse some additional cells. 
However, each
collapsing step will have one cell among those being collapsed at that step
such that all others being collapsed at that step are  in its closure.
\end{remark}

 \begin{definition}\label{curve-def}
 {\rm 
 Given the triple $(x(F_m), \{ x_{i_{l_m}},i_{i_{r_m}} \} )$ specifying the $m$-th collapse, call
the collections  of points 
 $$  \{ (t_1,\dots ,t_d) | t_{l_m} + t_{r_m} = k
 \hspace{.07in} {\rm and } \hspace{.07in}
 t_j = k_j \hspace{.07in} {\rm for} \hspace{.07in} j\not\in \{ l_m,r_m \} \} $$
 in $\overline{F}_m $ 
 for the various sets of constants $\{ k,k_1,\dots ,\hat{k}_{l_m},\dots \hat{k}_{r_m},
 \dots , k_d \} \in [0,1]$ summing to 1 the 
 {\it level curves} of $\overline{F_m}$.
 }
 \end{definition}

\begin{notation}
{\rm 
If cells $G$ and $G'$ are identified during one of the first
 $m-1$ collapsing steps, denote this by $G\sim_m G'$.
}
\end{notation}

\begin{remark}
The collapse given by $(x(F_m), \{ x_{i_{l_m}},x_{i_{r_m}} \} )$ will also collapse those 
cells in $\overline{F}_m $ given by subexpressions of $x(F_m)$ 
having both  $t_{i_{l_m}}>0$ and $t_{i_{r_m}}>0$.
 In this manner, the collapse will identify faces having $t_{i_{l_m}}=0$
with ones instead having $t_{i_{r_m}}=0$.
\end{remark}

To keep track combinatorially of which faces are identified by the collapses giving rise to 
$\sim_C$,
define a {\it slide-move}, or simply a {\it slide}, to be the replacement of $S = \{ j_1,\dots ,j_s \} $ 
by $S' = \{ k_1,\dots ,k_s \}$ for $j_1<\cdots < j_s$ and $k_1<\cdots < k_s$ with $j_i=k_i$ for
$i\ne r$ for some fixed $r$ and $i_{j_r} = i_{k_r}$.  An example in type A  for 
$(i_1,\dots ,i_d) = (1,2,3,1,2)$  is $S = \{ 1,5\} $ and $S' = \{ 4,5\} $. 
An {\it exchange } is the replacement of one letter by another letter than can be 
accomplished by a series of slide-moves and commutation moves.

Now we use combinatorics 
to verify that the hypotheses needed for topological collapses introduced in 
Section ~\ref{topol-collapse-section} are indeed met.  
 Condition 1 of Definition ~\ref{parallel-like-def} follows immediately from our
 set-up.   The next two lemmas check conditions 2 and 3, respectively, by checking the 
 distinct endpoints condition (DE) and  distinct initial points condition (DIP).
 
\begin{lemma}\label{end-lemma}
The collapses inducing $\sim_C $ satisfy Condition ~\ref{distinct-endpoint-condition} (DE).
\end{lemma}

\begin{proof}
What we must prove is that the two endpoints of any nontrivial level curve 
across which a cell $F_i$ is collapsed
live in distinct cells just prior to the collapse.
Suppose $G_1\subseteq \overline{F_i}$ 
with $t_{l_i}>0$ and $t_{r_i}=0$ had been identified already with the face 
$G_2\subseteq \overline{F_i} $ instead having
$t_{r_i}>0$ and $t_{l_i}=0$.  This would have required a series of earlier
slides, including one of the form $r_i\rightarrow r$ for some $r<r_i$.  Our 
collapsing order implies $r_i-r < r_i - l_i$.
By our dimension maximizing assumption in our collapsing order, 
the last of these slide moves shifting the right element of our deletion pair to the left
would have also 
collapsed $F_i$, by virtue of identifying it with a face already collapsed, a contradiction.
\end{proof}

\begin{lemma}\label{G_1-ends}
Suppose a cell $H_1 \subseteq G_1 $ is collapsed prior to the collapse of $F_j$, 
where $H_1$ is identified with $H_2\subseteq G_2$ 
in the collapsing step given by $(F_j,\{ x_{i_r},x_{i_s} \} )$ 
by an exchange of $x_{i_r}\in H_1 \subseteq G_1 $ for $x_{i_s}\in H_2
\subseteq G_2$ with $r<s$.  Then  
the least upper bound $H_1\vee H_2$ will have also been  
collapsed prior to the collapse of $F_j$, and in such a way that any two level
curves in $H_1\vee H_2$ having the same endpoint in $H_1$ are identified
prior to the collapse $(F_j,\{ x_{i_r},x_{i_s} \} )$  with coinciding
parametrizations. 
\end{lemma}

\begin{proof}
The fact that $H_1$ is collapsed before $F_j$ means that 
$x(H_1)|_{x_{i_1}\cdots x_{i_s}}$ contains an omittable pair.  However,  
$x(H_1)|_{x_{i_1}\cdots x_{i_s}} 
= x(H_1)|_{x_{i_1}\cdots x_{i_{s-1}}}$,
implying $x(H_1)|_{x_{i_1}\cdots x_{i_{s-1}}}$ contains an omittable pair based upon which
$H_1$ is collapsed.   By our prioritization of higher dimensional faces in our collapsing
order, the face $H_1\vee \{x_{i_s} \} = H_1\vee H_2$
will have been collapsed in the same way, yielding the desired result. 
\end{proof}
 
 Next we verify (LUB), after 
 first giving an example  showing the idea. 
  \begin{example}
 {\rm 
 The cell  $F$  with $x(F) = x_1x_1x_1$ is collapsed 
based on the deletion pair comprised of its leftmost two $x_1$'s, identifying $x_1\cdot 1\cdot x_1$
with $1\cdot x_1 \cdot x_1$ and in the process also identifying $x_1\cdot 1\cdot 1$ with 
$1\cdot x_1\cdot 1$.   The region with
expression $1\cdot x_1\cdot x_1$ is collapsed  
based on its pair of $x_1$'s,  identifying 
$1\cdot x_1\cdot 1$ with $1\cdot 1\cdot x_1$.  Composing face identifications 
based on these two steps causes 
$x_1\cdot 1 \cdot 1 $ to be identified with
$1\cdot 1 \cdot x_1$, potentially causing the attaching map for 
the face given by 
$x_1\cdot 1 \cdot x_1$ no longer to be injective; however, this face will itself have
been collapsed by this time, by virtue of 
having already been  identified 
with the face $1\cdot x_i \cdot x_i$ which was already collapsed.
}
\end{example}

 \begin{lemma}\label{inject-lemma}
The collapses inducing $\sim_C$ satisfy condition ~\ref{lub-condition} (LUB).
\end{lemma}

\begin{proof}
Suppose that $G_1$ and $G_2$ are identified during
the collapse of  $F$ via deletion pair $\{ x_{i_{l_j}},x_{i_{r_j}}\} $, for 
$x_{i_{l_j}}\in x(G_1)$, $x_{i_{r_j}}\in x(G_2)$, 
and  $F'$ is any face that is a least upper bound for $G_1$ and $G_2$ in the 
closure poset  just prior to the collapse of  $F$.  We must show  that 
$F'$ is  also collapsed by the end of the step collapsing $F$ or is identified with $F$ 
prior to this collapse.
By virtue of our set-up,  $x(\overline{F'} )$ must have subexpressions 
$x(G_1')$ and $x(G_2')$ with $G_1'\sim_j G_1$ and $G_2'\sim_j G_2$.
Consider $x(F')$ and its earliest subexpressions (in our collapsing ordering on
triples) which are $x$-expressions for some  such
$G_1' \sim_j G_1$ and $ G_2'\sim_j G_2$.

Suppose $x_{i_{r_j}} \not\in x(F') $.  This 
implies that $x_{i_{r_j}}$ must
have been exchanged with a letter $x_{i_r}$ to its left during an earlier identification step.
Then $x(F')$ will  have an omittable pair $\{ x_{i_l},x_{i_r} \} $  for some $r< r_j$,
causing $F'$ to have been collapsed strictly 
before the collapse of $F$. 

Now consider the case  $x_{i_{r_j}} \in x(F')$.  
$F'$ will again be collapsed during or prior to the collapse of $F$ unless $x_{i_{r_j}}$ is
the right letter in the leftmost deletion pair of $F'$ and $x_{i_{l_j}}$ has been 
exchanged for a letter $x_{i_l}\in x(F')$ 
to its left appearing instead in $F'$.  
But then by the fact that our collapsing order maximizes dimension among faces with the same
omittable pair, 
this exchange $x_{i_{l_j}}\rightarrow x_{i_l}$  (or each such exchange in a series)
would extend to a face including $x_{i_{r_j}}$, thereby identifying faces including
$x_{i_l}$ and $x_{i_{r_j}}$ with ones instead including $x_{i_{l_j}}$ and $x_{i_{r_j}}$.  
In particular, this identifies $F'$ with $F$ prior to the collapse of $F$.
\end{proof}

Combining the above results will   yield:

\begin{theorem}\label{homeom-lemma}
$(\RR_{\ge 0}^d \cap S_1^{d-1} )/\!\!\sim_C $ is a 
regular CW complex homeomorphic to a ball. 
\end{theorem}
 
 \begin{proof}
 We will use Theorem ~\ref{topol-collapse2}  to 
 prove that each  collapse on $\RR_{\ge 0}^d \cap S_1^{d-1}$  
 may be accomplished in turn by a map that
 preserves homeomorphism type and regularity, provided that all earlier 
 collapses also met the requirements of Theorem ~\ref{topol-collapse2}.  This will 
 imply that the series of collapses producing $(\RR_{\ge 0}^d \cap S_1^{d-1})/\sim_C$ 
 yields a regular CW complex homeomorphic to a ball.

  The parallel-like curves that we will use
 for the $(i+1)$-st collapsing step will be of the form 
 given in Definition ~\ref{curve-def}; they are the  images
 under $g_i\circ\cdots g_1$ of parallel line segments covering a closed cell of 
 $\RR_{\ge 0}^d\cap S_1^{d-1}$.
 To see that
  $g_i\circ \cdots\circ g_1 $ acts on each
 level curve either homeomorphically or by sending it to a point, notice that by 
 definition the interior of any
 nontrivial level curve lives entirely in some open cell $F\subseteq \overline{F_i}$,
 hence a cell upon which all earlier collapses act homeomorphically.
   Lemmas ~\ref{end-lemma}, ~\ref{G_1-ends}, and ~\ref{inject-lemma}
 confirm the distinct endpoints condition (DE), distinct initial points condition (DIP), 
 and least 
 upper bound condition (LUB), 
 respectively, i.e. the requirements of Theorem ~\ref{homeom-lemma}.
  \end{proof}

\begin{proposition}\label{commutation-equiv}
Suppose $x(R_S)$ and $x(R_T)$ are not commutation equivalent to stuttering expressions.
Then $R_S \sim_C R_T$ iff $S$ and $T$ differ 
from each other by a series of commutation moves and slide moves.
\end{proposition}

\begin{proof}
Let $S = \{ j_1,\dots ,j_s \}$ and $T = \{ k_1,\dots k_s \}$.
We begin with pairs of words $x(R_S), x(R_T)$ differing by a single slide, so 
$S\cap T = S\setminus \{ j_r\} = T\setminus \{k_r\} $ for some $r$ with $i_{j_r} = i_{k_r}$.
But then $x(R_{S\cup T})$ is stuttering, implying  $R_{S\cup T}$ was collapsed
by $\sim_C$.  The fact that $x(R_S),x(R_T)$ are not commutation equivalent to 
stuttering expressions implies
$R_{S\cup T}$ could have only been collapsed by 
identifying $R_S$ with $R_T$.
 By transitivity of $\sim_C$, $S$ and $T$ differing by a series of slide
moves  likewise give rise to $R_S, R_T$ with $R_S \sim_C R_T$.
Applying commutation moves to $x(R_S)$ as well to produce $\sigma (x(R_{\sigma (S)})$ 
which is slide equivalent to $x(R_T)$  ensures  $x(R_{S\cup T})$ also admits the 
same commutation moves
leading to a stuttering word, and again $x(R_{S\cup T})$ does not admit any other stuttering
pairs, so again $R_S \sim_C R_T$.
\end{proof}

\begin{example}
(a)
For $(i_1,\dots ,i_d) = (1,2,1)$ in type A,
$R_{\{ 1\} } $ 
is identified with $R_{\{ 3\} }$ 
during the collapse of 
$R_{\{ 1,3\} } $. 
(b) If $(i_1,\dots ,i_d) = (1,3,1)$, then $R_{\{ 1,2 \} }  $
 is identified with $R_{ \{ 2,3\} } $
 during the collapse of $R_{ \{ 1,2,3\} } $.
\end{example}

\subsection{Collapsing a simplex to obtain
$(\RR_{\ge 0}^d \cap S^{d-1}_1)/\!\!\sim $}
\label{explicit-section}

Now we turn to the identifications $\sim $ induced by a series of collapses 
which  collapse 
all faces whose words are nonreduced, i.e. starting afresh so as now to 
incorporate those collapses
requiring long braid moves. 

\begin{definition} 
Given a deletion pair $\{ x_{i_r},x_{i_s}\} $ with $r<s$  in $x(F)$, let 
$c(\{ x_{i_r},x_{i_s}\} ; x(F))$ be
the smallest number of long braid moves needed in a series of braid 
moves applied to $x_{i_r}\cdots x_{i_{s-1}}$ yielding an expression whose last letter
comprises a stutter with $x_{i_s}$. 
\end{definition}

  Lemma ~\ref{del-pair-lemma} combined 
with Theorem ~\ref{word-property}
guarantees existence and finiteness of $c(\{ x_{i_r},x_{i_s} \} ; x(F) )$.

 \begin{example}\label{long-move-example}
 In type A,  we have $c(\{ x_{i_1},x_{i_4} \} ;x_1x_2x_1x_2) = 1$, because we may 
 apply the relation $x_1x_2x_1 \rightarrow x_2x_1x_2$ to obtain $x_2x_1x_2x_2$.
\end{example}
 
 Now let us  order  triples $(x(F), \{ x_{i_l}, x_{i_r} \} )$ where $\{ x_{i_l},x_{i_r}\} $ is a 
deletion pair of $x(F)$ in preparation for our choice of  
a collapsing order on non-reduced  faces. 
By convention, say $l<r$.
Letting  the statistics listed earliest take highest priority, with later statistics used to 
break ties,  order the triples $(x(F),\{ x_{i_l},x_{i_r} \} )$ by: 
(1) linear order on 
$r$, 
(2) linear order on $r-l$,
(3) linear order on $c(\{ x_{i_{l}},x_{i_{r}} \} ; x(F))$, 
and (4) reverse linear order on $\dim F$. 
We may break any remaining ties arbitrarily.

We repeatedly choose the earliest triple
$(x(F),\{ x_{i_l},x_{i_r}\} )$ among those for faces $F$ not yet collapsed.   
Denote the $k$-th such triple chosen by 
$(x(F_k),\{ x_{i_{l_k}},x_{i_{r_k}} \} )$.  We will use 
Theorem ~\ref{topol-collapse2} to 
accomplish the collapse of $F_k$ by a collapsing map $g_k$ applied to 
$(\RR^d_{\ge 0} \cap S_1^{d-1})/\sim_k$, letting 
$\sim_k $ be the equivalence relation 
comprised of the identifications 
that result from  the first $k-1$ collapsing steps, described shortly.
First we will need some results regarding change of coordinate maps.

  Denote by $\sim^s $ the set of all possible identifications $(t_1,\dots ,t_d) \sim (t_1', \dots ,
  t_d')$ under $\sim $ where  $f_{(i_1,\dots ,i_s)}(t_1,\dots ,t_s) = f_{(i_1,\dots ,i_s)}(t_1',
  \dots ,t_s')$ and $t_j = t_j'$ for all $j>s$.  That is, $\sim^s$ consists of all possible 
  identifications based on the leftmost $s$ letters.
 By our inductive hypothesis based on length, we will be able to assume 
 $\sim^s$  is exactly the identifications accomplished by collapses based on 
 deletion pairs  involving only the
 leftmost $s$ letters.  The way our collapsing order was chosen implies then that 
  for each $s$ there will be some $k$ such that 
  $\sim^s = \sim_k $, meaning that $\sim^s$ is the equivalence relation resulting from the first
  $k-1$ collapses.
    For this same pair $k$ and $s$, let us also introduce the notation
  $g^s$ for the composition of maps $g^s =  g_{k-1}\circ \cdots \circ g_1$.
  Let us also establish the notation $g^s_{(i_1,\dots ,i_d)}$ as the map $g^s$ given by
  the reduced word $(i_1,\dots ,i_d)$.  
  Let  $\sim_{(i_1,\dots ,i_d)}$ denote the 
  equivalence relation given by reduced 
  word $(i_1,\dots ,i_d)$ after all possible collapses, and   
  let $\sim^s_{(i_1,\dots ,i_d)}$  denote 
  the equivalence relation $\sim^s $ consisting of identifications based on the leftmost 
  $s$ letters again with respect to initial choice of reduced word 
  $(i_1,\dots ,i_d)$.

\begin{definition}\label{general-level-def}
{\rm 
Given the triple $(x(F_k), \{ x_{i_{l_k}},x_{i_{r_k}} \} )$, 
choose a sequence of 
braid moves on $(i_{l_k},\dots ,i_{r_k-1})$ 
using exactly  $c(x(F), \{ x_{i_{l_k}},x_{i_{r_k}}\} )$ long braid
moves to transform $(i_1,\dots ,i_d)$ into $(j_1,\dots ,j_d)$ with a stutter $j_{r_k-1} = j_{r_k}$.
Obtain new coordinates $(u_1,\dots ,u_d)$ on $\overline{F_k}/\sim_k$ as the  unique 
solution (up to equivalence relation $\sim_k$ )  to 
$$f_{(i_1,\dots ,i_d)}(t_1,\dots ,t_d) = f_{(j_1,\dots ,j_d)}(u_1,\dots ,u_d)$$ which has 
$u_i = t_i $ for $i\not\in \{ l_k,\dots ,r_k -1\} $, as justified by
Lemmas  ~\ref{single-braid} and  ~\ref{collapse-order}.  The {\it level
curves} for this triple  are the collections of points 
$$ \{ (u_1,\dots ,u_d) | u_{r_k-1}+ u_{r_k} = c\hspace{.05in} {\rm and} \hspace{.05in}  u_{m}=c_m \hspace{.05in} {\rm for}\hspace{.05in}{\rm all}\hspace{.05in}  m \not \in \{ r_k-1,r_k \}  \} $$ for the various choices of  
constants $c,c_1,\dots ,\hat{c}_{r_k-1},\hat{c}_{r_k},\dots ,c_d \ge 0$. 
}
\end{definition}

We will eventually prove that these level-curves are parallel-like, in the sense of 
Definition ~\ref{parallel-like-def}.

\begin{example}
{\rm 
Applying  
braid moves  to $x_1x_2x_1x_3x_2x_3 $ yields the expression
$x_2x_1x_3x_2(x_3x_3)$.
Collapsing based on the resulting stuttering pair will cause the proper face
$1\cdot x_2x_1x_3x_2x_3$ to be identified with the face $x_1x_2x_1x_3x_2\cdot 1$.
The proper face $x_1x_2\cdot 1 \cdot 1\cdot 1 \cdot x_3$ is neither collapsed nor identified
with another face in the process since the first and last letters do not form a deletion pair.
On th e
other hand, the face $x_1x_2\cdot 1 \cdot x_3 x_2 x_3 $ would have already been collapsed
at an earlier step, hence need not be considered in the next lemma as part of the 
boundary of the cell indexed by $x_1x_2x_1x_3x_2x_3$.
}
\end{example}

\begin{lemma}\label{single-braid}
Consider the reduced expression $s_is_j\dots $ of length $m(i,j)$ comprised of 
alternating $s_i$'s and $s_j$'s.  
Then the resulting regular CW complex 
$\Delta = (\RR_{\ge 0}^{m(i,j)} \cap S_1^{m(i,j)-1})/\!\!\sim_C $ given by $(i,j,\dots )$
is homeomorphic via the map $f_{(j,i,\dots )}^{-1}\circ f_{(i,j,\dots )}$
 to the regular CW complex
$\Delta' = (\RR_{\ge 0}^{m(i,j)} \cap S_1^{m(i,j)-1})/\!\!\sim_{C'}$ given by $(j,i,\dots )$.
\end{lemma}

\begin{proof}
We may use the fact that $f_{(i,j,\dots )}$ and $f_{(j,i,\dots )}$
act homeomorphically on the interior of the big cell for $\Delta $ and $\Delta '$, 
respectively.
Each point  $x \in \Delta $ not in the interior of the big cell must instead 
belong to a region $R_{\{ i_{j_1},\dots ,i_{j_k} \} }$ 
whose associated Coxeter group element $w(x_{i_{j_1}}\cdots x_{i_{j_k}} )$
has a unique reduced expression, namely one with  the appropriate alternation of
$s_i$'s and $s_j$'s.
Thus,  $x$ must be sent
to a point in $\Delta '$ having this same reduced expression, so that by Proposition
~\ref{commutation-equiv}
the only choices to be made are equivalent to each other under $\sim_{C'}$. 
This map from $\Delta $ to $\Delta'$ 
is therefore a composition of 
two homeomorphisms, namely $f_{(i,j,\dots )}$ and $f_{(j,i,\dots )}^{-1}$, and hence is
itself a homeomorphism.  
\end{proof}

\begin{example}
The type A relation $s_is_{i+1}s_i = s_{i+1}s_is_{i+1}$ 
gives rise to the map $(t_1,t_2,t_3) \rightarrow (t_1',t_2',t_3')$ for 
$(t_1',t_2',t_3') = 
(\frac{t_2t_3}{t_1+t_3},t_1+t_3,\frac{t_1t_2}{t_1+t_3})$ on the interior of
$\{ f_{(1,2,1)}(t_1,t_2,t_3) | t_1,t_2,t_3\ge 0 \} $.  
The above proposition shows that
this map extends to the boundary, i.e. extends to a map
from $(\RR_{\ge 0}^d \cap S_1^{d-1})/\sim_C$
to $(\RR_{\ge 0}^d \cap S_1^{d-1})/\sim_{C'}$.  For instance, it sends
$(t_1,t_2,0)$
to $(0,t_1,t_2)$ for $t_1,t_2>0$, and it 
sends $(0,t_2,0)$  
to the  $\sim_{C'}$-equivalence class 
 $\{ (t_1',0,t_3') |  
t_1'+t_3' = t_2 \} $. 
\end{example}

\begin{lemma}\label{induct-on-d}
 Given a reduced word $(i_1,\dots ,i_d)$ and some $d'<d$ such that 
  the series of collapses for 
 $(i_1,\dots ,i_{d'})$ 
 successfully applies to $\RR_{\ge 0}^{d'}\cap S_1^{d'-1}$, with each collapse 
 preserving regularity and homeomorphism type, then the extension to 
 $\RR_{\ge 0}^d \cap S_1^{d-1}$ of this same series of collapses may be
 successfully performed on $\RR_{\ge 0}^d \cap S_1^{d-1}$.
 \end{lemma}
 
 \begin{proof}
 First apply each collapse to the subcomplex of $\RR_{\ge 0}^d \cap S_1^{d-1}$
 in which $t_{d'+1}=\cdots = t_d = 0$, since this is exactly $\RR_{\ge 0}^{d'}\cap S_1^{d'-1}$.
 Then extend continuously to the family of subspaces with $t_{d'+1}= \cdots = t_d = k$ for 
 $0\le k \le \frac{1}{d-d'+1}$, using that each of these for $k<\frac{1}{d-d'+1}$  is also isomorphic to 
 $\RR_{\ge 0}^{d'}\cap S_1^{d'-1}$ and that the slice given by $k=\frac{1}{d-d'+1}$ is a 1-point
 space.  Geometrically, we are adding a cone point and extending
 the collapse to the coned complex.  Continuity follows from the fact that the level curves we
 collapse across each hold the values 
  $t_{d'+1},\dots ,t_d$ constant, hence the collapses apply to all (nontrivial) 
 cross-sectional slices in the  same way.
 \end{proof}

\begin{lemma}
\label{ioa}
Suppose that the $k$-th collapsing step uses deletion pair $\{ x_{i_{l_k}},x_{i_{r_k}} \}$ for 
$l_k<r_k$.  Then $\sim^{r_k-1} = \sim_m$ where $m$ is the largest possible positive 
integer such that the 
triple $(x(F_m), \{ x_{i_{l_m}}, x_{i_{r_m}} \} )$ has $r_m < r_k$.
\end{lemma}

\begin{proof}
Consider the map $f_{(i_1,\dots ,i_{r_k-1})}$ on just the leftmost $r_k-1$ positions
in our reduced word.  By induction on wordlength, specifically our assumption of all 
results in the paper for all $d' < d$, we may use Theorem ~\ref{main-theorem}
to deduce that  $\overline{f_{(i_1,\dots ,i_{r_k-1})}}$ is a homeomorphism from 
$(\RR_{\ge 0}^{r_k-1}\cap S_1^{r_k-2})/\sim $ to $Y_{s_{i_1}\cdots s_{i_{r_k-1}}}$.  By Lemma
~\ref{induct-on-d}, this 
means in particular that the collapses based on deletion pairs using only positions
$1,\dots ,r_k-1$ are enough to accomplish all the desired identifications in $\sim^{r_k-1}$.
\end{proof}

\begin{lemma}\label{collapse-order}
For each $s<d$, there is a cell structure preserving homeomorphism $ch$ from
$(\RR_{\ge 0}^d \cap S_1^{d-1})/\sim^s_{(i_1,\dots ,i_d)}$ to 
$(\RR_{\ge 0}^d \cap S_1^{d-1})/\sim^s_{(j_1,\dots ,j_d)}$ where 
$(j_1,\dots ,j_d)$ is obtained from $(i_1,\dots ,i_d)$ by  braid moves involving only
the leftmost $s$ letters in $(i_1,\dots , i_d)$.  
\end{lemma}

\begin{proof}
The case of short braid moves is obvious, since we just switch the order of the parameters.
Now by induction on $d$,  we may assume the main results of the paper for all $s<d$ in our
proof of all these results for our given $d$, 
provided we check the base case of the induction.  In particular, this inductive usage of 
Theorem ~\ref{main-theorem} for $s<d$ together with Lemma ~\ref{induct-on-d} and the
present lemma for $s<d$ gives 
that the two complexes under consideration, namely  
$(\RR_{\ge 0}^d \cap S_1^{d-1})/\sim^s_{(i_1,\dots ,i_d)}$ and 
$ (\RR_{\ge 0}^d \cap S_1^{d-1})/\sim^s_{(j_1,\dots ,j_d)}$, both   are homeomorphic
by cell preserving homeomorphisms to the join 
$Y_u * \Delta_{d-s-1}$ where $(i_1,\dots ,i_s)$ and $(j_1,\dots ,j_s)$
are both reduced words for Coxeter group element $u$ and $\Delta_{d-s-1}$ is a simplex
having $d-s$ vertices.   Thus, by composing homeomorphisms they are also 
homeomorphic to each other.
The base case of this induction follows immediately from Lemma ~\ref{single-braid}.
\end{proof}

\begin{remark}
For any braid relation $(s_is_j)^{m(i,j)} = e$ in $W$ and  
any positive reals $t_2,\dots ,t_{m(i,j)} > 0$, there
is a unique $(t_1',\dots ,t_{m(i,j)}')\in \RR^{m(i,j)}_{\ge 0}$ satisfying
$$x_i (0) x_j(t_2)x_i(t_3)\cdots  = x_j(t_1')x_i(t_2') x_j (t_3')\cdots ,$$ namely
$t_1' = t_2, t_2' = t_3 , \dots  , t_{m(i,j) -1}' = t_{m(i,j)},$ and $ t_{m(i,j)}' = 0$.  
\end{remark}

\begin{remark}\label{single-face-braid}
The proof of Lemma ~\ref{collapse-order}  immediately implies
for each face $F$ with 
associated word $(i_1',\dots ,i_{d'}')$  that is a subword of $(i_1,\dots ,i_d)$
and for any word $(j_1',\dots ,j_{d'}') = ch(i_1',\dots ,i_{d'}')$  obtained by
braiding the leftmost $s$ letters in $(i_1',\dots ,i_{d'}')$   for $s \le \min (d-1,d') $  using
the same braid moves as above 
that  $g^s_{(i_1',\dots ,i_{d'}')}(\overline F_k)$ is homeomorphic to 
$(\RR_{\ge 0}^{d'}\cap S_1^{d'-1})/\sim^s_{(j_1',\dots ,j_{d'}')}$.
\end{remark}

Next we check a condition that will be helpful for verifying
our various requirements for performing collapses, i.e. for checking the hypotheses of
Corollary ~\ref{more-general-topol-collapse2}, assuming that all earlier collapses
were performed successfully.
Note that in the Lemma ~\ref{same-dim-lemma} 
we do not require  $F$ to be the maximal face that is collapsed in this step.

\begin{lemma}\label{same-dim-lemma}
If a cell $F$ is collapsed across level curves each having one endpoint in $H_1$
and the other endpoint in $H_2$, then $\dim (H_1) = \dim (H_2) = \dim (F) -1$, 
with neither $H_1$ nor $H_2$ collapsed earlier.
\end{lemma}

\begin{proof}
Let $\{ x_{i_r}, x_{i_s} \} $ be the deletion pair inducing the collapse of $F$, with 
$x_{i_r}\in x(H_1)$ and $x_{i_s}\in x(H_2)$.
Then $x(H_1)|_{x_{i_1}\cdots x_{i_s}}$  must be reduced, 
since otherwise $H_1$ and likewise 
$F$ would have been collapsed earlier.   This implies
$x(H_2)|_{x_{i_1}\cdots x_{i_s}}$ must also be reduced. 
Thus, 
neither $H_1$ nor $H_2$ will have been collapsed earlier, from which the result
follows by comparing wordlengths.
\end{proof}

Now we check  the various hypotheses of  
Theorem ~\ref{topol-collapse2} (as extended in Corollary ~\ref{more-general-topol-collapse2}),
using our collapses and parallel-like curves from Definition ~\ref{general-level-def} and
from  the discussion just after Example ~\ref{long-move-example}
as applied to our framework.  To check these conditions 
for the  $k$-th collapse, we assume by induction that all  earlier collapses were performed
successfully.  
In particular, this means that we assume that we had a regular CW complex after
each earlier  collapsing
step and hence that (LUB) held after each earlier collapsing step.  However, 
we do not need to check (DIP) or (DE) for the 
curves used in the 
$k$-th collapse at each earlier step.  Rather, 
it suffices to check these conditions just before the 
$k$-th collapse for the curves used to accomplish the $k$-th collapse. 

The fact that the level curves to be used in our collapses, 
defined in Definition ~\ref{general-level-def}, 
satisfy the first requirement for parallel-like curves is immediate
from our set-up.    Now let us confirm that the level curves used for the $k$-th 
collapse also satisfy the third requirement of parallel-like curves, the (DIP) condition, just prior to 
the $k$-th collapsing step, using that earlier collapses were all performed successfully.

When a collapse identifies  cells
$A$ and $A'$ via a deletion pair $\{ x_{i_u},x_{i_v} \} $, we say that $x_{i_u}$ is {\it exchanged} 
for $x_{i_v}$, denoted $x_{i_u}\rightarrow x_{i_v} $.

\begin{lemma}\label{G_2-ends}
If a cell $H_1 \subseteq G_1 $ is collapsed prior to the collapse of $F_j$,
where $H_1$ is to be 
identified with $H_2$ in the collapsing step given by $(F_j,\{ x_{i_r},x_{i_s} \} )$ 
by an exchange of $x_{i_r}\in H_1 $ for $x_{i_s}\in H_2$ for $r<s$, 
then $H_1\vee H_2$ is also 
collapsed prior to $F_j$, and in such a way that any two level
curves with the same endpoint in $H_1$ will have already been identified with 
each other in a manner that preserves the parametrization.
\end{lemma}

\begin{proof}
Given the two $x$-expressions 
$x(H_2)|_{x_{i_1}\cdots x_{i_s}} = x_{i_{j_1}}\cdots \hat{x}_{i_r}\cdots x_{i_s}$ and 
$x(H_1)|_{x_{i_1}\cdots x_{i_s}} = x_{i_{j_1}}\cdots x_{i_r}\cdots \hat{x}_{i_s}$, the fact
that $H_1$ has already been collapsed means  $x_{i_{j_1}}\cdots x_{i_r}\cdots \hat{x}_{i_s}$
is not reduced, implying that $x_{i_{j_1}}\cdots x_{i_r} \cdots x_{i_s} = x(H_1 \vee H_2)$ 
also is not reduced.   Moreover, $H_1$ will have been collapsed based on a deletion pair
strictly to the left of $x_{i_s}$, which implies the same for $H_1\vee H_2$.
Our collapsing order ensures 
that $H_1 \vee H_2$ will also have been collapsed prior to the collapse
of $F_j$,  using the same deletion pair and the same series of braid moves
as in $H_1$, hence the same parametrization for each curve in $H_1$ as in the curves with 
which it is identified in $H_1\vee H_2$.
\end{proof}

Next we verify that  the second requirement for parallel-like curves, Condition (DE), holds just 
prior to the $k$-th collapse for the curves used to accomplish the $k$-th collapse.
 
\begin{lemma}\label{distinct}
Condition ~\ref{distinct-endpoint-condition} (DE) holds at the $k$-th collapsing step, 
provided that the earlier collapses were performed successfully.
\end{lemma}

\begin{proof}
Suppose $F$ is  collapsed during the $k$-th collapsing step
via deletion pair $\{ x_{l_k},x_{r_k} \} $ having parameters
$\{ t_{l_k},t_{r_k} \} $.  We are not assuming $F$ is necessarily the maximal face $F_k$
among those faces which are collapsed at this step.    Let 
$\{ u_{r_k-1},u_{r_k} \} $ be the new parameters for the stuttering pair obtained from the deletion
pair $\{ x_{l_k}, x_{r_k} \} $  by a   suitable coordinate change as in Lemma
~\ref{collapse-order}.  Let $G_1$ and $G_2$ be the closed faces 
containing opposite endpoints of  the   curves  across 
which $F$ is collapsed, i.e. the curves introduced in Definition 
~\ref{general-level-def}.  Thus,   $G_1$ has $t_{r_k}=0$ and
$G_2$ has $t_{l_k}=0$, with $x(G_1), x(G_2)$ and  $x(F)$ agreeing  at all other
positions. 
In the new coordinates, $G_1$ has $u_{r_k}=0$ while $G_2$ instead has $u_{r_k-1}=0$.

What we must  prove is
that $G_1$ and $G_2$ are not identified in an earlier collapse.  
Suppose otherwise.  By induction, we may assume  
that the complex is regular immediately after each earlier collapsing
step.  This  precludes $G_1$ and
$G_2$ from being incomparable in the closure poset just prior to their
identification, unless
there is a face $G$ having $G_1,G_2\subseteq \overline{G} $ which is also identified with 
both of them at that same earlier step by collapsing $G$ across 
parallel-like  curves each
having one  endpoint in $G_1$ and the other endpoint in $G_2$.
But  then we may use (LUB) for that earlier collapse  
to deduce that all least upper bounds for 
$G_1$ and $G_2$ must have been collapsed at this earlier step, which would
necessarily include a face in $\overline{F}$ other than $F$ itself.  But this implies that
$G_1$ and $G_2$ both have dimension at least two less than $F$, contradicting
what we just proved in Lemma ~\ref{same-dim-lemma}. 
Suppose on the other hand we have 
$G_1\subseteq \overline{G_2}$ or $G_2\subseteq \overline{G_1}$ just prior to their
identification.   But this implies $\dim (G_1) \ne \dim (G_2)$, again contradicting our result from 
Lemma ~\ref{same-dim-lemma} that $\dim (G_1) = \dim (G_2) = \dim (F)-1$.
\end{proof}

In the proof of the next lemma, 
we often speak of a cell $A \in (\RR_{\ge 0}^d \cap S_1^{d-1})/\sim_k$,  
whereby  we mean the equivalence class  under $\sim_k$
that  contains the cell $A$ from 
$\RR_{\ge 0}^d \cap S_1^{d-1}$.  Referring in this manner to  particular
equivalence class representatives has the advantage that there is a unique associated
 $x$-expression $x(A)$ for that representative and also allows us to make sense of
 $A\vee B$ as being the equivalence class under $\sim_k$ that includes the join of the
 representatives $A$  and $B$ where join is taken in the original  simplex 
 $\RR_{\ge 0}^d \cap S_1^{d-1}$.  Thus, 
 $x(A)\vee x(B) = x(A\vee B)$ is the $x$-expression comprised of the union of the 
 $x$-expressions for $A$ and $B$.   In particular, $x(F_k)$ denotes the $x$-expression for
 the $\sim_k$-equivalence class representative  with respect to which $F_k$ is  collapsed, 
 namely one yielding optimal $x$-expression in terms of our collapsing order.

\begin{lemma}\label{stay-regular}
Condition ~\ref{lub-condition} (LUB) holds at the $k$-th collapsing step, provided that 
all  earlier collapsing steps were performed successfully.
\end{lemma}

\begin{proof}
Suppose that the collapsing step given by a triple $(F_k,\{ x_{i_{l_k}},x_{i_{r_k}}\} )$  causes a 
pair of cells $A$ and $B$ to be identified where neither cell has been collapsed yet and 
neither is in the closure of the other (the latter of
which would make the result trivial).  By virtue of our collapsing process,
this implies that we may choose $x(A)$ and $x(B)$ to be 
subexpressions of $x(F_k)$ which coincide except in that $x(A)$ includes the letter 
$x_{i_{r_k}}$ whereas $x(B)$ instead includes the letter $x_{i_{l_k}}$.  In other words, 
there must exist $\sim_k$ equivalence class representatives with this property.  
Let $F$ be the cell
such that $x(F) =  x(B)\vee \{ x_{i_{r_k}} \}  = x(A) \vee \{ x_{i_{l_k}} \} $.  We note that
$F$ might equal $F_k$,  or $F$  might  be a lower-dimensional
cell contained  in the  closure $\overline{F}_k$.    
By definition, $F$ has dimension exactly one more than $A$ and 
$B$.  Our collapsing order ensures that $F$ 
could not have been collapsed prior to 
our current step collapsing $F_k$ unless it were done 
through the earlier identification of $F$ with
another cell  $F'$ designated for collapse earlier than $F$,  where $F'$ is obtained from
$F$ by replacing $x_{i_{r_k}}$ with some $x_{i_{l'}}$ with $l_k < l' < r_k$;  
however, such a step would have also identified $A$ with $B$ earlier, a contradiction. 
Thus, we are assured of the existence of such an $F$ which is not collapsed prior 
to step $k$ and which satisfies  
$x(F) = x(A)\vee x(B) $ and $\dim (F) = \dim (A) + 1 = \dim (B) + 1$. 

By virtue of our collapsing process, 
the collapse of $F_k$ will induce at the same time the collapse of $F$ 
across  curves  as defined in Definition ~\ref{general-level-def} (which we will 
eventually prove are parallel-like) with 
each curve having one endpoint in $A$ and the other endpoint in $B$.
Lemma ~\ref{ioa} allows us to assume that all possible identifications  present in 
$\sim^{r_k-1}$   have already been accomplished  before the $k$-th
collapse.  Let $x(F')$ be an $x$-expression for 
any other cell besides $F$  that is a least upper bound for $A$ and $B$ just 
prior to the  $k$-th collapse.  If there is such an $F'$ distinct from $F$, then there must be
cells $A_u, B_v \subseteq \overline{F'}$ having $x(A_u)$ and $x(B_v)$ as subexpressions
of $x(F')$ with $x(F') = x(A_u)\vee x(B_v)$, 
$A\sim_k A_u$ and $B\sim_k B_v$.  
Our plan now is
to obtain from this a contradiction so as to deduce that no such $F'$ exists, thereby 
proving (LUB).  
To this end, we will prove  either (1) that $F'$ was already identified with
$F$ at an earlier step,  (2) that $F'$ is not  a least upper bound for $A$ and $B$ just prior to the 
$k$-th collapse due to the presence of an upper bound strictly contained in it, or (3) that
$A_u\sim_k  B_v$.

The fact that we have only done identifications
based on deletion pairs whose right endpoint is at or to the left of $x_{i_{r_k}}$ 
prior to the $k$-th collapsing step
implies  that $A, A_u, B, $ and $B_v$ all must coincide with  each other to the right
of $x_{i_{r_k}}$.  We also know that $w(A_u|_{x_{i_1}\cdots x_{i_{r_k}} }) = w(A|_{x_{i_1}\cdots x_{i_{r_k}}})$
and  $w(B_v|_{x_{i_1}\cdots x_{i_{r_k}}}) = w(B|_{x_{i_1}\cdots x_{i_{r_k}} })$, since 
$A\sim_k A_u$ 
and $B\sim_k B_v$.
By definition of deletion pair, we have 
$w(A|_{x_{i_1}\cdots x_{i_{r_k}}}) = w(B|_{x_{i_1}\cdots x_{i_{r_k}}})$, which implies
$w(A_u|_{x_{i_1}\cdots x_{i_{r_k}}}) = w(B_v|_{x_{i_1} \cdots x_{i_{r_k}}})$ as well.

Now consider the case where $x_{i_{r_k}}\not \in  x(F')  = x(A_u)\vee x(B_v)$.  
It follows immediately from 
$w(x(A)|_{x_{i_1}\cdots x_{i_{r_k}}})  = w(x(B)_{x_{i_1}\cdots x_{i_{r_k}}})$ that we also
have $w(A_u|_{x_{i_1}\cdots x_{i_{r_k-1}}}) = w(B_v|_{x_{i_1}\cdots
x_{i_{r_k-1}}})$.  But now we may use Lemma ~\ref{ioa}
to conclude that $A_u$ 
must get identified
with $B_v $  under one of the steps leading to $\sim^{r_k-1}$ 
and hence prior to step $k$ .
But by Lemma ~\ref{induct-on-d}, this implies $A_u\sim_k B_v$, completing
this case.

The remainder of the proof deals with the case 
$x_{i_{r_k}} \in x(F') =  x(A_u\vee B_v) = x(A_u)\vee x(B_v)$.  
First suppose that both $x(A_u)$ and $x(B_v)$ include the letter $x_{i_{r_k}}$.  Notice that $x(A_u)|_{x_{i_1}\cdots x_{i_{r_k-1}}}$ and $x(B_v)|_{x_{i_1}\cdots x_{i_{r_k-1}}}$ must both be reduced, since neither of these cells  has  been collapsed yet.  
We also have $w(A_u|_{x_{i_1}\cdots x_{i_{r_k}}}) = w(B_v|_{x_{i_1}\cdots x_{i_{r_k}}})$.  Hence, 
if  $x_{i_{r_k}}$ is nonredundant in both $A_u$ and $B_v$, then this together implies $w(A_u|_{x_{i_1}\cdots x_{i_{r_k-1}}}) = w(B_v|_{x_{i_1}\cdots x_{i_{r_k-1}}})$.  These two words
$x(A_u)$ and $x(B_v)$ also coincide on the subexpression consisting of $x_{i_{r_k}}$ and all  letters to its right, implying $A_u\sim^{r_k-1} B_v$
  and hence  $A_u \sim_k B_v$, a 
  contradiction.   Next suppose that $x_{i_{r_k}}$ is present but redundant in 
$x(B_v)$, that is, suppose $w(x(B_v)|_{x_{i_1}\cdots x_{i_{r_k-1}}}) = 
w(x(B_v)|_{x_{i_1}\cdots x_{i_{r_k}}})$;  also suppose $x_{i_{r_k}}$ is  
present in $x(A_u)$.  
Then consider $x(B_v')$ obtained from $x(B_v) $ by deleting $x_{i_{r_k}}$.
Then $A_u\vee B_v = A_u \vee B_v'$, which means that it suffices to 
show that $A_u\vee B_v'$ is collapsed by the end of the $k$-th collapsing step so as to deduce that $A_u \vee B_v'$ is also  collapsed by the end of the $k$-th collapsing step.  Our remaining arguments will cover the case of such $A_u\vee B_v'$.

Henceforth, we 
assume  $x_{i_{r_k}}\in x(A_u)$ and $x_{i_{r_k}}\not\in x(B_v)$.  Notice that $x(A_u)\vee x(B_v)|_{x_{i_1}\cdots x_{i_{r_k-1}}}$ must be reduced, since otherwise $A_u\vee B_v$ would have already been collapsed by the series of collapses yielding $\sim^{r_k-1}$, by virtue of our collapsing order and our use of induction on length.  
Thus, we have two cases left to consider, depending whether (a) $(x(A_u \vee B_v)) |_{x_{i_1}\cdots x_{i_{r_k}}}$ is also reduced, or 
(b) $x_{i_{r_k}}$ forms a deletion pair with a letter to its left in $x(A_u \vee  B_v)$.  

Now to (a),  namely  the case where $x(A_u\vee B_v)|_{x_{i_1}\cdots x_{i_{r_k}}}$ is 
reduced.  The Coxeter group  element $w(x(A_u)|_{x_{i_1}\cdots x_{i_{r_k}}}) = w(x(B_v)|_{x_{i_1}\cdots x_{i_{r_k}}}) = w(x(B_v)|_{x_{i_1}\cdots x_{i_{r_k-1}}})$ must in this case  
be strictly less than $ w((A_u\vee B_v)|_{x_{i_1}\cdots x_{i_{r_k}}}) $ in Bruhat order.  It 
follows from this that reading 
$x(A_u\vee B_v) |_{x_{i_1}\cdots x_{i_{r_k}}}$ from left to right, we must encounter a
leftmost letter 
$x_{i_j}$ whose associated reflection is not one of the associated reflections for 
$x(A_u)|_{x_{i_1}\cdots x_{i_{r_k}}}$ and likewise for $x(B_v)|_{x_{i_1}\cdots x_{i_{r_k}}}$, hence a letter $x_{i_j}$ which may be deleted to obtain a new $x$-expression $x(A_u\vee B_v)|_{x_{i_1}\cdots x_{i_{r_k}}} \setminus x_{i_j}$  where just the letter $x_{i_j}$ has been deleted; 
this expression will have the property that its 
associated Coxeter group element is again greater than or equal to  both $w(x(A_u)|_{x_{i_1}\cdots x_{i_{r_k}}})$ and  $w(x(B_v)|_{x_{i_1}\cdots x_{i_{r_k}}})$ in Bruhat order.   
Choose the leftmost such letter $x_{i_j}$

If  $j\ne r_k$, 
we must  therefore have  $A_u\sim^{r_k-1} A'$ for some $A'$ such that $x(A')$ is a subexpression of  $x(A_u\vee B_v)  \setminus x_{i_j}$, 
with $x(A')$  including $x_{i_{r_k}}$, omitting $x_{i_j}$,  and satisfying $w(x(A')|_{x_{i_1}\cdots x_{i_{r_k-1}}}) = w(x(A_u)|_{x_{i_1}\cdots x_{i_{r_k-1}}})$,  by virtue of 
Lemma ~\ref{ioa} together with the definition  of Bruhat order as an order based on taking subwords.  Similarly, there must exist $B'$  with $x(B')$ also a subexpression of 
$x(A_u\vee B_v)  \setminus x_{i_j}$ 
such that  $B_v\sim^{r_k-1} B'$ for $x(B')$ omitting both  $x_{i_j}$ and  $x_{i_{r_k}}$, with $w(x(B_v)|_{x_{i_1}\cdots x_{i_{r_k-1}}}) = w(x(B')|_{x_{i_1}\cdots x_{i_{r_k-1}}})$.
Thus, $x(A_u\vee B_v)  \setminus x_{i_j}$ 
gives an upper bound for $x(A_u)$ and $x(B_v)$ which is strictly contained in $x(F') = x(A_u \vee B_v)$ with $w(x(A_u \vee B_v))\ne w(x(A_u\vee B_v)\setminus x_{i_j})$, contradicting $F'$ being a least upper bound for $A_u$ and $B_v$.   

If $j=r_k$, then 
$l(w(x(A_u\vee B_v)|_{x_{i_1}\cdots x_{i_{r_k}}})) = l(w(x(A_u)|_{x_{i_1}\cdots x_{i_{r_k}}}))+1 = l(w(x(B_v)|_{x_{i_1}\cdots x_{i_{r_k}}}))+1$ where $l$ denotes  Coxeter-theoretic length.  But since $x(A_u\vee B_v)|_{x_{i_1}\cdots x_{i_{r_k}}}$ is reduced, we cannot delete from it one letter to obtain 
$x(A_u)|_{x_{i_1}\cdots x_{i_{r_k}}}$ and a different individual  letter to obtain $x(B_v)|_{x_{i_1}\cdots x_{i_{r_k}}}$, since that would imply $w(x(A_u)|_{x_{i_1} \cdots x_{i_{r_k}}}) \ne  w(x(B_v)|_{x_{i_1}\cdots x_{i_{r_k}}})$ by the exchange axiom for Coxeter groups, a contradiction.  This completes our proof in case  (a).

Now to case (b), i.e. the case where $x_{i_{r_k}}$ is redundant in $x(A_u\vee B_v)|_{x_{i_1}\cdots x_{i_{r_k}}}$ and where $x(A_u\vee B_v)|_{x_{i_1}\cdots x_{i_{r_k-1}}}$ is reduced.   
We then have  $A_u\vee B_v \sim^{r_k-1}  A\vee B$ unless 
$w(x(A_u\vee B_v)|_{x_{i_1}\cdots x_{i_{r_k-1}}}) \ne w(x(A\vee B) |_{x_{i_1}\cdots x_{i_{r_k-1}}})$, since $x_{i_{r_k}}$ appears  both in $x(A\vee B)|_{x_{i_1}\cdots x_{i_{r_k}}}$ and in $x(A_u\vee B_v)|_{x_{i_1}\cdots x_{i_{r_k}}}$ and is redundant in both.  But 
$w(x(A\vee B)|_{x_{i_1}\cdots x_{i_{r_k}}})$ 
also equals $w(B_v|_{x_{i_1}\cdots x_{i_{r_k-1}}})$ which is a subword of 
$w(x(A_u\vee B_v)|_{x_{i_1}\cdots x_{i_{r_k-1}}})$, so this implies that $w(x(A\vee B)|_{x_{i_1}\cdots x_{i_{r_k}}})$  must be strictly less than
$w(x(A_u\vee B_v)|_{x_{i_1}\cdots x_{i_{r_k-1}}})$   in Bruhat order.   But this means there is a letter $x_{i_j}$ we may delete from $x(A_u\vee B_v)|_{x_{i_1}\cdots x_{i_{r_k-1}}}$ to obtain a word whose associated Coxeter group element is still greater than or equal to $w(x(B_v)|_{x_{i_1}\cdots x_{i_{r_k-1}}})$ in Bruhat order, and hence is also greater than or equal to $w(x(A_u)|_{x_{i_1}\cdots x_{i_{r_k-1}}})$ in Bruhat order since 
$$w(x(A_u)|_{x_{i_1}\cdots x_{i_{r_k}-1}}) \le_{Bruhat} w(x(A_u)|_{x_{i_1}\cdots x_{r_k}}) = w(x(B_v)|_{x_{i_1}\cdots x_{i_{r_k-1}}}).$$  But then deleting this $x_{i_j}$ from $x(A_u\vee B_v)$ yields an upper bound for $A_u$ and $B_v$ just prior to the $k$-th collapse whose $x$-expression is strictly contained in $x(A_u\vee B_v)$ with distinct associated Coxeter group elements, i.e. with $w(x(A_u\vee B_v)\setminus x_{i_j}|_{x_{i_1}\cdots x_{i_{r_k-1}}}) \ne w(x(A_u\vee B_v)|_{x_{i_1}\cdots x_{i_{r_k}-1}})$.  
Thus, we get a cell that is an upper bound for $A_u$ and $B_v$ that is strictly contained in $A_u\vee B_v$, contradicting $A_u\vee B_v$ being a least upper bound.   This completes case (b).
\end{proof} 

Next, we give a projection map  $\pi^v_u$ from any closed cell $\overline{\sigma}_v $ 
 in $(\RR_{\ge 0}^d \cap S_1^{d-1})/\sim $
onto an open cell  $\sigma_u$ in its boundary, choosing our notation to reflect that
$\overline{\sigma_v}$ is mapped by 
$\overline{f_{(i_1,\dots ,i_d)}}$ to $Y_v$ while $\sigma_u$  is mapped by 
$\overline{f_{(i_1,\dots ,i_d)}}$  to $Y_u^o$.  To be more precise, $\pi_u^v$ applies to 
the union of open cells contained in $\overline{\sigma}_v$ 
that have $\sigma_u$ in their closure.
Once equipped with this projection map, we may  define the links of cells using the 
notion of link provided  in Definition ~\ref{strat-link-def}
which is based on ideas from stratified Morse theory (cf. \cite{GM}).
It is most natural to define this projection map using all of $\RR_{\ge 0}^d$ at once 
and its quotient spaces $\RR_{\ge 0}/\!\!\sim $ as well as $\RR_{\ge 0}/ \!\!\sim_k$ 
induced by the collapsing process that we have 
developed over the course of this section.   
This is the approach that we take.

\begin{definition}\label{proj-map}
{\rm
Consider any $x = (t_1,\dots ,t_d)$ in the aforementioned domain of $\pi^v_u$, 
choosing the representative for the equivalence class of $x$ under $\sim $
whose only nonzero parameters appear in the positions of letters 
in  the rightmost reduced word 
for $v$ appearing as a subword of $(i_1,\dots ,i_d)$.  
Let $(i_{j_1},\dots ,i_{j_r})$ be this subword. 
The {\it projection map} given by $u\subseteq \overline{v}$, 
denoted $\pi^v_u$, 
sets to 0 each $t_i$ not appearing in either the rightmost 
subword of $(i_{j_1},\dots ,i_{j_r})$
that is a reduced word for $u$ or obtained by reading our word from right to left, including
additionally just those letters which  are redundant (in the sense of 
not increasing the $0$-Hecke algebra theoretic  length)
when appended
to the word comprised of  those
letters to its right that have already been chosen. 
}
\end{definition}

For intermediate
stages in the collapsing process, instead of projecting onto a cell indexed by $u\in W$, we 
project to a $\sim_k$-equivalence class of faces of the simplex, each of which  maps to $Y_u^o$ under 
$\overline{f_{(i_1,\dots ,i_d)}}$.  The  corresponding intermediate projection maps are defined 
completely analogously to $\pi_u^v$,  but now using only those 
subwords belonging to the allowed equivalence
classes mapping under $\overline{f_{(i_1,\dots ,i_d)}}$  
 to cells $Y_{u'}^o$ with $u\le u' \le v$.

We now prove the main theorem of this section, which is largely a matter of 
pulling together the various lemmas we have 
proven already.  This will require some further notation.
 Let $(i'_1,\dots ,i'_{d'})$ be the  subword of  
 $(i_1,\dots ,i_d)$  associated to 
  the maximal cell $F_k$ to be  collapsed at the $k$-th collapsing step, 
 so $d'$ is the 
wordlength of $x(F_k)$.  Let $(j'_1,\dots ,j'_{d'})$ be 
the word obtained from $(i'_1,\dots ,i'_{d'})$ by
applying the series of braid moves giving rise to the change of coordinates 
homeomorphism $ch$, i.e. the chosen 
series of (long and short) braid moves  applied to $(i_1',\dots ,i_{d'}')$ yielding a stutter 
between positions $r-1$ and $r$ to be used to induce our collapse of $F_k$.  Note that
there must be some $r$ such that 
$(i_r',\dots ,i_{d'}') = (i_{r_k},\dots ,i_d) = (j_r',\dots ,j_{d'}')$.   The restriction of
$\sim^{r_k-1}_{(i_1,\dots ,i_d)}$ to $\overline{F}_k$ equals the equivalence relation 
$\sim^{r-1}_{(i_1',\dots ,i_{d'}')}$.  This  
also carries out 
exactly the same identifications as $\sim^{r-1}_{(j_1',\dots ,j_{d'}')}$ after suitable change
of coordinates.

\begin{theorem}\label{homeom2-lemma}
$(\RR_{\ge 0}^d \cap S^{d-1}_1 )/\!\!\sim $ is a regular CW complex homeomorphic
to a ball.  Moreover, the link of any cell is also a regular CW complex homeomorphic
to a ball.
\end{theorem}

\begin{proof}
We start our proof with the simplex 
$K_0 = \RR_{\ge 0}^d \cap S^{d-1}_1$.  We perform on $K_0$  a series of collapses, using 
Theorem ~\ref{topol-collapse2} and its generalization described in 
Corollary ~\ref{more-general-topol-collapse2} 
to justify that each of these collapses preserves homeomorphism type as well
as the property of having a regular CW complex.  These  collapses use the families 
of curves introduced in Definition ~\ref{general-level-def}.
We assume inductively that all earlier 
collapses were performed successfully in order to justify that the $k$-th collapse preserves
homeomorphism type, regularity, and all the requisite properties for our inductive step.  
We also assume inductively all results in the paper for all strictly
smaller $d$.   

We use the cell collapsing order given just
after Example ~\ref{long-move-example}.
The $k$-th collapsing step is specified by a deletion pair $\{ x_{i_{l_k}},x_{i_{r_k}} \} $ in an 
$x$-expression $x(F_k)$, with the cell 
$\overline{F}_k/\sim_k  \in (\RR_{\ge 0}^d \cap S^{d-1}_1)/\sim_k$ 
collapsed across curves each having one endpoint in the 
closed cell $\overline{G}_1 /\sim_k $  which consists of the points of $\overline{F_k}/\sim_k$
with  $t_{l_k}=0$ and the other endpoint in the closed
cell $\overline{G}_2/\sim_k$ instead having $t_{r_k}=0$.  The series of lemmas we have
just proven will allow us to verify all the requirements for our curves to be a 
parallel-like family of curves (in the extended sense of Corollary ~\ref{more-general-topol-collapse2}) and to check all of the hypotheses of Theorem ~\ref{topol-collapse2} (again as  
extended in Corollary  ~\ref{more-general-topol-collapse2}),
once we show how to incorporate the requisite 
change of coordinates homeomorphisms $ch$  into the picture.  This change of 
coordinates which will help us accomplish the $k$-th collapse
will be done on $(\RR_{\ge 0}^d \cap S_1^{d-1})/\sim^{r_k-1}$ to create a stutter, 
so in particular will be  done
prior to the $k$-th collapse.

The fact that $i_{r_k} = i_r'$ is the leftmost 
right endpoint of a deletion pair in $x(F_k)$ implies that $(i_1',\dots ,i_{r-1}')$ is reduced 
and hence that  $(j_1',\dots ,j_{r-1}')$ is also reduced with
$w(i_1',\dots ,i_{r-1}') = w(j_1',\dots ,j_{r-1}')$ and hence $w(i_1',\dots ,i_{d'}') = w(j_1',\dots ,
j_{d'}')$. 
Let us also now choose a  reduced word 
$(j_1,\dots ,j_d)$ such that $j_i = j_i'$ for all $i\le d'$.   The point of extending 
$(j_1',\dots ,j_{d'}')$ to this longer reduced
word $(j_1,\dots ,j_d)$  is to  have at our disposal  a regular CW 
ball $(\RR_{\ge 0}^d\cap S_1^{d-1})/\sim^{r-1}_{(j_1,\dots ,j_d)}$ 
given by $(j_1,\dots ,j_d)$ with this ball  of the correct dimension so as to be 
homeomorphic to the complex $(\RR_{\ge 0}^d \cap S_1^{d-1})/\sim_m$ given by
$(i_1,\dots ,i_d)$  for each $m > 0$ 
and which for $m\ge k$ will also have $(\RR_{\ge 0}^{d'}\cap S_1^{d'-1})/\sim_{(j_1',\dots ,j_{d'}')}$ as a subcomplex.
 
The image under $ch $
of each curve in our family covering $\overline{F}_k/\sim^{r-1}_{(i_1',\dots ,i_{d'}')} $
will be a collection of points $(t_1',\dots ,t_d')$ with 
$t_{r_k-1}' + t_{r_k}' $ held constant and each $t_i'$ for $i\not\in \{  r_k-1,r_k \} $ also held
constant, all within  $(\RR_{\ge 0}^{d'}\cap S_1^{d'-1})/\sim^{r-1}_{(i_1',\dots ,i_{d'}')}$.  
This will imply 
that $ch$ maps our family of curves to the images under a series of collapsing
maps of a family of parallel line segments covering the cell having the structure
$(\RR_{\ge 0}^{d'}\cap S_1^{d'-1})/\sim^{r-1}_{(j_1',\dots ,j_{d'}')}$ 
within the regular CW complex  
$(\RR_{\ge 0}^d \cap S_1^{d-1})/\sim^{r-1}_{(j_1,\dots ,j_d)}$.  
In particular, $ch^{-1}$ also thereby will 
induce a transfer of  a parametrization function to each nontrivial curve in the closed cell
$\overline{F}_k/\sim_{(i_1,\dots ,i_d)}^{r_k-1}$ 
in $(\RR_{\ge 0}^d \cap S_1^{d-1})/\sim^{r_k-1}_{(i_1,\dots ,i_d)}$ which is a continuous
function to $[0,1]$ on the union of these nontrivial curves.  Moreover, condition (DE) together
with  the fact that 
collapses subsequent to those producing $\sim^{r_k-1}$ and leading to $\sim_k$  
restrict to homeomorphisms on each closed cell that 
they do not collapse will imply that  these collapses
will also carry forward these curve  parametrizations for all curves that  stay nontrivial under
the intermediate collapses to the quotient complex
$(\RR_{\ge 0}^d\cap S_1^{d-1})/\sim_k$  from the quotient complex
$(\RR_{\ge 0}^d\cap S_1^{d-1})/\sim^{r_k-1}$.   Thus, our 
parallel-like curves for $\overline{F}_k/\sim_k$ will be the images under $ch^{-1}$ of these 
parallel-like curves from $(\RR_{\ge 0}^d \cap S_1^{d-1})/\sim_{(j_1,\dots ,j_d)}^{r-1}$, 
pushed forward by the collapses yielding $\sim_k$ from $\sim^{r_k-1}_{(i_1,\dots ,i_d)}$.
Now let us carefully define the change of coordinates map $ch$ (and thereby $ch^{-1}$.

This  map  $ch$ is  most naturally 
defined on the closed cell  
$\overline{F}_k/\sim^{r_k-1}_{(i_1,\dots ,i_d)}$.
Lemma ~\ref{collapse-order} (as explained in 
Remark ~\ref{single-face-braid}) proves  $ch$  is a cell structure preserving
homeomorphism from $\overline{F}_k/\sim^{r-1}_{(i_1',\dots ,i_{d'}')} = 
\overline{F}_k /\sim_{(i_1,\dots ,i_d)}^{r_k-1}$  to 
$(\RR_{\ge 0}^{d'} \cap S_1^{d'})/\sim^{r-1}_{(j_1',\dots ,j_{d'}')}$.
Let us now show how $ch$ (or equivalently $ch^{-1}$) may be extended to a 
homeomorphism from  the
entire quotient complex $(\RR_{\ge 0}^d \cap S_1^{d-1})/\sim^{r_k-1}_{(i_1,\dots ,i_d)}$ 
to the quotient complex
$(\RR_{\ge 0}^d \cap S_1^{d-1})/\sim^{r_k-1}_{(j_1,\dots ,j_d)}$ given by the word
$(j_1,\dots ,j_d)$.
It will be  a necessity for our upcoming approach to collapsing $\overline{F}_k/\sim_k$ 
via transfer of parallel-like curves carried out by the 
map $ch^{-1}$ 
that $f_{(i_1',\dots ,i_{d'}')}(x)
= f_{(j_1',\dots , j_{d'}')}(ch(x))$ for each $x\in \overline{F}_k/\sim^{r_k-1}_{(i_1,\dots ,i_d)}$; 
however,  
this relationship between $x$ and $ch(x)$  will not 
be needed   for  
the extension of  $ch^{-1}$ to points $x$ 
outside of $\overline{F}_k/\sim_k$,  
because  the collapsing  map we will use to collapse 
$F_k$ will restrict to an injection
outside of $\overline{F}_k/\sim_k$. 
Thus,  we may simply 
extend  $ch^{-1}$ to a neighborhood $N$ of $\overline{F}_k/\sim_k$ 
by thickening the boundary of
$\overline{F}_k/\sim_k$ to $\partial (\overline{F}_k )/\sim_k \times [0,\epsilon) $ and letting $ch^{-1}(x,t) = (ch^{-1}(x),t)$ for
each $t\in [0,\epsilon)$. 
This thickening is possible since the closed complement of $\overline{F}_k$ within the boundary
of a cell of dimension one higher is a topological manifold with boundary, by Lemma 
~\ref{topol-collapse2} applied to the earlier collapse, hence has a collar by Theorem 
~\ref{collar-theorem}.

It is worth emphasizing   that the transfer of curve parametrizations resulting from $ch^{-1}$
will not in any way actually modify the collapsing procedure that we already indicated we 
would use and which has been described and analyzed in detail 
earlier in  this section of the paper.  Rather, the map $ch^{-1}$ on 
$\overline{F}_k$ is used to justify  that we 
indeed have parallel-like curves in a suitable sense (i.e. as in Corollary  ~\ref{more-general-topol-collapse2}) to enable the collapse of $F_k$  even when  long braid moves are needed 
to create a stutter in $x(F_k)$; 
we do this by giving an alternate way 
that we could have obtained the closed cell $\overline{F_k}/\sim^{r-1}_{(i_1',\dots ,i_{d'}')}$.  
This in turn gives an alternate way we 
could have  obtained the closed cell $\overline{F_k}/\sim_k$ by performing exactly the 
collapses  on $(\RR_{\ge 0}^d \cap S_1^{d-1})/\sim^{r_k-1}_{(i_1,\dots ,i_d)}$ (and in the process 
also on $\overline{F_k}/\sim^{r-1}_{(i_1',\dots ,i_{d'}')}$ that  we had planned for word 
$(i_1,\dots ,i_d)$ to get from the equivalence relation
$\sim^{r_k-1}$ to the equivalence relation $\sim_k$, but having incorporated the transfer 
map $ch^{-1}$  at the 
step where we had equivalence relation $\sim^{r_k-1}$.  In other words, we regard
$(\RR^d_{\ge 0}\cap S_1^{d-1})/\sim^{r_k-1}_{(i_1,\dots ,i_d)}$ 
as our starting point complex $K_0$
for purpose of justifying the collapse of $\overline{F}_k/\sim_k$.

Having now realized our curves from Definition ~\ref{general-level-def}
as the images of parallel line segments in a suitable sense,
i.e. with a transfer from one complex to another complex potentially involved, 
we now turn to the remaining requirements for parallel-like curves.  
Lemmas  ~\ref{distinct},  ~\ref{stay-regular} and ~\ref{lub-condition}
confirm that the requirements  (DIP),  (DE)  (as defined in  Conditions ~\ref{dip-def} and 
~\ref{distinct-endpoint-condition}, respectively) and condition (LUB)  hold for the 
family of curves to be collapsed in the $k$-th collapse just prior to this collapse, assuming
all earlier collapses were done successfully. 
Thus, all of the requirements  of  Theorem ~\ref{topol-collapse2}  (or at least 
their relaxations as in Corollary ~\ref{more-general-topol-collapse2}) are met, enabling
us to repeatedly collapse cells  until all cells given by non-reduced subwords of 
$(i_1,\dots ,i_d)$ have been eliminated, preserving homeomorphism type and regularity 
at each step.    Thus, the end result is a regular CW complex homeomorphic to a ball.
Regularity  and homeomorphism type for links, as defined in  
Definition ~\ref{strat-link-def},   
follow from Lemma ~\ref{same-dim-lemma}, since  it allows us to invoke
Lemma  ~\ref{regular-links}, noting that the transversality requirements follow easily 
from the definition of  our series of projection maps.
\end{proof}

 \subsection{Regularity and homeomorphism type of $Y_w$}
\label{final-section}

We now finally turn to studying  the topological structure of 
$Y_w$ itself.  We  will use  the fact  that 
$(\RR_{\ge 0}^d \cap S_1^{d-1})/\sim $ is a regular CW complex homeomorphic to a 
ball  to prove now  that $Y_w$ is  as well. 
First we must
 verify Condition 4 of Theorem ~\ref{sufficient-theorem} for the characteristic maps 
 $\overline{f_{(i_1,\dots ,i_d)}}: (\RR_{\ge 0}^d \cap S_1^{d-1})/\sim \rightarrow Y_w$ 
 and their restrictions to the closed cells of $(\RR_{\ge 0}^d \cap S_1^{d-1})/\sim $.
 This is done in Lemma 
 ~\ref{injective-condition} below, which 
will only require the following properties of $\sim $, which are immediate from the
definition of  $\sim $:
\begin{enumerate}
\item
Each  $p \in \partial (\RR_{\ge 0}^d \cap S^{d-1}_1 )$ whose $x$-expression
is not reduced
is identified by $\sim $ with a point having more parameters set to $0$
\item
$p\sim q$ implies
$w(f_{(i_1,\dots ,i_d)}(p))=w(f_{(i_1,\dots ,i_d)}(q))$.
\end{enumerate}

The points
in a cell boundary, i.e. the preimage of one of the attaching maps, 
are obtained by setting
 some positive parameters  to 0.

 \begin{lemma}\label{injective-condition}
 Given a reduced word $(i_1',\dots ,i_{d'}')$ which is a subword of  reduced word $(i_1,\dots ,i_d)$, 
 then $\overline{f_{(i_1,\dots ,i_d)}}$ restricted to the codimension one faces of
 $F = \overline{R_{\{ i_1',\dots ,i_{d'}'\} }}/\sim $ 
 is an injection into  $Y_{s_{i_1'}\cdots s_{i_{d'}'}}$.
 \end{lemma}
 
 \begin{proof}
Notice first
that $\overline{f_{(i_1,\dots ,i_d)}}|_F = 
\overline{f_{(i_1',\dots ,i_{d'}')}}$. 
    By Lemma
 ~\ref{red-lemma}, this means that
 $w(x_{i_1'}\cdots \hat{x}_{i_r'}\cdots x_{i_{d'}'}) \ne w(x_{i_1'}\cdots \hat{x}_{i_{r'}'}
 \cdots x_{i_{d'}'})$ for $r\ne r'$, provided  that both expressions are reduced, since 
 that implies  that the map $w$ 
 just replaces each $x_i$ by simple reflection $s_i$.
 Consequently,   boundary points 
 obtained by sending distinct single parameters to 0 to obtain
 reduced expressions of length one shorter must belong  to  distinct cells, hence
 have distinct images under $\overline{f_{(i_1', \dots ,i_{d'}')}}$.
 On the other hand, varying values of the nonzero parameters while keeping fixed which 
 parameters are   0 and which are nonzero with the subexpression of nonzero parameters
 a  reduced expression must also yield  points with distinct images under
 $\overline{f_{(i_1',\dots ,i_{d'}')}}$, 
 since Lusztig proved that  $f_{(i_1',\dots ,i_{d'}')}$ acts homeomorphically on $\RR_{>0}^s$ 
 for $(i_1',\dots ,i_{d'}')$ reduced.
 Combining yields that  $\overline{f_{(i_1',\dots ,i_{d'}')}}$ is injective upon restriction to the codimension  one cells, as desired.
\end{proof}

Now to our main result,  Theorem ~\ref{main-theorem}.  
It is phrased in a somewhat technical way so as to enable proof by induction on the 
length $d$, and also to overcome the challenge that it was not known 
previously even that   
$\overline{f_{(i_1,\dots ,i_d)}}((\RR_{\ge 0}^d \cap S_1^{d-1})/\sim )$
 was a CW complex.  Theorem ~\ref{main-theorem}
 is immediately followed by corollaries with more natural statements.
  
 To use Theorem ~\ref{sufficient-theorem} in proving Theorem ~\ref{main-theorem}, 
 we will need
the preimages of the various characteristic maps to be closed cells in 
$(\RR_{\ge 0}^d \cap S_1^{d-1})/\sim $, since this will give condition 5 of our regularity
criterion.
It is not clear that taking  the closure in $(\RR_{\ge 0}^d \cap S_1^{d-1})/\sim $
of an open cell which is sent by $\overline{f_{(i_1,\dots ,i_d)}}$ 
to $Y_{\sigma }$ for  a Coxeter group element
$\sigma $ of length $d' < d$ is the same as constructing a complex 
$(\RR_{\ge 0}^{d'} \cap S_1^{d'-1})/\sim $ for $\sigma $ itself directly; we overcome 
this issue by allowing  flexibility in the choice of $\sim $ in the statement of
the next theorem, in particular allowing
the collapsing maps for the closure of a cell which is not top-dimensional to be 
induced from the collapsing maps on the entire complex.
 
\begin{theorem}\label{main-theorem}
Let $(i_1,\dots ,i_d)$ be a reduced word for $w\in W$.
Let $\sim $ be the identifications given by any series of face 
collapses (cf. Definition ~\ref{collapse-map-defn}) on $\RR_{\ge 0}^d \cap  S_1^{d-1}$ 
such that (1) $x\sim y$ implies $f_{(i_1,\dots ,i_d)}(x) = f_{(i_1,\dots ,i_d)}(y)$,
and (2) the series of collapses eliminates all regions whose words are not reduced.
Then 
 $\overline{f_{(i_1,\dots ,i_d)}} : 
(\RR_{\ge 0}^d \cap S_1^{d-1})/\!\!\sim\hspace{.05in} \rightarrow Y_w$ is a homeomorphism which preserves
cell structure. 
\end{theorem}

\begin{proof}
The proof is by induction on $d$, with the case $d=1$ being trivial.  Therefore, 
we may assume the result for all finite Coxeter group elements of 
length strictly less than $d$.  Remark ~\ref{ident-remark} enables us to deduce
continuity of $\overline{f_{(i_1,\dots ,i_d)}}$ from continuity of $f_{(i_1,\dots ,i_d)}$.
Notice that $f_{(i_1,\dots ,i_d)}$ restricts to
any region obtained by setting some $t_i$'s 
to 0 since $x_i(0)$ is the identity matrix.  Whenever the 
resulting subword is reduced, results in \cite{Lu}
guarantee that $f_{(i_1,\dots ,i_d)}$ acts homeomorphically on this open cell.  
The requirements of Corollary ~\ref{regular-build}
regarding  closures of cells in the $(d-1)$-skeleton 
follow from our inductive hypothesis, along with 
the fact that any series of face collapses  will restrict to one
on the closure of any cell.
Thus, we may 
apply Corollary ~\ref{regular-build} to deduce that 
$\overline{f_{(i_1,\dots ,i_d)}}((\RR_{\ge 0}^d \cap S^{d-1}_1 )/\!\!\sim )$
is a finite CW complex with the restrictions of $\overline{f_{(i_1,\dots ,i_d)}}$ 
to the various cell closures in $(\RR_{\ge 0}^d \cap S_1^{d-1} )/\!\!\sim $ 
giving the characteristic maps, and that this CW complex structure
satisfies conditions 1,2 and 5 of Theorem ~\ref{sufficient-theorem}.
 Lemma   ~\ref{injective-condition}
confirmed condition 4 of Theorem ~\ref{sufficient-theorem}, while the 
result of \cite{BW} that Bruhat order is shellable and thin gives condition 3.
Thus, by Theorem ~\ref{sufficient-theorem}, 
$\overline{f_{(i_1,\dots ,i_d)}}((\RR_{\ge 0}^d \cap S_1^{d-1} )/\!\!\sim )$ is a 
regular CW complex with characteristic maps given by the restrictions of 
$\overline{f_{(i_1,\dots ,i_d)}}$ to the various
cell closures, which is exactly what is needed. 
\end{proof}

\begin{corollary}\label{lk-cor}
If $(i_1,\dots ,i_d)$ is a reduced word for $w$, then 
$f_{(i_1,\dots ,i_d)}$ induces a homeomorphism 
$\overline{f_{(i_1,\dots ,i_d)}}: (\RR_{\ge 0}^d \cap S_1^{d-1})/\!\!\sim \hspace{.05in} \rightarrow Y_w $
which preserves  cell structure.  Hence, $Y_w$ is a 
regular CW complex homeomorphic to a ball with Bruhat interval $(1,w]$ as its closure poset. 
\end{corollary}

\begin{proof}
By Theorem  ~\ref{homeom2-lemma}, 
$K = (\RR_{\ge 0}^d \cap S^{d-1}_1)/\!\!\sim $ is a regular CW complex
homeomorphic to a ball.
We chose $\sim $ so that  $f_{(i_1,\dots ,i_d)}(x)=f_{(i_1,\dots ,i_d)}(y)$
whenever $x \sim y $.  Combining 
with Lusztig's result that $f_{(i_1,\dots ,i_d)}$ is  continuous on
$\RR_{\ge 0}^d \cap S_1^{d-1}$ and a homeomorphism on  $\RR_{>0}^d \cap S_1^{d-1}$
(see \cite{Lu}, Section 4]), as well as the fact that 
our collapsing maps are identification maps, 
it follows that $\overline{f_{(i_1,\dots ,i_d)}}$ is also continuous on $K$. 
In particular, $K$ meets all the requirements of Theorem ~\ref{main-theorem}.
\end{proof}

\begin{corollary}
For $(t_1,\dots ,t_d), (t_1',\dots ,t_d') \in 
\RR_{\ge 0}^d$ and $(i_1,\dots ,i_d)$ any reduced word, 
$$x_{i_1}(t_1)\cdots x_{i_d}(t_d) = x_{i_1}(t_1')\cdots x_{i_d}(t_d')
\hspace{.05in} {\rm iff}
\hspace{.05in} (t_1,\dots ,t_d)\sim (t_1',\dots,t_d').$$
\end{corollary}

Finally we consider  $lk(u,w)$.  
Theorem ~\ref{homeom2-lemma} proved
regularity and determined homeomorphism type for 
the links of the cells in  $(\RR_{\ge 0}^d \cap S_1^{d-1})/\sim $, using the natural projection 
map in that context  (cf. Definition ~\ref{proj-map}).
This immediately transfers  to  yield analogous results for links in $Y_w$ via our  cell-preserving 
homeomorphism $\overline{f_{(i_1,\dots ,i_d)}}$.

\begin{corollary}
For each $u < w$ in Bruhat order, the subcomplex $lk(u,w)$ of $Y_w$ 
obtained as the image under the cell-preserving homeomorphism 
$\overline{f_{(i_1,\dots ,i_d)}}$ of the corresponding link within  $\RR_{\ge 0}^d /\!\sim $  
is a regular CW complex homeomorphic to a ball having the Bruhat interval 
$(u,w]$ as its poset of closure relations.
\end{corollary}

We expect  that 
this notion of link 
should  coincide, at least up to homeomorphism, with the notion of link given by  
Fomin and Shapiro for $Y_w$ in \cite{FS}  by virtue of 
Thom's first isotopy lemma (cf. \cite{GM}).   

 It would  be interesting  to understand better how $lk(u,w)$ relates both to subword
complexes (cf.\cite{Kn},  \cite{KM},  \cite{KM2}) 
and also to the synthetic  CW complexes for Bruhat intervals
 studied by Reading in \cite{Re}.

\section{Acknowledgments and sources of support}

The author was funded by NSF grants DMS-0500638, DMS-0757935, 
DMS-1002636 and DMS-1200730 
as well as the Ruth I. Michler Prize of the Association for 
Women in Mathematics,  the last of which provided a research semester at Cornell.
She  thanks the Cornell mathematics department for its warm 
hospitality and intellectually stimulating atmosphere during her stay there.

The author is  grateful to Sara Billey, Anders Bj\"orner,
Jim Davis, Sergey Fomin, Mark Goresky, Mike Hopkins,
Nets Hawk Katz,  Allen Knutson, Chuck Livingston, Mark McConnell, Sergey Melikhov,
Ezra Miller, Nathan Reading, David
Speyer and Lauren Williams  for very helpful discussions and comments on various
versions of the paper.   She also thanks 
the anonymous referees
for  insightful  questions and  comments 
which led to substantial improvements in the paper.

\end{document}